\documentclass[leqno]{amsart}

\usepackage[foot]{amsaddr}

\usepackage{color}
\usepackage{geometry}

\usepackage[numbers]{natbib}
\usepackage{amsmath}
\usepackage{mathtools}

\usepackage{subfiles}
\usepackage{graphicx}
\usepackage[section]{placeins}
\usepackage{hyperref}
\usepackage{stmaryrd}
\usepackage{amsmath,amssymb,amsfonts,amsthm,textcomp}
\usepackage{mathtools}
\usepackage{tensor}
\usepackage{epstopdf}
\usepackage{multicol}
\usepackage{empheq}
\usepackage{tikz,pgfplots}
\usetikzlibrary{calc}
\usepackage{enumitem}
 \newcommand{\lJump}{[\![}
\newcommand{\rJump}{]\!]}

\newtheorem{remark}{Remark}
\usepackage[linesnumbered,ruled,vlined]{algorithm2e}
\usepackage{subfigure}
\usepackage{fancyhdr}
\DeclareMathOperator{\spn}{span}

\graphicspath{{./Figures/}
              {./Figures/Convergence/} 
 			{./Figures/Impact/} 
			{./Figures/TwistedBar/}
			{./Figures/TossRule/}}

\newfont{\tenbfsl}{cmbxti9 scaled 1200}
\newfont{\tenbbb}{msbm10}
\newfont{\svnbbb}{msbm8}

\newtheorem{thm}{Theorem}

\theoremstyle{remark}

\theoremstyle{definition}

\newcounter{syn}[section] 
\renewcommand{\thesyn}{\arabic{section}.\arabic{syn}}

\makeatletter
\def\threevdots{\mskip+4mu\vbox{\baselineskip2.25\p@ \lineskiplimit\z@
  \kern4.9\p@\hbox{.}\hbox{.}\hbox{.}}\mskip+3.8mu}
\makeatother

\begin{document}

\title[Adaptive explicit predictor/multicorrector time marching for elastodynamics]{An explicit predictor/multicorrector time marching with automatic adaptivity for finite-strain elastodynamics}
\author{Nicol\'as A. Labanda$^{\diamondsuit}$}
\address{$^{\diamondsuit}$ Curtin Institute for Computation \& School of Electrical
	Engineering, Computing and Mathematical Sciences, Curtin
	University, P.O. Box U1987, Perth, WA 6845, Australia}
\email{nlabanda@facet.unt.edu.ar (N.A. Labanda)}
\author{Pouria Behnoudfar$^{\flat}$}
\address{$^{\flat}$ Mineral Resources, Commonwealth Scientific and Industrial Research Organisation (CSIRO), Kensington, Perth, WA 6152, Australia}
\author{Victor M. Calo$^{\natural}$}
\address{$^{\natural}$ School of Electrical
	Engineering, Computing and Mathematical Sciences, Curtin
	University, P.O. Box U1987, Perth, WA 6845, Australia}

\date{\today}

 \begin{abstract}
 \noindent
 We propose a time-adaptive predictor/multi-corrector method to solve hyperbolic partial differential equations, based on the generalized-$\alpha$ scheme that provides user-control on the numerical dissipation and second-order accuracy in time. Our time adaptivity uses an error estimation that exploits the recursive structure of the variable updates. The predictor/multicorrector method explicitly updates the equation system but computes the residual of the system implicitly. We analyze the method's stability and describe how to determine the parameters that ensure high-frequency dissipation and accurate low-frequency approximation. Subsequently, we solve a linear wave equation, followed by non-linear finite strain deformation problems with different boundary conditions. Thus, our method is a straightforward, stable and computationally efficient approach to simulate real-world engineering problems. Finally, to show the performance of our method, we provide several numerical examples in two and three dimensions.  These challenging tests demonstrate that our predictor/multicorrector scheme dynamically adapts to sudden energy releases in the system, capturing impacts and boundary shocks. The method efficiently and stably solves dynamic equations with consistent and under-integrated mass matrices conserving the linear and angular momenta as well as the system's energy for long-integration times.

\textbf{Keywords:} Time - adaptivity, Stability, Finite-strain, Dissipation control, Hyperbolic equations.

 \end{abstract}

\maketitle




 \section{Introduction}
\label{section:Introduction}

Time-adaptive schemes are crucial for the efficient solution of hyperbolic partial differential equations, specifically in solid mechanics.  These highly nonlinear problems with changing boundary conditions require computationally robust adaptive numerical methods that control the time-step size to keep their cost tractable. Among several time-marching method classifications, the most conventional is the division into the explicit and implicit schemes. Some implicit schemes are unconditionally stable, allowing larger time steps, but are computationally more expensive and usually unsuitable for highly nonlinear mechanical problems~\cite{ hughes2012finite}.  Nonetheless, most time integrators for structural mechanics and deformable bodies belong to the implicit Newmark family~\cite{ Newmark1959}. Even though these implicit time-marching schemes are linearly stable, they still require numerical damping to control the detrimental impact of the poor resolution of the high-frequency modes on the problem dynamics. The engineering community continually develops implicit schemes with controlled numerical dissipation; for example,  classical methods with second-order accuracy exist, such as the Wilson-$\theta$ method~\cite{ wilson1968}, the HHT-$\alpha$ method~\cite{ Hilber1977}, the WBZ-$\alpha$~\cite{ wood1980}, the $\rho$~method~\cite{ Bazzi1982} and the $\theta_{1}$~method~\cite{ HOFF1988367, HOFF198887, HOFF198987}. In the early 90s, Chung and Hulbert proposed the generalized-$\alpha$ method, which generalizes these methods and comprises them as particular cases, to solve the hyperbolic partial differential equations that model structural dynamics problems~\cite{ chung1993time}. The generalized-$\alpha$ method controls the high-frequency dissipation while preserving accuracy in the low-frequency regions. This time-marching method successfully tackled a vast number of applications such as fluid-structure interactions~\cite{Bazilevs2006,Bazilevs2008} and finite-strain deformation~\cite{ NOELS2005358, ROSSI2016208, ScovazziIJMNE2016, BONET2015689, GIL2016146, LAVRENCIC2020107075}. Recent extensions of the method include splitting techniques~\cite{ BEHNOUDFAR2020100021} and higher-order generalized-$\alpha$ methods~\cite{ behnoudfar2019higher,BEHNOUDFAR2021113725,LOS2020109}.

Alternatively, Miranda et al.~\cite{ Miranda1989} proposed an explicit HHT-$\alpha$ method that extends the implicit HHT-$\alpha$ method, improving prior work based on the Newmark method~\cite{ Hughes1978, Hughes1978b}.  Two decades later, Hulbert and Chung~\cite{ HULBERT1996175} proposed an explicit predictor-corrector version of the generalized-$\alpha$ method, capable of controlling the numerical dissipation and maintaining the second-order accuracy, while~\cite{ Daniel2003} rewrote the classical generalized-$\alpha$ method to obtain an explicit/implicit version. Alternative explicit predictor-corrector methods in the literature are~\cite{ BONELLI2002695, Bonelli2005, Tripodi2016, lopez2020}.

Most authors use a uniform time step size to fulfil the initial CFL condition, assuming the solution behavior remains unchanged throughout the simulation, a critical weakness when dealing with nonlinear problems. The first attempts to address this limitation focused on quasi-static structural problems using continuation methods~\cite{ Thomas1973, Schmidt1978, PADOVAN1982365}. The critical time step, the maximum time-step size that ensures convergence, is related to the maximum eigenvalue of the amplification matrix. In large-scale simulations, the eigenvalue computation is extremely onerous; additionally, the stiffness matrix may change between time steps and iterations as the nonlinear process evolves. Thus, cheap convergence and error estimates are extremely valuable. Furthermore, critical time-step estimates based on eigenvalue inequalities are reliable only when the problem has a smooth response and fail when subject to sudden energy releases or impulsive loads~\cite{ Belytschko2014}. 

Other authors estimate the temporal error by comparing the discrete solution against higher-order methods; most of these methods modify Runge-Kutta integration schemes, using Butcher tables to reduce the computational cost of higher-order calculations~\cite{ Gustafsson1991, SODERLIND2006225, Butcher2016}. In structural and solid mechanics, the central difference integration with variable-step size selection was a seminal contribution in time adaptivity~\cite{ PARK1980241, UNDERWOOD1980259}, where the apparent highest frequency predicts the time-step size. Later, several authors used truncation error estimates for Newmark integrators for simple dynamics problems~\cite{ Zienkiewicz1991, Zeng1992, Wiberg1993} and extended in~\cite{ Romero2006}. Hulbert then introduced the concept of local error estimation~\cite{ HULBERT1995155} for both the explicit and implicit generalized-$\alpha$ schemes for hyperbolic partial differential equations with a focus on structural mechanics.

This paper presents an automatic time-step adaptive algorithm that uses a built-in error estimate that arises from the recursive update that the time-marching scheme uses. Our time-marching approach builds on the generalized-$\alpha$ method, allowing dissipation control in the higher frequencies while accurately approximating the low frequencies in the discrete spectrum. This method computes every update explicitly, while the system's residual corresponds to the implicit residual, ensuring the explicit dynamics track the proper nonlinear behavior of the underlying physical system. By comparing the norm of each iteration component, we evaluate the local convergence of the problem and adapt the time-step size to reach a bounded result. We apply our method to finite-strain problems to show its ability to deal with changing scenarios and non-smooth solutions. Furthermore, our scheme delivers a stable, robust, and computationally efficient formulation to solve highly nonlinear problems, preserving energy for long-integration periods even when the updating procedure uses lumping techniques.  

The paper organization follows: Section~\ref{section:insight} explains our method, while Section~\ref{section:stability} defines the method's parameters and analyses its stability. We solve the dynamics for a hyperelastic material in Section~\ref{section:nonlinear} and  detail the method's computational implementation in Section~\ref{section:Implement}. Section~\ref{section:Examples} presents theoretical and practical numerical examples; finally, we draw conclusions in Section~\ref{section:Conclusions}.

 \section{Predictor/multicorrector time integration}
\label{section:insight}

In this section, we introduce a predictor/multicorrector scheme based on the generalized-$\alpha$ philosophy. For this, we discuss the strong and weak forms of the time-dependent problem. Then, we introduce a time discretization to update the displacement and acceleration.

\subsection{Strong and weak formulations of a time-dependent problem}

Let $\mathcal{B}\subset \mathbb{R}^{d}$, where $d \in \left\lbrace 1,2,3\right\rbrace$ is the spatial dimension and $I \in \left[0,T \right]$ is the time interval. We consider the following wave equation:
\begin{equation}
  \left\{
    \begin{aligned}
      \ddot{{u}} - \Delta {u} &= {f}\left({x},t\right), &&
                                                               \text{in} \ \ \mathcal{B}\times I, \\ 
      {u}\left({x},t\right)  & = {u}_D,  && \text{on} \ \  \partial\mathcal{B} \times I,  \\
      {u} \left({x}, 0 \right) & = {u}_0,\quad\quad 
      \dot{{u}}\left({x}, 0 \right)  = \dot{{u}}_0,  && \text{on} \ \  \mathcal{B},
    \end{aligned}\right.
  \label{eq:HyperbolicPDE}
\end{equation}
where $x \in \mathcal{B}$ is the body's coordinate system, ${u} : \mathcal{B} \times I \rightarrow \mathbb{R}^{m}$, with $m \geq 1$, is the unknown field, $\ddot{{u}} \coloneqq \frac{\partial^{2} {u}}{\partial t^{2}} $ is its second time derivative and $\Delta u=\frac{\partial^2 u}{\partial x_i^2}$ its spatial Laplacian, ${f} : \mathcal{B} \times I \rightarrow \mathbb{R}^{m}$ is the source term, and $\partial\mathcal{B}$ denotes the spatial boundary of $\mathcal{B}$. The partial differential equation (PDE)~\eqref{eq:HyperbolicPDE}$_1$ is subject to ${u}_D$, ${u}_0$, and $\dot{{u}}_0$, the Dirichlet boundary~\eqref{eq:HyperbolicPDE}$_2$ and initial conditions ~\eqref{eq:HyperbolicPDE}$_3$, respectively.

Multiplying~\eqref{eq:HyperbolicPDE}$_1$ by a smooth weighting function ${v}$ in $\mathcal{V}$ that vanishes on $\partial \mathcal{B}$, and integrating by parts in space, we obtain:
\begin{align}
    \displaystyle \int_{I} \int_{\mathcal{B}} {v}  {f}  \
  d\mathcal{B}\ dt
  &=\displaystyle \int_{I} \int_{\mathcal{B}} \left( {v} \ddot{{u}}  -
    {v} \Delta {u}  \right) d\mathcal{B}\ dt,  
    = \int_{I} \int_{\mathcal{B}} \left( {v}  \ddot{{u}} + \nabla {v}
    \cdot \nabla {u} \right) d\mathcal{B}\ dt. \label{eq:HyperbolicPDEWeak2}
\end{align}
Using standard notation, we state that the above formulation is sufficient to ensure that all terms in~\eqref{eq:HyperbolicPDEWeak2} belong to  ${H}^1$ in space and ${L}^2$ in time. In this sense, the test space is
\begin{equation}
  \mathcal{V} \coloneqq {H}^{1}_{0}\left( \mathcal{B}\right) \coloneqq
  \left\lbrace  {v} \in {L}^2\left( \mathcal{B}\right)  |   \nabla {v}
    \in {L}^2\left( \mathcal{B}\right), {v} = 0 \in
    \partial\mathcal{B} \right\rbrace  
  \label{eq:HyperbolicPDEWeak3}
\end{equation}
Given~\eqref{eq:HyperbolicPDEWeak2}, then $\ddot{{u}}$ and ${f}$ belong to the dual space of the test space $\mathcal{V}^{\ast} = {H}^{-1}\left( \mathcal{B}\right)$. Thus, the trial space is
\begin{equation}
  \mathcal{U} \coloneqq \left\lbrace {u} \in \mathcal{V} \ \ | \ \
    \ddot{{u}} \in \mathcal{V}^{\ast},  {u} \left({x}, 0 \right) = {u}_0,
    \dot{{u}}\left({x}, 0 \right)  = \dot{{u}}_0,
      {u} = {u}_D \in  \partial\mathcal{B} 
  \right\rbrace . 
  \label{eq:HyperbolicPDEWeak4}
\end{equation}
where $  {u}_0, \dot{{u}}_0 \in {L}^2\left( \mathcal{B}\right).$

Finally, the weak problem reads: \textit{Given $ {f} \in \mathcal{V}^{\ast}$ a source term, find ${u} \in \mathcal{U}$ such that}
\begin{align}
  \displaystyle \int_{I} \left\langle  {v} , \ddot{{u}}  \right\rangle
  dt
  &= \displaystyle \int_{I} \left\langle {v} , {f}   \right\rangle dt
    - \int_{I} a\left(  {v} , {u} \right) dt,
  &     \forall v \in  \mathcal{V}
    \label{eq:HyperbolicPDEproblem}
\end{align}
where $\left\langle \cdot , \cdot \right\rangle $ is the duality pairing between $\mathcal{V}$ and $\mathcal{V}^{\ast}$, $a\left( {v} , {u} \right) = \left( \nabla {v} , \nabla {u} \right)$ is a linear form with $\left( \cdot, \cdot \right)$ the inner product in ${L}^{2}\left( \mathcal{B}\right)$.

\subsection{Predictor/multicorrector time-marching method}

Herein, we describe our time-marching method. For this aim, we discretize the time domain $t_0<t_1<...<t_n<...<t_f$, defining a time-step $\Delta t_n=t_n-t_{n-1}$. Approximating ${u}\left(t_n\right)$ and $\ddot{{u}} \left( t_n \right)$, respectively, using ${u}_n$ and $\ddot{{u}}_n$, we can state (for more details, see,~\cite{chung1993time}): 
\begin{equation}
  \begin{aligned}
    t_{n+\alpha_f} &= t_{n} + \alpha_f \Delta t, \\
    {u}_{n+\alpha_f} &=  \displaystyle {u}_{n} + \alpha_{f}  \lJump {u} \rJump ,  \\
    \ddot{{u}}_{n+\alpha_m} &=  \displaystyle \ddot{{u}}_{n}
    + \alpha_{m}  \lJump \ddot{{u}} \rJump, 
  \end{aligned}
  \label{eq:EqGenAlpha1}
\end{equation}
where $\lJump  \bullet \rJump = \left(\bullet \right)_{n+1} - \left( \bullet \right)_{n}$ is the time increment of a certain field $\bullet$. Using a Taylor series in time, we have:
\begin{align}
  {u}_{n+1}  &=  \displaystyle \sum_{j = 0}^{\phi+1} \frac{{\Delta
               t}^{j}}{j !} \left.{\frac{\partial^j {u}}{\partial
               t^j}}\right|_{n} ,  
               \label{eq:EqGenAlpha2}
\end{align}
where $\phi$ is the expansion order. Substituiting~\eqref{eq:EqGenAlpha1} into~\eqref{eq:HyperbolicPDEproblem}, the semi-discrete problem reads: 
\begin{equation}
\left\langle {v}  ,  \ddot{{u}}_{n+\alpha_m} \right\rangle   =  \displaystyle \left\langle {v}  , {f}_{n+\alpha_f} \right\rangle  -  a \left( {v}  ,  {u}_{n+\alpha_f}  \right)    , \ \ \ \  \forall {v} \in  \mathcal{V}  \ .
\label{eq:HyperbolicPDEproblemGalpha}
\end{equation}
\begin{remark}
  We denote the last-term constant in~\eqref{eq:EqGenAlpha2} as $\gamma_{\phi}$, while its increment is
  \begin{align}
    \lJump {u} \rJump &= \displaystyle \sum_{j = 1}^{\phi}
                        \frac{{\Delta t}^{j}}{j !}
                        \left.{\frac{\partial^j {u}}{\partial
                        t^j}}\right|_{n} +  \ \gamma_{\phi} {\Delta
                        t}^{\phi} \left[\!\!\left[
                        {\frac{\partial^{\phi} {u}}{\partial
                        t^{\phi}}} \right]\!\!\right]  
                        \label{eq:EqGenAlpha3}
  \end{align}
  For hyperbolic PDEs, we let $\phi = 2$,  
  \begin{align}
    \lJump {u} \rJump  &=  \Delta t \left( \dot{{u}}_{n} + \Delta t
                         \left( \frac{\ddot{{u}}_{n} }{2}  +
                         \gamma_{2} \lJump \ddot{{u}}  \rJump \right)
                         \right),
                         \label{eq:taylorhyper}
  \end{align}
  while we compute the velocity as
  \begin{align}
    \lJump \dot{{u}} \rJump  &=  \Delta t \left( \ddot{{u}}_{n} +
                               \gamma_{1} \lJump \ddot{{u}}  \rJump
                               \right). 
                               \label{eq:taylorhyper11}
  \end{align}
  In the literature, the coefficients $\gamma_1$ and $\gamma_2$ are usually denoted in the bibliography as $\gamma$ and $\beta$, respectively~\cite{ chung1993time}. We use this notation in the remainder of the paper. 
\end{remark}

\subsection{Update in terms of ${u}$}

We substitute~\eqref{eq:EqGenAlpha1} into~\eqref{eq:HyperbolicPDEproblemGalpha} and solving in terms of $u$, we have:
\begin{equation}
\displaystyle \left\langle {v}  , \ddot{{u}}_{n} + \alpha_{m}  \lJump \ddot{{u}} \rJump \right\rangle   =  \left\langle {v}  , {f}_{n+\alpha_f}  \right\rangle  -  a \left(   {v} ,  {u}_{n} + \alpha_{f}  \lJump {u} \rJump \right).  
\label{eq:Galphau1}
\end{equation}
We deduce a set of useful equations for different problems, let $a \left( {v} , {u}_{n} + \alpha_{f} \lJump {u} \rJump \right)$ be a nonlinear functional. In this sense, we can approximate it using a second-order Taylor series expansion as
\begin{equation}
a \left( {v} , {u}_{n} + \alpha_{f}  \lJump {u} \rJump  \right) \approx a \left(   {v} , {u}_{n} \right) + \alpha_{f} \ a' \left(   {v}  , \lJump {u} \rJump ;  {u}_{n} \right),
\end{equation}
where $a' \left(  {v} , \lJump \bullet \rJump  ;  \bullet_{n} \right) $ is the Gateaux derivative of $a \left(   {v}  , \bullet  \right)$, defined as
\begin{equation}
a' \left(  {v} , \lJump {u} \rJump  ;  {u}_{n} \right) = \frac{d}{d\epsilon} a \left( {v} , {u}_{n} + \epsilon  \lJump {u} \rJump  \right) \vert_{\epsilon=0},
\end{equation}
Then, combining this equation with~\eqref{eq:Galphau1} leads to:
\begin{equation}
  \displaystyle  \alpha_{m}  \left\langle {v} ,  \lJump  \ddot{{u}}
    \rJump  \right\rangle = \displaystyle - \alpha_{f} \  a' \left(
    {v}  ,  \lJump {u} \rJump ; {u}_{n}  \right) + \left\langle  {v}  ,
    {f}_{n+\alpha_f}  \right\rangle  - \left\langle  {v} ,
    \ddot{{u}}_{n} \right\rangle -  a \left(  {v} , {u}_{n}   \right), 
\end{equation}
which we can express in terms of $\lJump {u} \rJump$ using~\eqref{eq:EqGenAlpha2} as
\begin{equation}
  \displaystyle \left\langle  {v}  , \lJump  {u} \rJump \right\rangle
  =  \displaystyle - \frac{\alpha_f \Delta t^2 \beta}{\alpha_m} a'
  \left(   {v} ,   \lJump {u} \rJump;  {u}_{n}  \right) + \ell \left(
    {v} \right),   
  \label{eq:Galphau1bis}
\end{equation}
with
\begin{equation}
  \ell \left( {v} \right) = \displaystyle \frac{\beta \Delta
    t^{2}}{\alpha_m}    \displaystyle \left[ \left\langle {v} ,
      {f}_{n+\alpha_f}  \right\rangle - \left\langle  {v} ,
      \ddot{{u}}_{n}  \right\rangle -  a \left( {v} ,  {u}_{n}  ,
    \right) \right]  +  \left\langle {v} , \sum_{j = 1}^{2}
    \frac{{\Delta t}^{j}}{j !} \left.{\frac{\partial^j {u}}{\partial
          t^j}}\right|_{n}  \right\rangle. 
\label{eq:Galphau1bisbis}
\end{equation}
We write the time increment as $\lJump {u} \rJump^{k+1} = \lJump {u} \rJump^{k} + \delta {u}^{k+1}$ in~\eqref{eq:Galphau1bis}, and considering the left term an iteration ahead respect to the right-hand side, we obtain
\begin{align}
    \displaystyle \left\langle {v}  ,  \lJump  {u} \rJump^{k+1}
    \right\rangle    &=  \displaystyle - \frac{\alpha_f \Delta t^{2}
                       \beta}{\alpha_m} a' \left(
                       {v} ,  \lJump {u} \rJump^{k} ; {u}_{n}  \right)
                       + \ell \left( {v} \right), 
\end{align}
which we can rewrite as
\begin{align}
    \displaystyle \left\langle  {v} , \delta {u}^{k+1}  \right\rangle
                     &= \displaystyle - \frac{\alpha_f \Delta t^{2}
                       \beta}{\alpha_m} {a}' \left(  {v}
                       , \lJump {u} \rJump^{k} ; {u}_n   \right) -
                       \left\langle {v}  ,  \lJump  {u} \rJump^{k}
                       \right\rangle  + \ell \left( {v} \right). 
\label{eq:Galphau2}
\end{align}
Now, we consider a spatial discretization for $\mathcal{V}$ and $\mathcal{U}$:
\begin{equation}
  {u}^{h} = \sum_{i=1}^{n} \Phi_i \mathbf{U}_{i}  \subset  \mathcal{U}
  , \ \ \text{and} \ \ {v}^{h} \in \spn \left\lbrace  \Psi_j
  \right\rbrace^{m}_{j=1} \subset  \mathcal{V}, 
\end{equation}
with $\Phi_i$ and $\Psi_j$ being appropriate trial and test functions, respectively, while $\mathbf{U}_{i}$ is the unknown vector. Therefore, we obtain the fully-discrete problem as:
\begin{equation}
  \begin{aligned}
    \displaystyle \mathbb{M} \cdot \delta  \mathbf{U}^{k+1}
    &= \displaystyle -  \frac{\alpha_{f} \Delta t^{2} \beta
    }{\alpha_{m} }  \mathbb{K} \cdot  \lJump   \mathbf{U} \rJump^{k} -
    \mathbb{M} \cdot \lJump  \mathbf{U} \rJump^{k}    +  \mathbf{L},
    \\ 
    \displaystyle \delta  \mathbf{U}^{k+1}
    &= \displaystyle \left( -  \frac{\alpha_{f}  \Delta t^{2}
        \beta}{\alpha_{m} } \mathbb{M}^{-1} \mathbb{K} - \mathbb{I}
    \right) \cdot \lJump  \mathbf{U} \rJump^{k}    +  \mathbb{M}^{-1}
    \cdot \mathbf{L}, \\ 
    \displaystyle \delta  \mathbf{U}^{k+1}
    &= \displaystyle \left( \mathbb{U} - \mathbb{I} \right) \cdot
    \lJump  \mathbf{U} \rJump^{k}    +  \mathbf{H}, 
  \end{aligned} 
  \label{eq:Galphau4}
\end{equation}
where $\mathbb{M}$ and $\mathbb{K}$ are the mass stiffness matrices defined as
\begin{align*}
  \mathbb{M}_{ij} &= \left\langle \Psi_i  , \Phi_j \right\rangle =
                    \int_{\mathcal{B}} \Psi_i \Phi_j  \ d \mathcal{B},
  &
    \mathbb{K}_{ij} &= a' \left( \Psi_i , \Phi_j ; {u}_n \right) \\
  \mathbf{H} &=  \mathbb{M}^{-1} \cdot \mathbf{L},
  & \mathbf{L} &= \ell \left( \Psi_j \right) \ ,
\end{align*}
$\mathbb{M}^{-1}$ is the inverse matrix of the mass matrix $\mathbb{M}$ and $\mathbb{U} = - \frac{\alpha_{f} \Delta t^{2} \beta}{\alpha_{m} } \mathbb{M}^{-1} \mathbb{K}$ is the updating matrix, with the identity matrix $\mathbb{I}$.  Starting from~\eqref{eq:Galphau4} and assuming that, for $k=0$ the variable increment is $\lJump \mathbf{U} \rJump^{0} = 0$, the first increment is
\begin{align}
    \delta \mathbf{U}^{1} =  \mathbf{H}  \qquad \Rightarrow \qquad
  \lJump \mathbf{U} \rJump^{1} = \mathbf{H} . 
  \label{eq:Galphau5}
\end{align}
Substituting $\lJump \mathbf{U} \rJump^{k}$ into~\eqref{eq:Galphau4}, one can readily compute the increment $\delta \mathbf{U}^{k+1} $ in a recursive manner as:
\begin{equation}
  \begin{aligned}
    \delta \mathbf{U}^{2} & = \mathbb{U} \cdot \mathbf{H}
    && \Rightarrow && \lJump \mathbf{U} \rJump^{2}
     = \left(\mathbb{I} + \mathbb{U}\right) \cdot \mathbf{H},\\
    \delta \mathbf{U}^{3} &= \mathbb{U}^{2} \cdot \mathbf{H}
    && \Rightarrow && \lJump \mathbf{U} \rJump^{3}
    = \left(\mathbb{I} + \mathbb{U} + \mathbb{U}^{2} \right) \cdot \mathbf{H}, \\
    &&& \ \, \vdots &&\\
    \delta \mathbf{U}^{N} &= \mathbb{U}^{N-1} \cdot \mathbf{H}
    && \Rightarrow && \lJump \mathbf{U} \rJump^{N}
    = \left(\mathbb{I} + \mathbb{U} + \mathbb{U}^{2} + ...
      + \mathbb{U}^{N-1} \right) \cdot \mathbf{H},
  \end{aligned} 
  \label{eq:recursiveu}
\end{equation}
Finally, the predictor/multicorrector updates the variable using an explicit generalized-$\alpha$ scheme as follows
\begin{equation}
  \boxed{
    \ \lJump \mathbf{U} \rJump^{N} = \left[ \sum^{N-1}_{j=0}
      \mathbb{U}^{j} \right] \cdot \mathbf{H} .\ 
  }
  \label{eq:Galphau7}
\end{equation}

\subsection{Update in terms of $\ddot{{u}}$}

Following a similar approach, we express the explicit update in terms of the field derivative $\ddot{{u}}$, starting from~\eqref{eq:Galphau1} and~\eqref{eq:EqGenAlpha2}. Thus, the update is 
\begin{equation}
  \boxed{
    \ \lJump \ddot{\mathbf{U}} \rJump^{N} = \left[ \sum^{N-1}_{j=0}
      \mathbb{U}^{j} \right] \cdot \widehat{\mathbf{H}} \ , 
  }
  \label{eq:Galphadotu1}
\end{equation}
where we compute the vector $\widehat{\mathbf{H}} =  \mathbb{M}^{-1} \cdot \widehat{\mathbf{L}}$ using the following definition of the right-hand side
\begin{align}
    \ell \left( {v} \right) = \displaystyle \frac{1}{\alpha_m}  \left[
  \left\langle  {v}  ,  {f}_{n+\alpha_f} \right\rangle - \left\langle
  {v}  ,  \ddot{{u}}_{n}  \right\rangle -  a \left(   {v}  ,  {u}_{n}
  \right) - \alpha_f \ {a}' \left(  {v}  ,  \sum_{j = 1}^{2}
  \frac{{\Delta t}^{j}}{j !} \left.{\frac{\partial^j {u}}{\partial
  t^j}}\right|_{n}  ; {u}_{n}  \right) \right].
  \label{eq:redefinitionofl}
\end{align}
These definitions complete the predictor/multicorrector scheme in terms of the time derivative $\ddot{{u}}$.
\begin{remark}
  When the update uses only the first term of the update scheme~\eqref{eq:Galphau7} or~\eqref{eq:Galphadotu1}, the update formula is similar to the explicit generalized-$\alpha$ proposed in~\cite{ HULBERT1996175}; nevertheless, our coefficients definitions are different from those used in~\cite{ HULBERT1996175}, as we discuss in the following section. Herein, we generalize the method, including the truncation error that employs higher-order terms, which allow us to predict the series convergence. 
\end{remark}

\begin{remark}
Proposed scheme does not require to update the stiffness matrix within the time increment.
\end{remark}

 \section{Stability and accuracy for hyperbolic problems}
\label{section:stability}

In this section, we analyze the stability analysis of our scheme following~\cite{ hughes2012finite, HULBERT1996175}. We use the spectral decomposition to simplify the set of equations to analyze
\begin{equation}
  \lJump  {\mathbf{U}}_{i} \rJump^{N} = \left[ \sum^{N-1}_{j=0} \left(
      -  \frac{\alpha_{f} \Delta t^{2} \beta}{\alpha_{m} }
      \omega^{2}_{i} \right)^{j} \right] \cdot
  \widetilde{\mathbf{H}}_{i}, 
\end{equation}
where $\omega_{i}$ are the eigenvalues of the system
\begin{equation}
  \left(\mathbb{M}^{-1} \mathbb{K} - \omega^{2}_{i}  \mathbb{I} \right)
  \Psi_{i} = \mathbf{0}, 
  \label{eq:freq}
\end{equation}
that represent the frequency of mode $i$ for the base vector $\Psi_{i}$. Furthermore, we use a compact notation for the sum term:
\begin{equation}
  \left[{P^{N}_{k}} \right]_{i} =  \sum^{N-1}_{j=k} \left( -
    \frac{\alpha_{f} \Delta t^{2} \beta}{\alpha_{m} } \omega^{2}_{i}
  \right)^{j} .
  \label{eq:modaldecomposition}
\end{equation}
In what follows, we spectrally decompose the matrix in~\eqref{eq:modaldecomposition}. This decomposition allows us to obtain a set of parameters that control the method's dissipation and stability region.

\subsection{Accuracy analysis}

Herein, before analyzing the method's stability behavior, we study its order of accuracy.  We seek to deliver second-order accuracy in time; therefore, we define the parameter $\gamma$ as follows
\begin{thm} \label{thm:3o}
  We obtain a scheme with second-order accuracy in time by setting
  \begin{equation} \label{eq:3ov1}
    \gamma=\frac{1}{2}-\alpha_{f}+\alpha_m.
  \end{equation}
\end{thm}
\begin{proof}
  We consider~\eqref{eq:EqGenAlpha2} and the displacement update formula~\eqref{eq:Galphau7}, and perform algebraic substitutions; then, the equation system at each time step $n+1$ in terms of a known right-hand side becomes
  \begin{equation}
    \mathbb{A} \mathbf{X}_{n+1} = \mathbb{B} \mathbf{X}_{n} + \mathbf{F}^{\sharp}_{n+\alpha_f},
  \end{equation}
  where the matrices and vectors are
  \begin{equation}
    \begin{aligned}\mathbb{A}&=
      \begin{bmatrix}
	1 & 0 & -\beta \\
	0 & 1 & -\gamma \\
	1& 0 & 0 \\
      \end{bmatrix},
      &
      \mathbb{B}&=
      \begin{bmatrix}
	1 & 1 & \displaystyle \frac{1}{2}-\beta \\
	0 & 1 & \displaystyle 1-\gamma \\
	\displaystyle   1 + \frac{1}{\alpha_f}   \left[{P^{N}_{1}}
        \right]_{i} & \displaystyle \left[{P^{N-1}_{0}} \right]_{i} &
        \displaystyle \left( \frac{1}{2} - \frac{\beta }{\alpha_m }
        \right) \left[{P^{N-1}_{0}} \right]_{i} 
      \end{bmatrix},\\
      \mathbf{X}_{n}&=
      \begin{bmatrix}
	\mathbf{U}_{n}  \\
	\Delta t \ \dot{\mathbf{U}}_{n}  \\
	\Delta t^2 \ddot{\mathbf{U}}_{n}\\ 
      \end{bmatrix}, &
      \mathbf{F}^{\sharp}_{n+\alpha_f} &=
      \begin{bmatrix}		
	0 \\
	0 \\
	\displaystyle \frac{\beta \Delta t^2}{\alpha_m}
        \left[{P^{N-1}_{0}} \right]_{i} \left[\mathbb{M}^{-1}
          \mathbf{F}_{n+\alpha_f}  \right]_{i}  
      \end{bmatrix}.
    \end{aligned}
  \end{equation}
  We update the matrices using the acceleration update formula~\eqref{eq:Galphadotu1}, 
  \begin{equation}
    \begin{aligned}\mathbb{A}&=
      \begin{bmatrix}
	1 & 0 & -\beta \\
	0 & 1 & -\gamma \\
	0& 0 & 1 \\
      \end{bmatrix},\qquad\qquad
      \mathbb{B}=
      \begin{bmatrix}
	1 & 1 & \displaystyle \frac{1}{2}-\beta \\
	0 & 1 & \displaystyle 1-\gamma \\
	\displaystyle   \frac{1}{\alpha_f  \beta } \left[{P^{N}_{1}}
        \right]_{i} & \displaystyle   \frac{1}{ \beta}
        \left[{P^{N}_{1}} \right]_{i} & \displaystyle   1 -
        \frac{1}{\alpha_m} \left[{P^{N-1}_{0}} \right]_{i} +
        \frac{1}{2 \beta} \left[{P^{N}_{1}} \right]_{i}  
      \end{bmatrix}\\
           \mathbf{F}^{\sharp}_{n+\alpha_f} &=
      \begin{bmatrix}		
	0 \\
	0 \\
	\displaystyle \frac{1}{\alpha_m} \left[{P^{N-1}_{0}}
        \right]_{i} \left[\mathbb{M}^{-1}  \mathbf{F}_{n+\alpha_f}
        \right]_{i}  
      \end{bmatrix}.
    \end{aligned}
  \end{equation}
  The amplification matrix  $\mathbb{G}$ for a hyperbolic problem is
  \begin{equation} \label{eq:ampm}
    \mathbb{G}=\mathbb{A}^{-1}\mathbb{B}.
  \end{equation}
  For an arbitrary $3\times 3$ invertible matrix, we can state~\cite{ hughes2012finite}:
  \begin{equation} \label{eq:a40}
    \mathbb{G}_0 \bold{U}_{n+1} - \mathbb{G}_1 \bold{U}_n +
    \mathbb{G}_2 \bold{U}_{n-1} - \mathbb{G}_3 \bold{U}_{n-2} = 0, 
  \end{equation}
  where the coefficients $\mathbb{G}_i$ for $i=1,2,3$ are the amplification matrix invariants; that is, $\mathbb{G}_0 = 1$, $\mathbb{G}_1 $ is the trace of $\mathbb{G}$, $\mathbb{G}_2$ is the sum of principal minors of $\mathbb{G}$, and $\mathbb{G}_3$ is the determinant of $\mathbb{G}$. Using a Taylor expansion analogous to the one in~\eqref{eq:EqGenAlpha2}, we obtain
  \begin{equation} \label{eq:Taypar}
    \begin{aligned}
      \bold{U}_{n+1} & = \bold{U}_n + \Delta t \dot{\bold{U}}_n +
      \dfrac{\Delta t^2}{2} \ddot{\bold{U}}_n+\mathcal{O}(\tau^3), \\ 
      \bold{U}_{n-1} & = \bold{U}_n - \Delta t
      \dot{\bold{U}}_n+\dfrac{\Delta t^2}{2}
      \ddot{\bold{U}}_n+\mathcal{O}(\tau^3), \\ 
      \bold{U}_{n-2} & = \bold{U}_n - 2 \Delta t \dot{\bold{U}}_n +
      4\frac{\Delta t^2}{2} \ddot{\bold{U}}_n +\mathcal{O}(\tau^3). 
    \end{aligned}
  \end{equation}
  Then, substituting~\eqref{eq:Taypar} into~\eqref{eq:a40} and setting the sum limit $N=2$ in~\eqref{eq:modaldecomposition}, we obtain
  \begin{multline}
    \frac{\left[{P^{1}_{1}} \right]^2 \left(4 \beta
        (\alpha_f-\alpha_m)+\alpha_m (-4 \alpha_f (\gamma-1)+6
        \gamma-5)\right)}{4 \beta^2}\bold{U}_n\\
    + 
    \left[{P^{1}_{1}} \right]^2 \alpha_m \Delta t
    \left(\frac{\alpha_f (3 \gamma-2) }{3 \beta^2}-\frac{7 \gamma
      }{6 \beta^2}+\frac{3 }{4 \beta^2}-\frac{\alpha_f }{\beta
        \alpha_m}+\frac{1}{\beta}\right)\dot{\bold{U}}_n \\+  
    \dfrac{\left[{P^{1}_{1}} \right] \Delta t}{\beta}
    \left(-\alpha_f -\gamma +\alpha_m
      +\frac{1}{2}\right)\ddot{\bold{U}}_n = 0 . 
  \end{multline}
  Next, we obtain second-order accuracy by canceling the third term in the last equation; thus, we set
  \begin{equation}
    \gamma=\frac{1}{2}-\alpha_{f}+\alpha_m.
  \end{equation}
  This completes our proof.
  \begin{remark}
    For values of $N>2$, we can follow a similar approach as the higher-order terms of $\left[{P^{N}_{1}} \right]$, noting that these terms do not affect the accuracy order.  
  \end{remark}
\end{proof}
\subsection{Stability analysis and CFL conditions}\label{sec:anal}

We can control the solution stability by bounding the spectral radius of the method; that is, the absolute value of the eigenvalues of the amplification matrix must be bounded by one. First, we calculate the eigenvalues of the amplification matrix for hyperbolic problems in~\eqref{eq:ampm} for the case $ \Delta t^{2}\omega_{i}^{2} \to 0$ as:
\begin{equation}\label{eq:t0}
  \lambda_1=\lambda_2=1,\qquad \lambda_3=\frac{\alpha_m-1}{\alpha_m}.
\end{equation}
The bounding  $\lambda_3$ in~\eqref{eq:t0} implies that $ {\alpha_m}\geq \frac{1}{2} $. 

We derive our method's CFL condition~\cite{ Courant56}. The implicit generalized-$\alpha$ method proposed in~\cite{ chung1993time, behnoudfar2019higher} are unconditionally stable. For these methods we study their high-frequency eigenvalues in the limit where $\Theta:= \Delta t^{2}\omega_{i}^{2} \to \infty$ where $\omega_{i}$ denotes the corresponding eigenvalue of the spatial discretization (eigenvalue of $\mathbb{M}^{-1}\mathbb{K}$). Thus, to guarantee unconditional stability, we set the parameters such that the eigenvalues are equal to a constant $\rho_\infty \in [0,1]$, which controls the numerical dissipation.

In our predictor multi-corrector method, we need to find a stable region due to its conditional behavior. Therefore, we only consider two of the roots of the characteristic polynomial as principal roots, while the third is spurious.  We maximize high-frequency dissipation by enforcing the principal roots to remain complex conjugate except in the high-frequency limit; the principal roots bifurcate only in the high-frequency limit.  Here, we denote the high-frequency limit as the critical stability limit, $\Omega_s$. Root bifurcation decreases the high-frequency dissipation; thus, our method's two limits of concern are $\Omega_s$ and the bifurcation limit $\Omega_b$.  The characteristic equation of the amplification matrix reads
\begin{equation}\label{eq:char}
  \sum_{j=0}^{3}\left(\tilde{a}^j+\tilde{y}^j \Theta^2\right)\lambda^{(3-j)}=0,
\end{equation}
where $\tilde{a}$ and $\tilde{y}$ depend on parameters such as $\gamma$ and $\beta$. The characteristic polynomial of $\mathbb{G}$ in~\eqref{eq:ampm} is
\begin{multline}\label{eq:charG}
  -\lambda^3+\lambda^2 \left(-\frac{\alpha_f^2 \left[{P^{1}_{1}}
      \right]}{\beta}+\frac{\alpha_f \alpha_m \left[{P^{1}_{1}}
      \right]}{\beta}+\frac{\alpha_f \left[{P^{1}_{1}}
      \right]}{\beta}-\frac{\alpha_f  \left[{P^{1}_{1}}
      \right]}{\alpha_m}-\frac{1}{\alpha_m}+\left[{P^{1}_{1}}
    \right]+3\right)\\ 
  +\lambda \left(\frac{2 \alpha_f^2 \left[{P^{1}_{1}}
      \right]}{\beta}-\frac{2 \alpha_f \alpha_m \left[{P^{1}_{1}}
      \right]}{\beta}-\frac{2 \alpha_f \left[{P^{1}_{1}}
      \right]}{\beta}+ 
      \frac{\alpha_m \left[{P^{1}_{1}}
      \right]}{\beta}+\frac{\left[{P^{1}_{1}} \right]}{\beta}+\frac{2
      \alpha_f \left[{P^{1}_{1}}
      \right]}{\alpha_m}+\frac{2}{\alpha_m}-2 \left[{P^{1}_{1}}
    \right]-3\right)\\ 
  -\frac{\alpha_f^2 \left[{P^{1}_{1}} \right]}{\beta}+\frac{\alpha_f
    \alpha_m \left[{P^{1}_{1}} \right]}{\beta}+\frac{\alpha_f
    \left[{P^{1}_{1}} \right]}{\beta}-\frac{\alpha_m \left[{P^{1}_{1}}
    \right]}{\beta}-\frac{\alpha_f \left[{P^{1}_{1}}
    \right]}{\alpha_m}-\frac{1}{\alpha_m}+\left[{P^{1}_{1}} \right]+1=0 
\end{multline}
We set two eigenvalues to $\rho_{b}$ and one becomes $\rho_{s}$. Thus, we rewrite the characteristic polynomial~\eqref{eq:char} as
\begin{equation}\label{eq:root}
\left(\lambda^2+\rho_{b}^2+2\lambda\rho_{b}\right)\left(\lambda+\rho_{s}\right)=0,
\end{equation}
where $\lambda$, $\rho_{b}$, and $\rho_{s}$ are the eigenvalues, and user-defined values, respectively.  The bifurcation limit $\Omega_{b}$ shows when the largest eigenvalue decreases from one and converges to the user-defined parameters $\rho_{b}$. Following a similar argument to~\cite{ HULBERT1996175}, having all three roots with real values at the bifurcation limit, we equate~\eqref{eq:root} and~\eqref{eq:char} and solve for the parameters $\alpha_m, \left[ P^{1}_{1} \right] ,\,\beta$. Then, we obtain two parameters
\begin{align}\label{eq:par}
\alpha_m=&\frac{-{\rho_b} {\rho_s}+{\rho_s}+2}{{\rho_b} {\rho_s}+{\rho_b}+{\rho_s}+1},
\end{align}
and
\begin{align}\label{eq:parbeta}
  \beta = \frac{
  \splitfrac{
  \alpha_f \left(\rho_b+1\right)
  \left\{\alpha_f \left( \rho_b+1\right)
  \left(\rho_s+1\right) \left[\left(\rho_b-1\right)
  \rho_s-2\right]
  -\rho_b^2
  \rho_s+2 \rho_b \left(\rho_s^2+\rho_s-1\right)-2 \rho_s^2-7
  \rho_s-6\right\}
  } { 
  + 2 \rho_b^2 \rho_s+\rho_b \left(\rho_s^2+2 \rho_s-3\right)-\rho_s^2-4 
  \rho_s-5 } }
  {\left(\rho_b+1\right)^2 \left(\rho_s+1\right)
  \left[\alpha_f (\rho_b+1) (\rho_s+1)+(\rho_b-1) \rho_s-2\right]
  }.  
\end{align}
We optimize the method's dissipation control using the available parametrization. We set $\rho_b=\rho_s$ and, for the sake of brevity, we omit many algebraic manipulations.  Therefore, our method reduces to a one-parameter technique with the following bifurcation limit $\Omega_{b}$
\begin{equation}
  \sum_{k=1}^{N}\left(\dfrac{\beta}{\alpha_m} \Omega_b\right)^k
  =-\frac{(\rho_{b}+1) \left[\alpha_f^2 (\rho_{b}-2) (\rho_{b}+1)^2+\alpha_f
      \left(-3 \rho_{b}^2+3 \rho_{b}+6\right)+3 \rho_{b}-5\right]}{\alpha_f
    \rho_{b}+\alpha_f+\rho_{b}-2}. 
\end{equation} 
We determine the stability region by finding the conditions under which the spectral radius grows to one after its bifurcation. For this, we propose the following
\begin{equation}\label{eq:stab}
  \sum_{k=1}^{N}\left(\dfrac{\beta}{\alpha_m} \Omega_s\right)^k
  =-\frac{\beta (4-8 \alpha_m)}
  {4 \beta (\alpha_f-\alpha_m)+2 \alpha_m
    \left[2\alpha_f- (\alpha_f-\alpha_m-1)+\alpha_m\right]+\alpha_m}. 
\end{equation}
Letting $N=1$, we substitute the definitions~\eqref{eq:par} into~\eqref{eq:stab} to obtain
\begin{equation}\label{eq:stabf}
  \Omega_s=-\frac{12 (\rho_{b}-2) (\rho_{b}+1)}{\rho_{b}^2-5 \rho_{b}+10}.
\end{equation}
\begin{remark}
  Following~\cite{ HULBERT1996175,  behnoudfar2018variationally, BEHNOUDFAR2021113656} and considering~\eqref{eq:stabf}, $\alpha_{f}$ is a free parameter with respect to spectral radius; thus, we obtain the optimal combination of low- and high-frequency dissipation by setting  $\alpha_{f}=0$ and letting the remaining parameters be:
  \begin{equation}
    \boxed{
      \begin{aligned}
        \alpha_m&= \frac{2-\rho_{b}}{\rho_{b}+1},\\
        \beta&= \frac{3 \rho_{b}-5}{(\rho_{b}-2) (\rho_{b}+1)^2},\\
        \Omega_b&= \frac{-3 \rho_{b}^2+2 \rho_{b}+5}{\rho_{b}-2}.
      \end{aligned}
    }
    \label{eq:Galphaparametersprev}
  \end{equation}
\end{remark}
\begin{remark}
  Setting $\rho_{b}(=\rho_b=\rho_s)=1$ leads to the largest stability region $\Omega_s=4$ which is equivalent to the second-order central difference method~\cite{ hairer2010solving}. This stability region imposes the following CFL condition $\Delta t^2\leq\dfrac{4}{\omega_{i}}$ on our method.
\end{remark}
We perform a similar analysis for higher-order methods, that is, $N>1$. Herein, we analyze the case of choosing $N=2$. We first determine the stability region as
\begin{equation}
  \Omega_s= \dfrac{\alpha_m-\alpha_m \sqrt{1-\dfrac{4\beta (4-8
        \alpha_m)}{4 \beta (\alpha_f-\alpha_m)+2 \alpha_m (2\alpha_f-
        (\alpha_f-\alpha_m-1)+\alpha_m)+\alpha_m}}}{2 \beta}. 
\end{equation}
Here, we cannot set $\alpha_f=0$; then, we have
\begin{equation}
  \boxed{
    \begin{aligned}
      \alpha_m&= \frac{2-\rho_{b}}{\rho_{b}+1},\\
      \beta&=\frac{\alpha_f^2 (\rho_{b}-2) (\rho_{b}+1)^2+\alpha_f \left(-3 \rho_{b}^2+3 \rho_{b}+6\right)+3 \rho_{b}-5}{(\rho_{b}+1)^2 (\alpha_f \rho_{b}+\alpha_f+\rho_{b}-2)},\\
      \alpha_f&=\dfrac{4}{1+4\rho_{b}}
    \end{aligned}
  }
  \label{eq:Galphaparameters}
\end{equation}

 \section{ Dynamics of hyperelasticity}
\label{section:nonlinear}

\begin{figure}[ht!]
  \centering
  {\includegraphics[height=9.5cm]{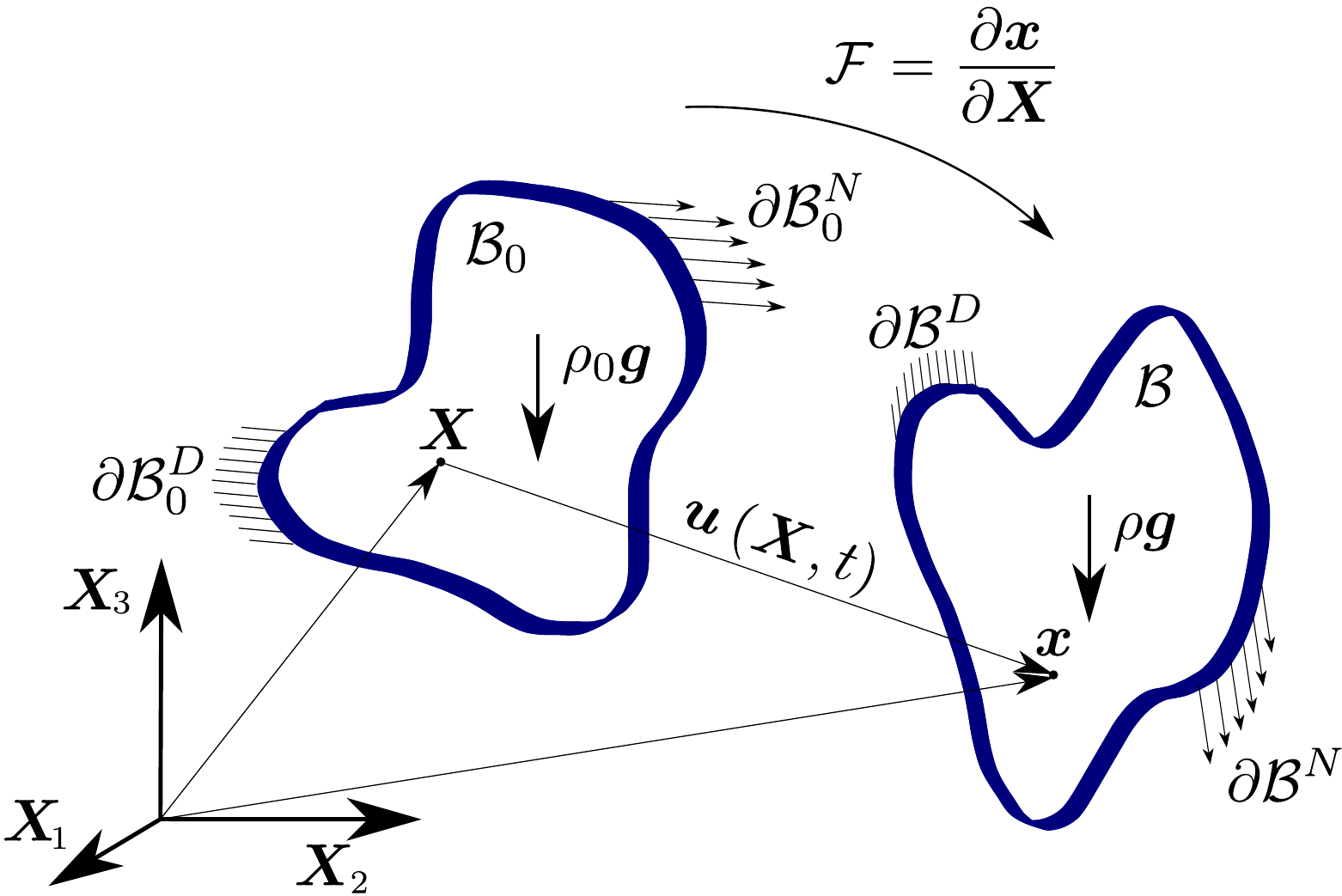}}
  \caption{Notation:  body undergoing large translations and deformations.}
  \label{fig:figure}
\end{figure}

In this section, bold variables are vector-valued functions. The hyperelasto-dynamical problem that Figure~\ref{fig:figure} sketches, reads: \textit{Find the one-to-one vector field going from the reference configuration $\boldsymbol{X} \in \mathcal{B}_{0}$  to the current configuration  $\boldsymbol{x} \in \mathcal{B}$
  \begin{equation}
    \boldsymbol{u}\left( \boldsymbol{X} , t \right) : \mathcal{B}_{0} \to \mathcal{B} \subset \mathbb{R}^{d},
  \end{equation}
  that fulfills the following partial differential equations as well as the intial and boundary conditions}
\begin{equation}
  \left\{
    \begin{aligned}
      \rho_{0} \ddot{\boldsymbol{u}}  - \nabla_{\boldsymbol{X}} \cdot
      \mathcal{P}\left( \boldsymbol{u} \right) &= \rho_{0}
      \boldsymbol{g}  && \text{in} \ \ \mathcal{B}_{0} \times I, \\ 
      \boldsymbol{u} \left( \boldsymbol{X} , t \right) & = \boldsymbol{u}_D
      \left( \boldsymbol{X} , t \right)  && \text{on}  \ \ \partial
      \mathcal{B}^{D}_{0} \times I,  \\ 
      \mathcal{P}\left( \boldsymbol{u} \right) \cdot \boldsymbol{n} &
      =\boldsymbol{t} \left( \boldsymbol{X} , t \right) && \text{on}  \ \
      \partial \mathcal{B}^{N}_{0} \times I, \\ 
      \boldsymbol{u} \left( \boldsymbol{X} , 0 \right)   &= \boldsymbol{u}_0
      \qquad
      \dot{\boldsymbol{u}}  \left( \boldsymbol{X} , 0 \right)   =
      \dot{\boldsymbol{u}}_0  && \text{on} \ \   \mathcal{B}_{0}, 
    \end{aligned}\right.
  \label{eq:elastodynamic}
\end{equation}%
where $\nabla_{\boldsymbol{X}} \left( \cdot \right)$ is the gradient in the reference configuration. In~\eqref{eq:elastodynamic}, $\rho_{0}$ is the initial material density, $\partial \mathcal{B}^{D}_{0}$ and $\partial \mathcal{B}^{N}_{0}$ are the Dirichlet and Neumann part of the body boundary $\partial \mathcal{B}_{0} = \partial \mathcal{B}^{D}_{0} \cup \partial \mathcal{B}^{N}_{0}$ respectively, $\boldsymbol{g}$ is the gravity's acceleration and $\mathcal{P}$ is the first Piola-Kirchhoff stress tensor defined as:
\begin{equation}
  \displaystyle  \mathcal{P} = \frac{\partial \,\mathsf{W}}{\partial \mathcal{F}},
  \label{eq:PiolaKir}
\end{equation}
with $\mathcal{F}$ being the deformation gradient tensor defined as:
\begin{equation}
  \mathcal{F} = \frac{\partial \boldsymbol{x}}{\partial\boldsymbol{X}}
  = \frac{\partial \left( \boldsymbol{X} + \boldsymbol{u}
    \right)}{\partial\boldsymbol{X}} = \mathcal{I} +
  \nabla_{\boldsymbol{X}} \boldsymbol{u} ,  
\end{equation}
where $\mathcal{I}$ is the second order identity and $\mathsf{W}$ is the stored-energy function \cite{SimoandHughes1998}, 
\begin{equation}
  \mathsf{W} = \mathsf{W}_{1} \left(\mathcal{F} \right) +
  \mathsf{W}_{2} \left(J \right), 
  \label{eq:storedenergy}
\end{equation}
that we express in terms of a  term that depends on the deformation gradient
\begin{equation}
  \mathsf{W}_{1} \left(\mathcal{F} \right) = \frac{\mu}{2} \left(
    \text{tr} \left( \mathcal{\mathcal{F}^{T} \mathcal{F}} \right) -3
  \right), 
\end{equation}
and another that depends on the volumetric deformation $J = \det \mathcal{F}$
\begin{equation}
  \mathsf{W}_{2} \left(J \right) = {\lambda}   \frac{\left(J^2 -1
    \right)}{4} -  \left( \frac{\lambda}{2} + \mu \right) \ln \left( J
  \right)  . 
\end{equation}
The stored-energy function $\mathsf{W}$ depends on two material parameters $\lambda , \mu > 0$, the Lamé constants,  while the second component that depends on the logarithm of the volumetric strain $J$ behaves like a penalization term that tends to prevent element flipping.  Combining~\eqref{eq:storedenergy} and~\eqref{eq:PiolaKir}, the first Piola-Kirchhoff stress tensor is
\begin{equation}
\displaystyle  \mathcal{P} \left( \boldsymbol{u} \right) = \lambda \frac{J^2 - 1}{2} \mathcal{F}^{-T} + \mu \left( \mathcal{F} - \mathcal{F}^{-T} \right).
\end{equation}
Our method solves the equation system explicitly in time; thus, we update our calculations using the non-symmetric stress tensor, which only affects the right-hand side. The formulation remains in terms of a mapping given by the deformation gradient.  Given this description of the physical problem, we adjust all these equations to fit them in our predictor multi-corrector time-marching scheme. Thus, for solving problem~\eqref{eq:elastodynamic}, we can rewrite the bi-linear form~\eqref{eq:HyperbolicPDEproblem} as:
\begin{equation}
  a \left( \boldsymbol{v}  ,  \boldsymbol{u}  \right) = \displaystyle
  \left( \nabla_{\boldsymbol{X}} \boldsymbol{v} , \mathcal{P} \left(
      \boldsymbol{u} \right)  \right)_{\mathcal{B}_{0}} ,  
  \label{eq:FuncAHyper}
\end{equation}
while the linear term is 
\begin{equation}
  \displaystyle  \left\langle  \boldsymbol{v} , \boldsymbol{f}
  \right\rangle =  \left\langle \boldsymbol{v}  , \rho_{0}
    \boldsymbol{g}  \right\rangle_{\mathcal{B}_{0}}  + \left(
    \boldsymbol{v}  , \boldsymbol{t} \right)_{\partial
    \mathcal{B}_{0}} .  
\label{eq:RHShyper}
\end{equation}
The first term in~\eqref{eq:HyperbolicPDEproblem}, the inertia term, is multiplied by the density as $\left\langle \boldsymbol{v} , \rho_{0} \ddot{\boldsymbol{u}}  \right\rangle$ in order to obtain a proper mass matrix. We obtain the bilinear functional using the Gateaux derivative~\cite{ HUGHES1978391, BatheRamm1975} on~\eqref{eq:FuncAHyper} to obtain:
\begin{equation}
  {a}' \left( \boldsymbol{v}, \lJump \boldsymbol{u} \rJump ;
    \boldsymbol{u}_{n} \right)  = \displaystyle \frac{d}{d\epsilon} a
  \left( \boldsymbol{v} , \boldsymbol{u}_{n} + \epsilon  \lJump
    \boldsymbol{u} \rJump  \right) \vert_{\epsilon=0} =  \left(
    \nabla_{\boldsymbol{X}} \boldsymbol{v},   \left[ \frac{\partial
        \mathcal{P} \left( \boldsymbol{u} \right)}{\partial \mathcal{F}}
    \right]_{n} \nabla_{\boldsymbol{X}} \lJump \boldsymbol{u} \rJump
  \right)_{\mathcal{B}_{0}} , 
  \label{eq:FuncAHypergateaux} 
\end{equation}
where the tangent constitutive hyperelastic tensor is 
\begin{equation}
  \displaystyle \frac{\partial \mathcal{P} \left( \boldsymbol{u}
    \right)}{\partial \mathcal{F}} =  \mathcal{F}^{-1} \left[ \lambda
    J^{2} \mathcal{I} \otimes \mathcal{I} + 2 \mu \left( 1 +
      \frac{\lambda}{2 \mu} \left( 1 - J^2 \right)\right)
    \mathfrak{1}  + \left( \mathcal{P} \left( \boldsymbol{u} \right)
      \mathcal{F}^{T} \right) \otimes \mathcal{I}\right]\mathcal{F}^{-T}, 
\end{equation}
where $\mathfrak{1} = \frac{1}{2} \left( \delta_{ac} \delta_{bd} + \delta_{ad} \delta_{bc}\right)$ denotes a symmetric fourth-order unit tensor with $ \delta_{ij}$ being a delta Kronecker function.  Finally, using~\eqref{eq:FuncAHyper},~\eqref{eq:RHShyper} and~\eqref{eq:FuncAHypergateaux}, we update~\eqref{eq:Galphau7} and~\eqref{eq:Galphadotu1} for the time-marching calculations.

 \section{Solution strategies and time-step size adaptivity}
\label{section:Implement}

\begin{figure} 
  \begin{center}
    \subfigure[(Interior) Gauss quadrature] {\includegraphics[width=.45\textwidth]{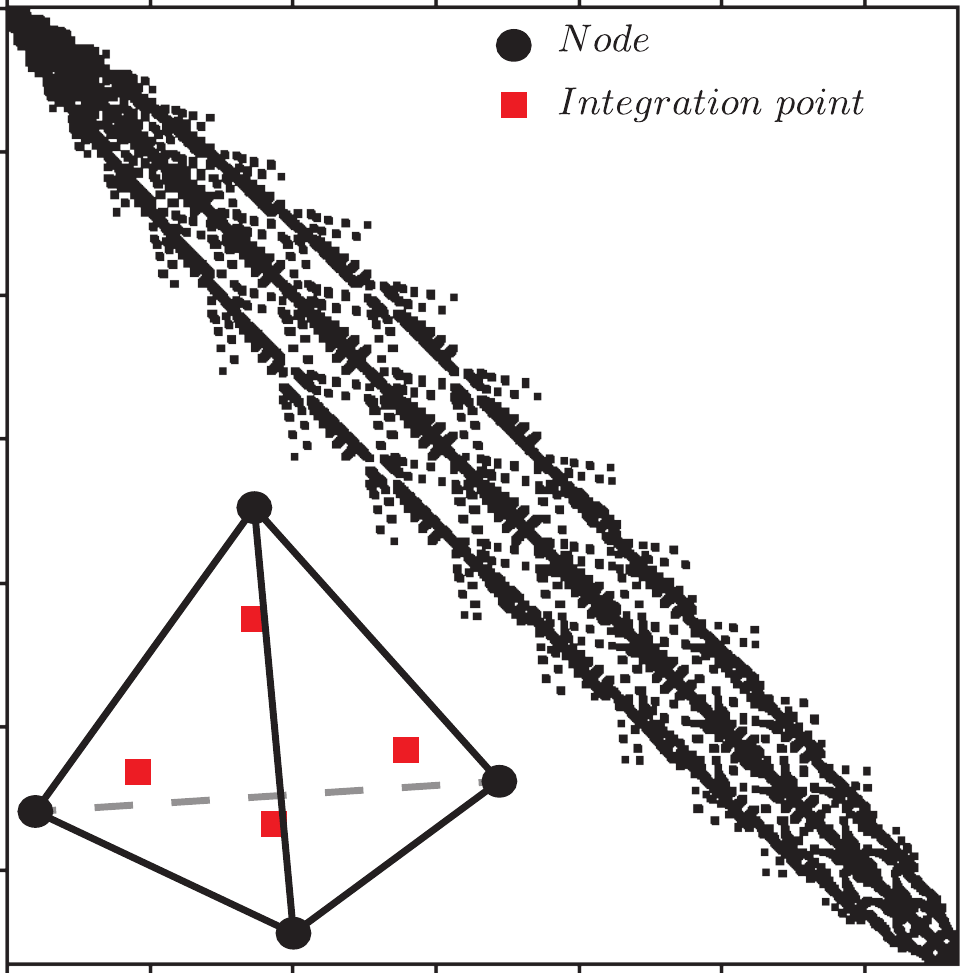}}
    \hfill
    \subfigure[(Nodal) Gauss-Lobatto quadrature] {\includegraphics[width=.45\textwidth]{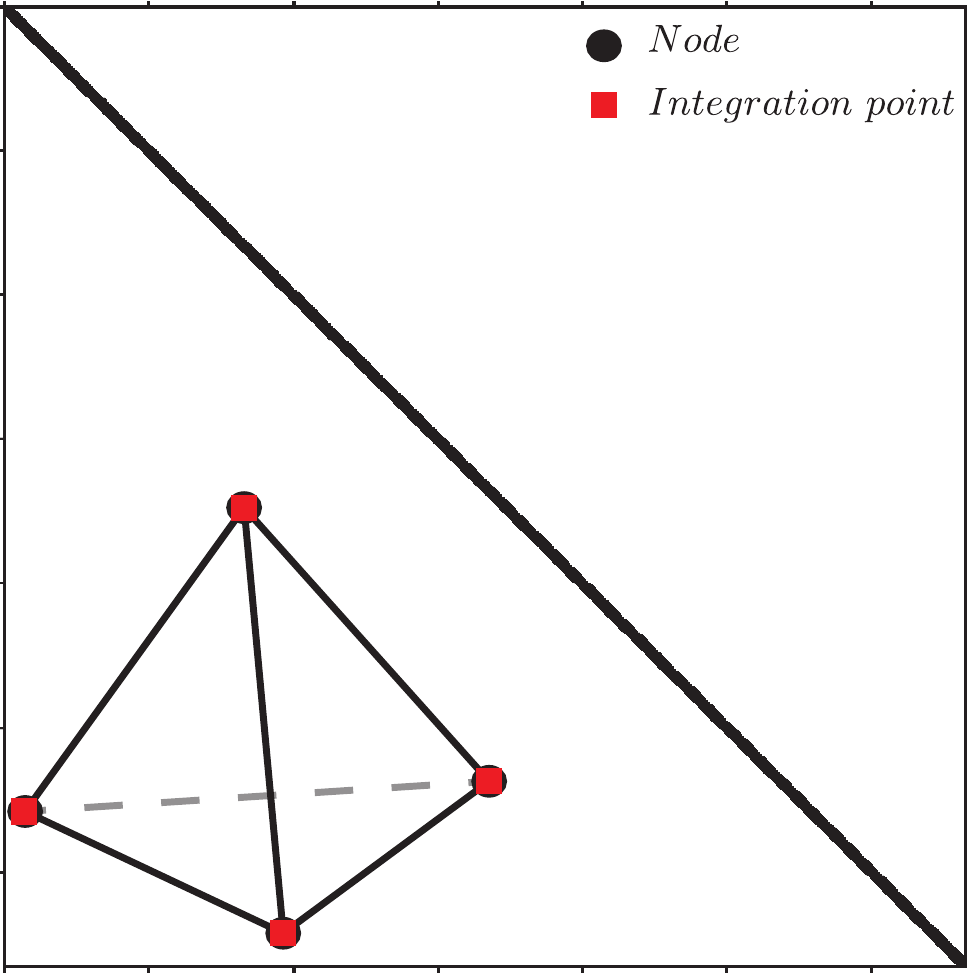}}
  \end{center}
  \caption{Mass matrix sparsity structure for a tetrahedral mesh using different integration rules}
\label{Quadrature}
\end{figure}

This section discusses time adaptivity and efficient variable updates for our method. We perform all the numerical experiments using the open-source FEniCS environment~\cite{ alnaes2015fenics} using the FIAT package~\cite{ Kirby2004} to integrate with different quadratures.  We estimate the truncation error to update the time-step size dynamically as the simulation progresses as in~\cite{ GOMEZ20084333}
\begin{equation}
  F\left(e;\Delta t; \text{tol} \right) = \rho_{\text{tol}} \left( \frac{\text{tol}}{e} \right)^{\frac{1}{d}} \Delta t,
  \label{TimeUpdate}
\end{equation}
where $ \rho_{\text{tol}}$ is a safety coefficient with the tolerance $\text{tol}$ taken from~\cite{ LANG1995223}, $d$ is the problem's spatial dimension and we estimate the truncation error $e$ using a summation expressed in terms of the update 
\begin{equation}
e = \frac{\| \delta \boldsymbol{u}^{k} \|_{L^2}}{\| \delta \boldsymbol{u}^{k-1}\|_{L^2}}  \ ,
\end{equation}
where $k$ is the term number of the updating formula.  In highly nonlinear problems, the time-step size can vary orders of magnitude to maintain the maximum eigenvalue of the amplification matrix lower than one. In our experience, the adaptive scheme avoids thousands of calculation steps while providing convergence at a low computational cost,  and no time-step bounding are needed since this is provided by the method itself. Algorithm~\ref{Alg1Adap} sketches the time marching scheme.  We make our explicit time-integrator attractive for engineering applications by using reduced integration for the mass matrix. Figure~\ref{Quadrature} shows the structure of a mass matrix calculated using different integration methods. Figure~\ref{Quadrature}(a) represents the sparsity obtained using Gauss quadrature, whereas Figure~\ref{Quadrature}(b) shows that, when we collocate integration points at the nodes; that is, we use a Gauss-Lobatto quadrature. Thus, the resultant mass matrix is diagonal. As a consequence, all updates only involve scalar-vector algebraic operations with a negligible computational cost. When we use the mass matrix integrated with Gauss-Lobatto quadrature as a preconditioner for the LGMRES solver, it converges in two iterations in serial computing. 	
\begin{algorithm}[h]  
\SetAlgoLined
\KwData{$\boldsymbol{u}_n$, $\dot{\boldsymbol{u}}_n$, $\ddot{\boldsymbol{u}}_n$ and $\Delta t_{n+1}$}
\KwResult{$\boldsymbol{u}_{n+1}$, $\dot{\boldsymbol{u}}_{n+1}$, $\ddot{\boldsymbol{u}}_{n+1}$ and $\Delta t_{n+1}$}
 Initialize $ \boldsymbol{u}_n  = \boldsymbol{u}_0$, $ \dot{\boldsymbol{u}}_n = \dot{\boldsymbol{u}}_0$ and $\ddot{\boldsymbol{u}}_n = \ddot{\boldsymbol{u}}_0$ for $t_n = 0$\;
 \While{$t_{n} \leq t_{f}$}{
  Update Matrix $\mathbb{K}^{*} =  -  \frac{\alpha_{f} \gamma_{2} \Delta t^{2}}{\alpha_{m} } \mathbb{K}$ and Mass Matrix $\mathbb{M}$\;
  Update Vector $\bold{H} $ \;
  Compute $\delta \boldsymbol{u}^{1}_{n+1} = \text{solve} \left(\mathbb{M},  \mathbb{K}^{*} \cdot \bold{H}\right)$ \;
  Compute $\delta \boldsymbol{u}^{2}_{n+1} = \text{solve} \left(\mathbb{M},  \mathbb{K}^{*} \cdot \delta \boldsymbol{u}^{1}_{n+1}\right)$ \;
  $\vdots$ \\
  Compute $\delta \boldsymbol{u}^{k}_{n+1} = \text{solve} \left(\mathbb{M},  \mathbb{K}^{*} \cdot \delta \boldsymbol{u}^{k-1}_{n+1} \right)$ \;
  Compute error $e_{n+1} = \frac{\| \delta \boldsymbol{u}^{k}_{n+1} \|_{L^2}}{\| \delta \boldsymbol{u}^{k-1}_{n+1} \|_{L^2}}$ \;
  \eIf{$e_{n+1} \leq \text{tol}$}{
   Update Current time-step $t_{n} \leftarrow t_{n} + \Delta t_{n+1}$\;
   Update $\lJump \boldsymbol{u} \rJump_{n+1} = \sum^{k}_{j=1} \delta \boldsymbol{u}^{j}_{n+1}$ \;
   Update $\boldsymbol{u}_{n+1}$, $\dot{\boldsymbol{u}}_{n+1}$ and $\ddot{\boldsymbol{u}}_{n+1}$\;
   Update $\boldsymbol{u}_{n} \leftarrow \boldsymbol{u}_{n+1}$, $\dot{\boldsymbol{u}}_{n} \leftarrow \dot{\boldsymbol{u}}_{n+1}$ and $ \ddot{\boldsymbol{u}}_{n} \leftarrow \ddot{\boldsymbol{u}}_{n+1}$\;
   \If{$e_{n+1} < \text{tol}_{min}$}{
   Increase time-step $\Delta t_{n+1} \leftarrow F\left(e_{n+1};\Delta t_{n+1}; \text{tol}_{min}\right)$}
   }{
   Reduce time-step $\Delta t_{n+1} \leftarrow F\left(e_{n+1};\Delta t_{n+1}; \text{tol}\right)$\;
  }
 }
 \caption{Explicit update with time-step size adaptativity}
 \label{Alg1Adap}
\end{algorithm}

  \section{Numerical Experiments}
\label{section:Examples}

In this section, we study several challenging applications that demonstrate the features of our predictor multi-corrector scheme. We analyze the method's convergence analysis in well-defined linear problems. We model finite deformation models, using the formulation of Section~\ref{section:nonlinear}.  In all numerical experiments, unless otherwise noted, we compute the residuals and the mass matrices using Gauss quadratures, which produce sparse matrices, but not diagonal ones.

\subsection{Convergence analysis}

\subsubsection{Wave equation}

\begin{figure}[h!]
  \begin{center}
    \subfigure[${L_2}$ error in ${u}$]{\includegraphics[height=5.0cm]{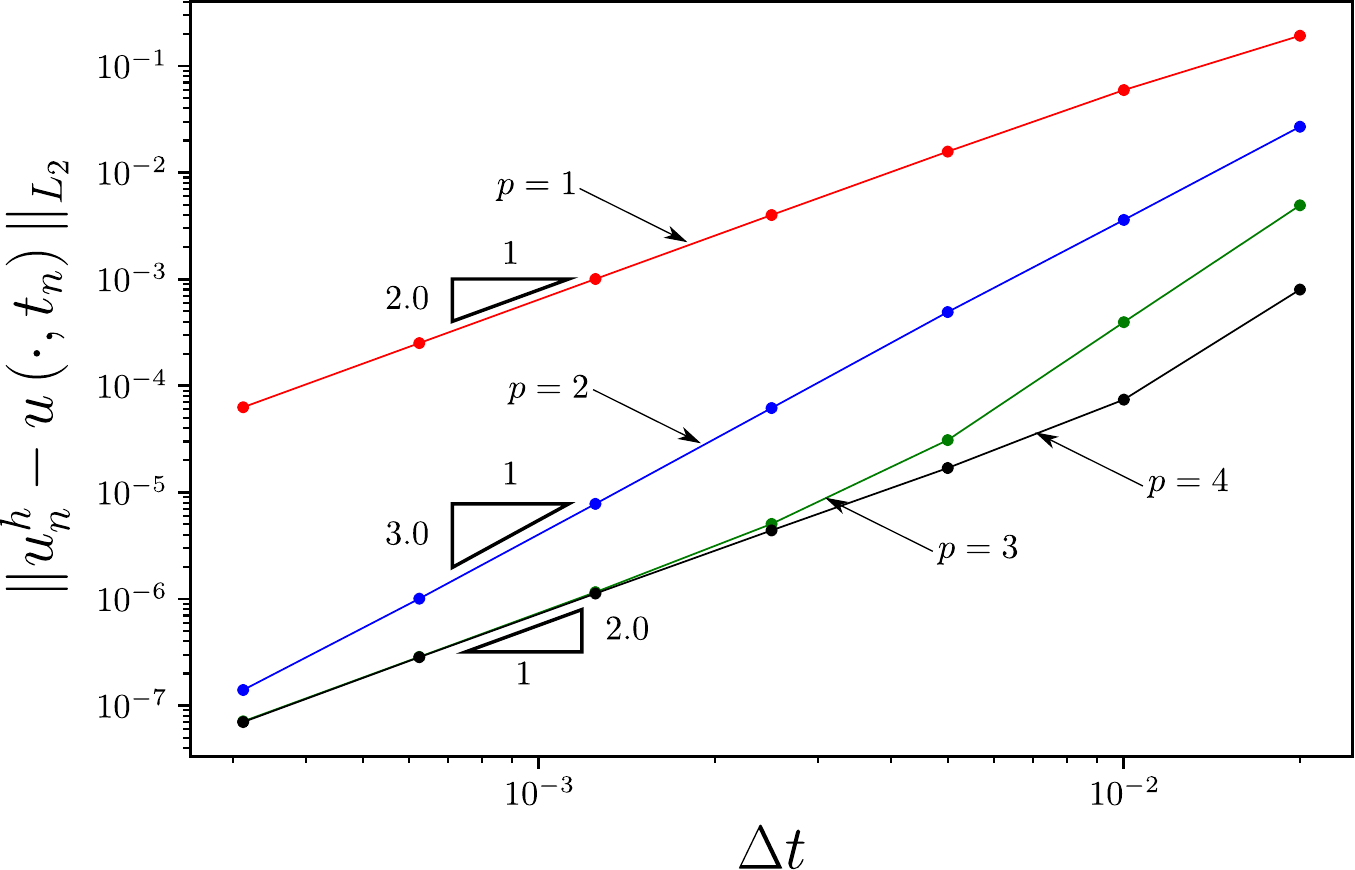}} \hspace{0.5cm}
    \subfigure[${L_2}$ error in $\dot{{u}}$]{\includegraphics[height=5.0cm]{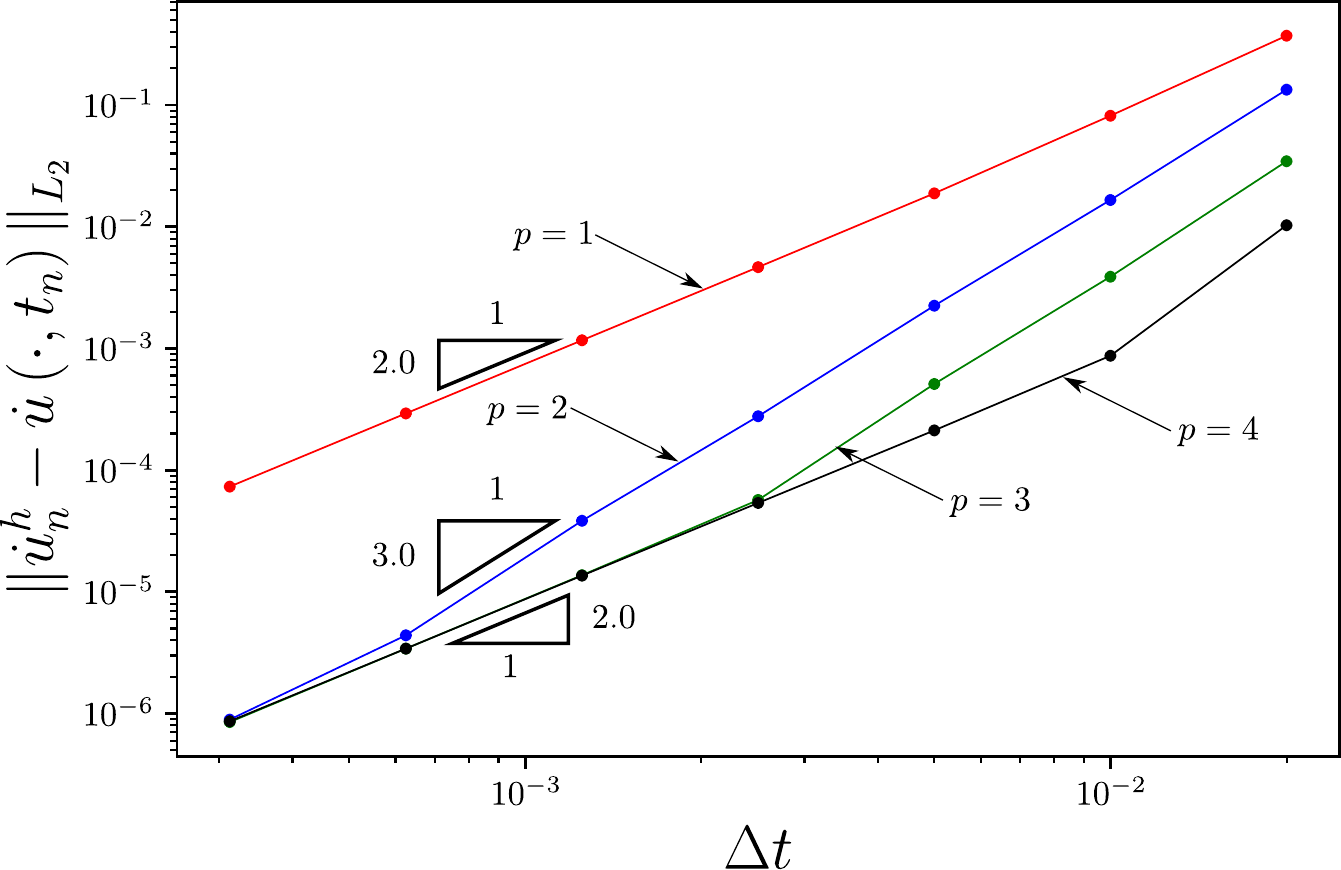}}
  \end{center}
  \caption{Wave equation temporal convergence for different polynomial orders}
  \label{fig:Waveconv}
\end{figure}

First, we present results for a linear wave-propagation equation. Assuming ${u}: \mathcal{B} \times I \rightarrow \mathbb{R}^{m}$ as a scalar field with $m=1$ in Eq.~\eqref{eq:HyperbolicPDE},  leads to the wave equation.  In this case, the forcing, initial, and boundary conditions ensure that the solution is:
\begin{equation}
  \left\{
    \begin{aligned}
      u \left(x,y,t\right) & = \cos \left( \pi x \right)\cos \left(
        \pi y \right)\cos \left(  \pi \sqrt{2} \ t \right), \\ 
      \dot{u} \left(x,y,t\right) &= - \pi \sqrt{2} \cos \left( \pi x
      \right)\cos \left( \pi y \right)\sin \left(  \pi \sqrt{2} \ t
      \right) ,\\ 
      \ddot{u} \left(x,y,t\right) &=- 2 \ \pi^2  \cos \left( \pi x
      \right)\cos \left( \pi y \right)\cos \left(  \pi \sqrt{2} \ t
      \right)  .
    \end{aligned}\right.
\end{equation}
In this example, we only use the first two terms of the updating formula~\eqref{eq:Galphau7}. Figure~\ref{fig:Waveconv} shows the temporal convergence for the scalar wave equation in the $L_2$ norm of ${u}$ and its time derivative $\dot{{u}}$ at the final time $t_f = 0.1$,  using different polynomial orders $p=1,...,4$ on the space discretization. The figures show second-order convergence in all cases.
 
\subsubsection{Hyperelasticity}

\begin{figure} [h!]
  \begin{center}
    \subfigure[${L_2}$ error in ${\boldsymbol{u}}$.]
    {\includegraphics[width=.6525\linewidth]{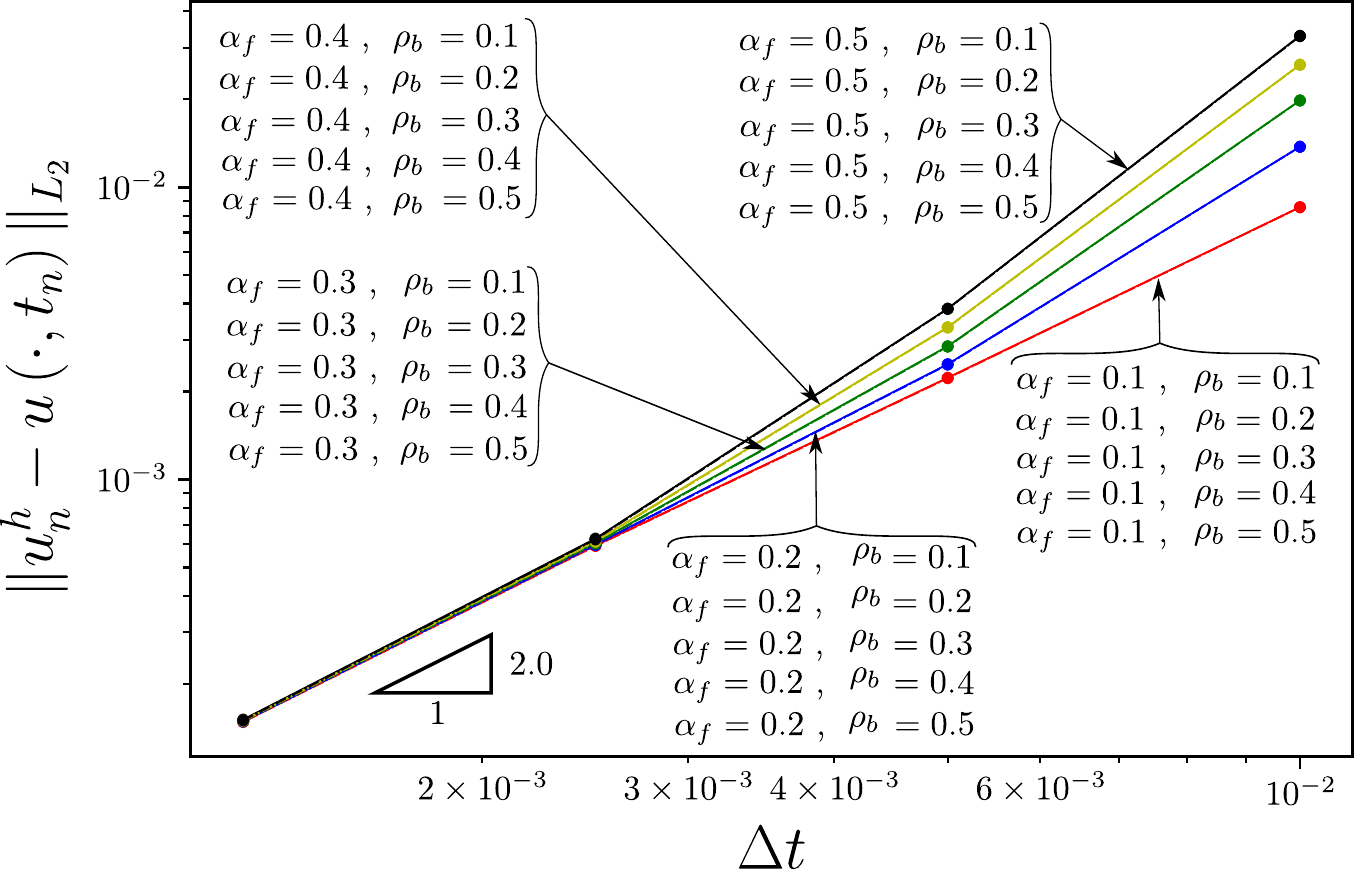}}
    \subfigure[${L_2}$ error in $\dot{\boldsymbol{u}}$.]
    {\includegraphics[width=.6525\linewidth]{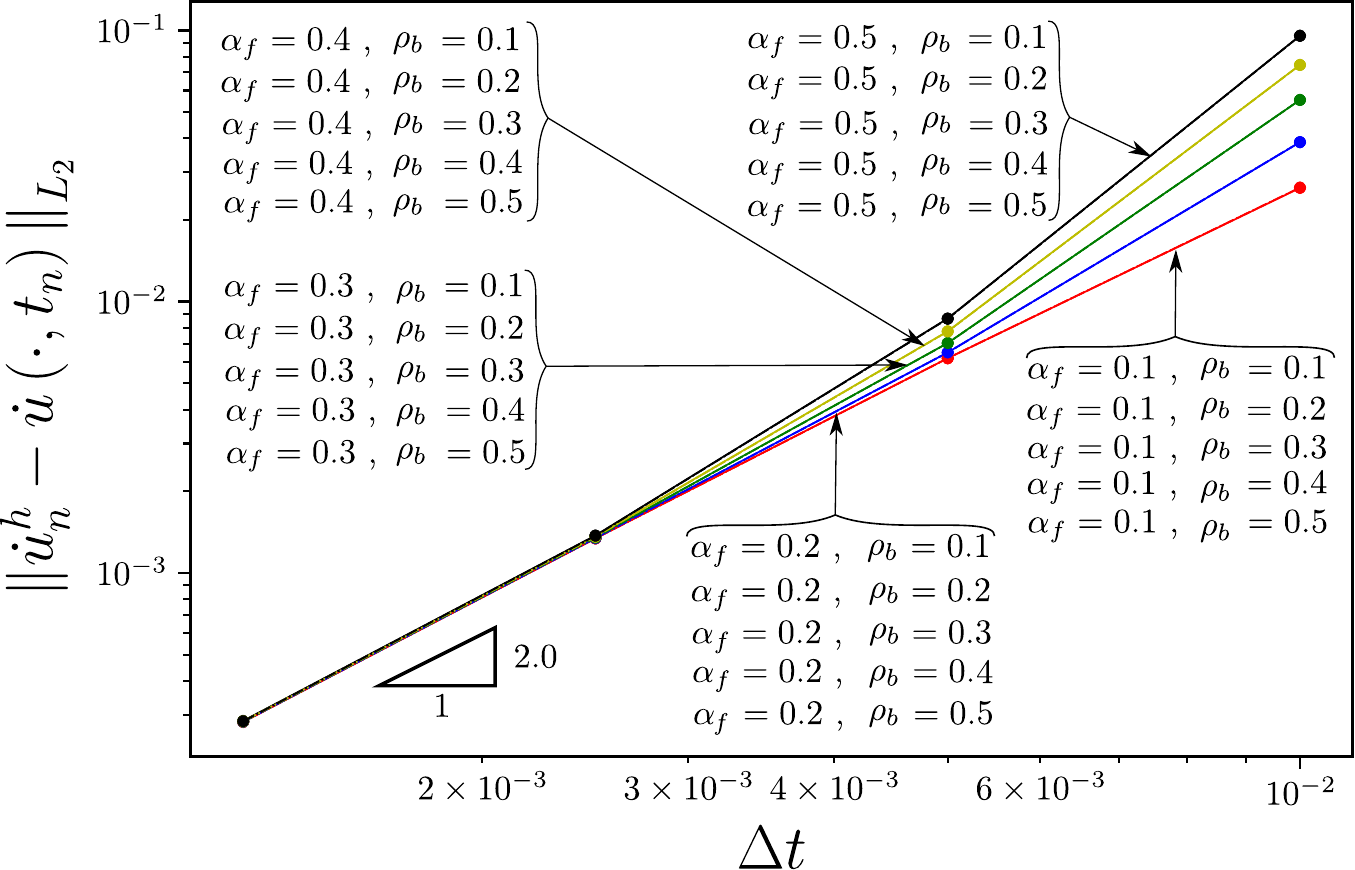}}
  \end{center}
  \caption{Convergence analysis for linear tetrahedra}
  \label{fig:firstorder_hyper}
\end{figure}

We analyze the convergence for the problem we describe in Section~\eqref{section:nonlinear}. Using the parameter definitions of~\eqref{eq:Galphaparametersprev}, we evaluate $\alpha_f$'s influence on the method convergence.  In this example, unitary material parameters, are used for young modulus and density, i.e.  $E = \rho_{0} = 1$, whereas the Poisson modulus is defined as $\nu = 0.2$.  The resulting system of partial differential equations is hyperbolic where the maximum eigenvalue of the amplification matrix is in the order of  $\max \lambda_{k} \approx \frac{h}{\Delta t}$. For the convergence analysis, we consider a cubic domain $\mathcal{B} = \left[0,1 \right] \times \left[0,1 \right] \times \left[0,1 \right] $ with the following exact solution 
\begin{equation}
  \boldsymbol{u} \left(x,y,z,t\right)   = U_0 \sin \left( \omega t \right)  \left[
    \begin{array}{c}
      -2 \sin \left({\pi} x \right) \cos\left( {\pi} y \right) \cos \left( {\pi} z \right)  \\
      \cos\left( {\pi} x \right) \sin\left( {\pi} y \right) \cos\left( {\pi} z \right)\\
      \cos\left( {\pi} x \right) \cos\left( {\pi} y \right) \sin\left( {\pi} z \right)
    \end{array}
  \right].
  \label{eq:exacthyperel}
\end{equation}
%
\begin{figure} [h!]
  \begin{center}
    \subfigure[$\| \cdot \|_{L_2}$ error in ${\boldsymbol{u}}$.]
    {\includegraphics[width=.6525\linewidth]{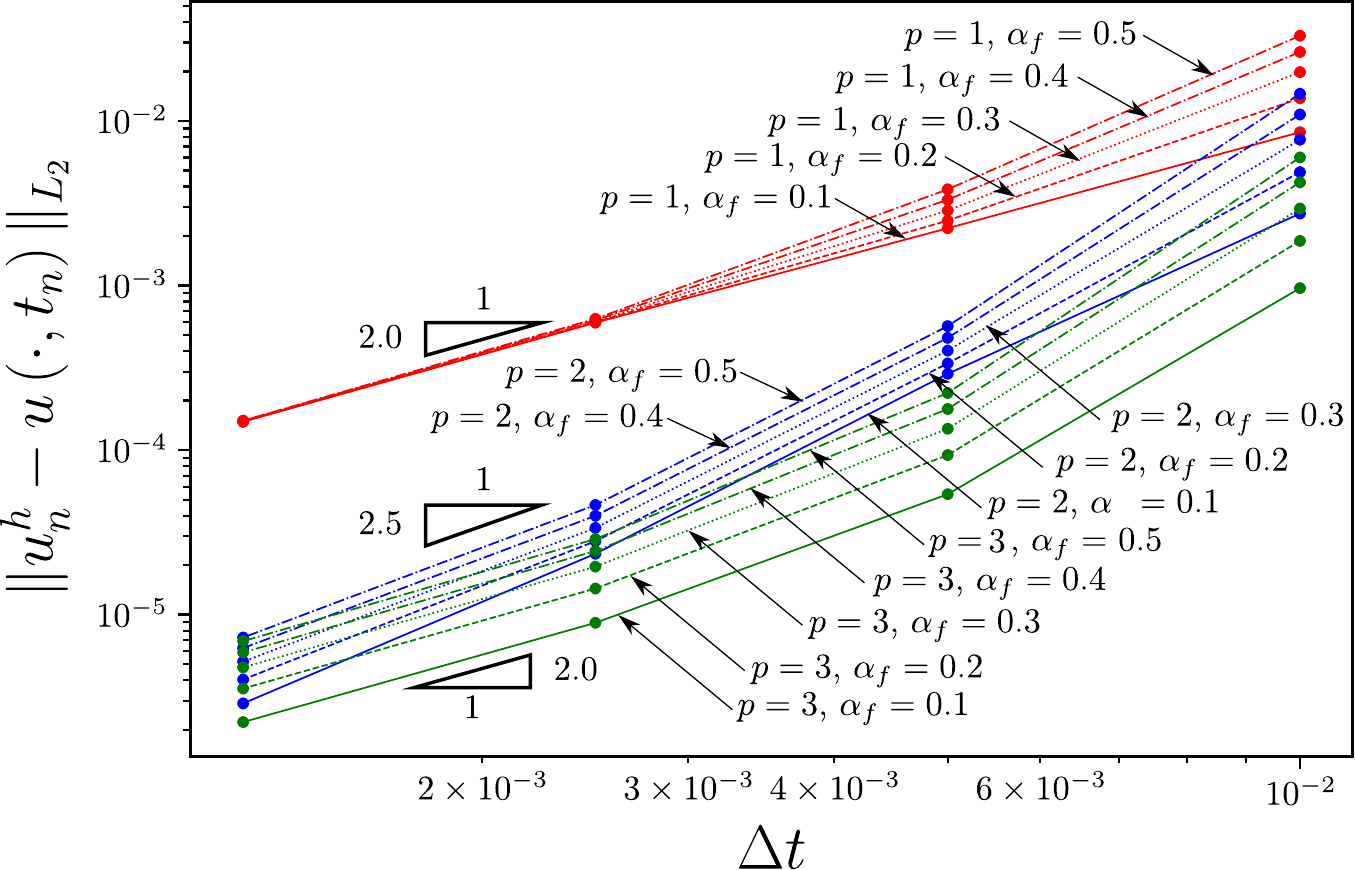}}
    \hspace{0.5cm} 
    \subfigure[$\| \cdot \|_{L_2}$ error in $\dot{\boldsymbol{u}}$.]
    {\includegraphics[width=.6525\linewidth]{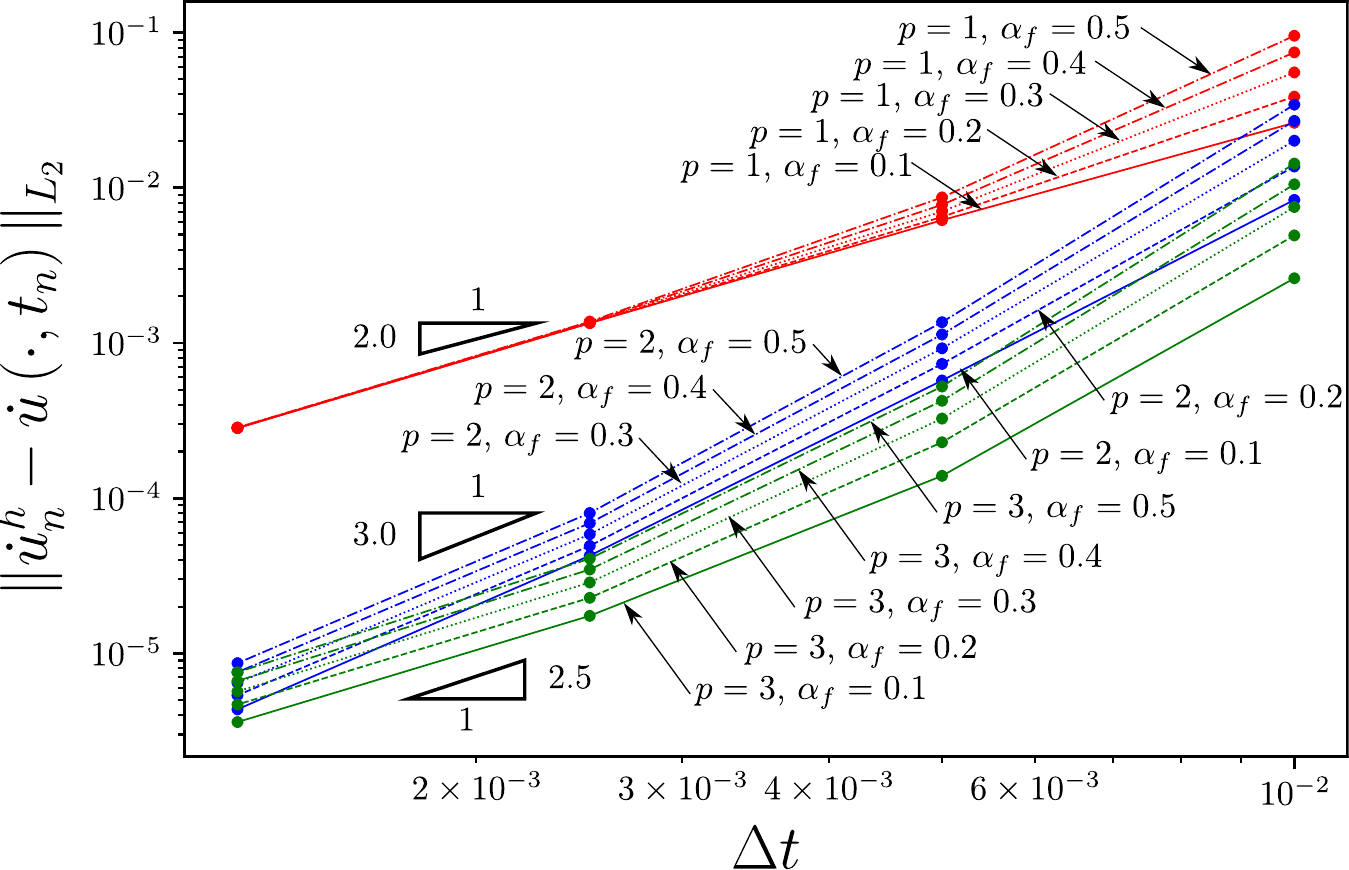}} 
  \end{center}
  \caption{Convergence analysis of linear, quadratice, and cubic tetrahedra for different $\alpha_f$ values}
  \label{fig:alphaf_hyper}
\end{figure}

Figure~\ref{fig:firstorder_hyper} shows the convergence results for linear tetrahedra with different $\rho_{b}$ and $\alpha_f$ values.  The method has second-order accuracy for both the displacement and velocity fields in the ${L_2}$ norm. Figure~\ref{fig:firstorder_hyper}  also shows that the dissipation parameter $\rho_{b}$ has a negligible influence on the error; the curves coincide as the time-step size decreases. Nevertheless, $\alpha_f$ influences the solution for large time-step sizes delivering a faster convergence in this region.

Figure~\ref{fig:alphaf_hyper} displays second-order accuracy in time for linear, quadratic, and cubic tetrahedra when using different $\alpha_f$ values. This results were computed using a fixed value $\rho_{b}$, since its influence in method's convergence is negligible.

\subsection{Dynamics of hyperelasticity}

\subsubsection{Release of a twisted-compressed bar under gravity loading}

\begin{figure} [h!]
  \begin{center}
    \subfigure[Problem setup]
    {\includegraphics[width=.325\linewidth,clip]{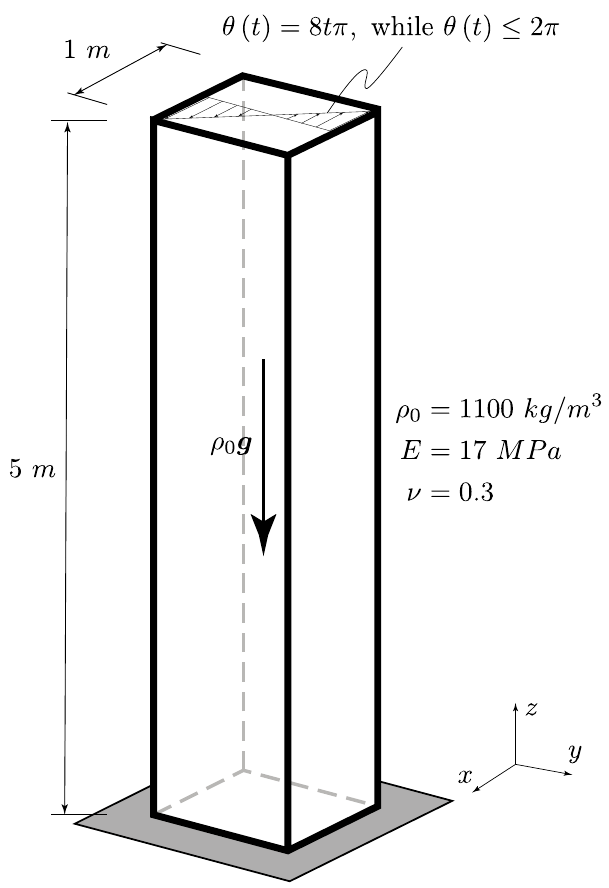}}
    \hspace{0.051cm}
    \subfigure[Potential, kinetic, and total energies evolution]
    {\includegraphics[width=.65\linewidth,clip]{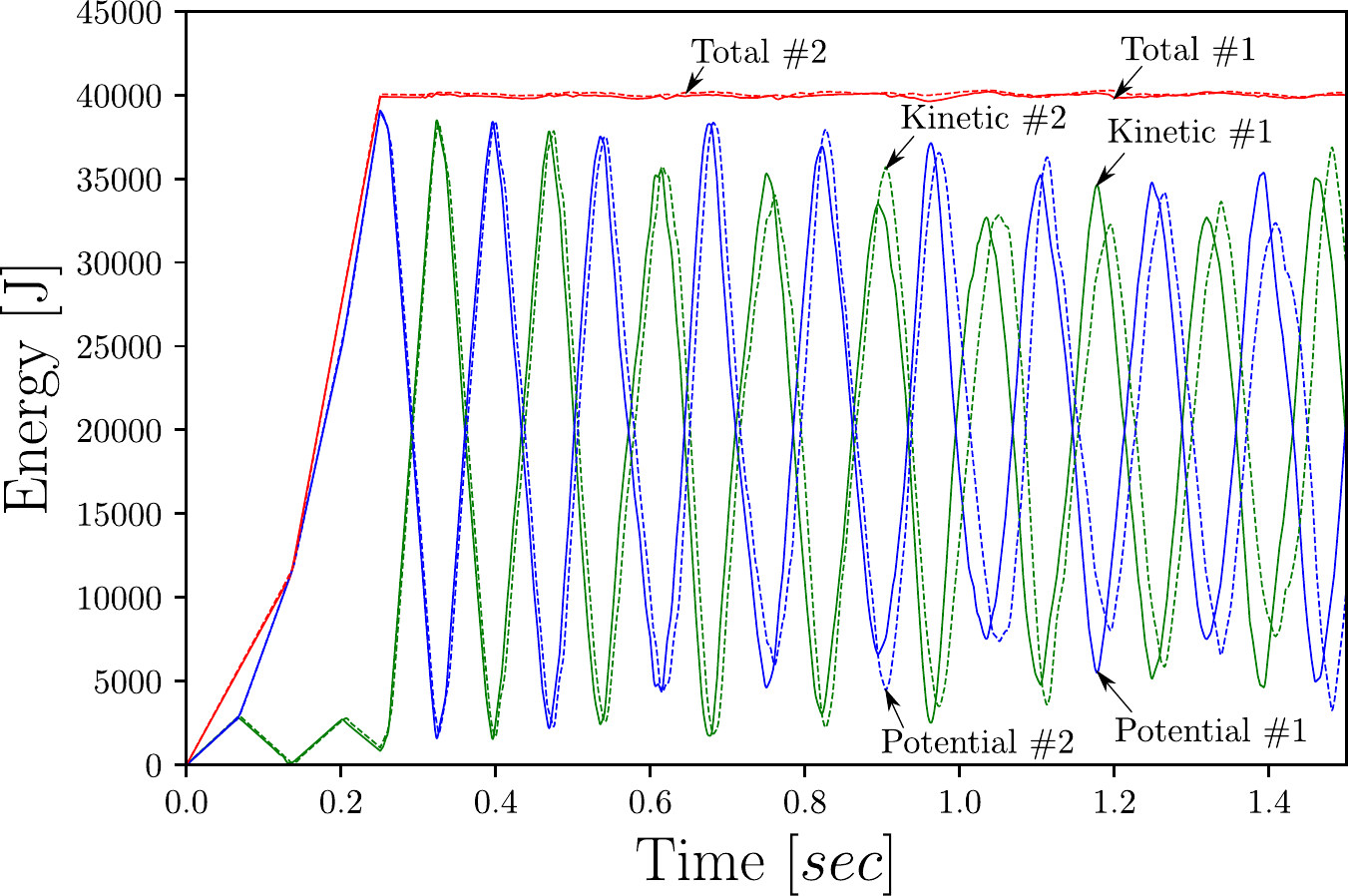}}
  \end{center}
  \caption{Twisted-compressed bar: Problem setup and energies
    evolution for Meshes $\# 1\ \&\ \# 2$}
  \label{fig:energyplot}
\end{figure}

\begin{figure} [h!]
  \begin{center}
    \subfigure[Mesh $\# 1$]
    {\includegraphics[width=.495\linewidth]{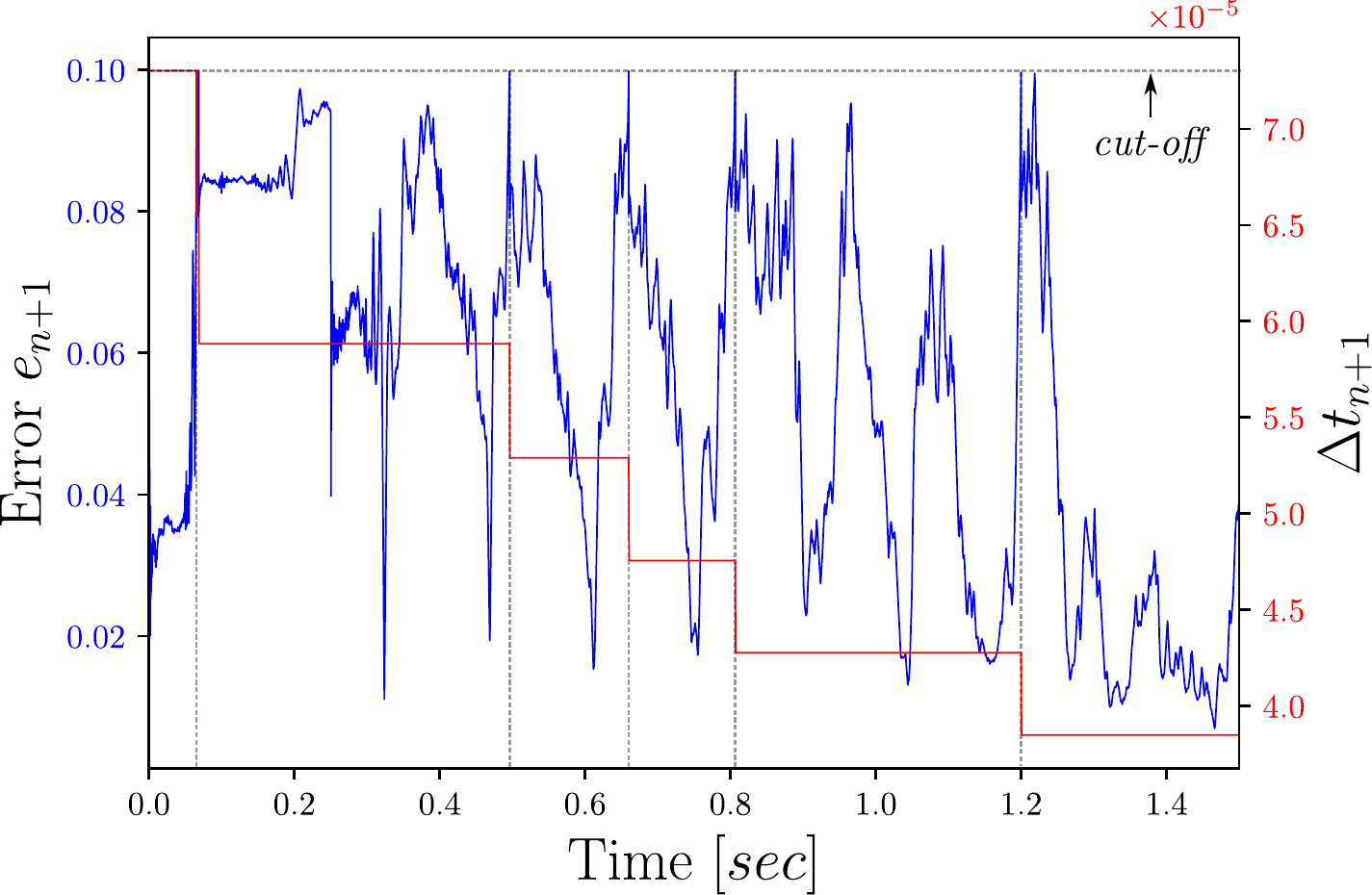}}
    \subfigure[Mesh $\# 2$]
    {\includegraphics[width=.495\linewidth]{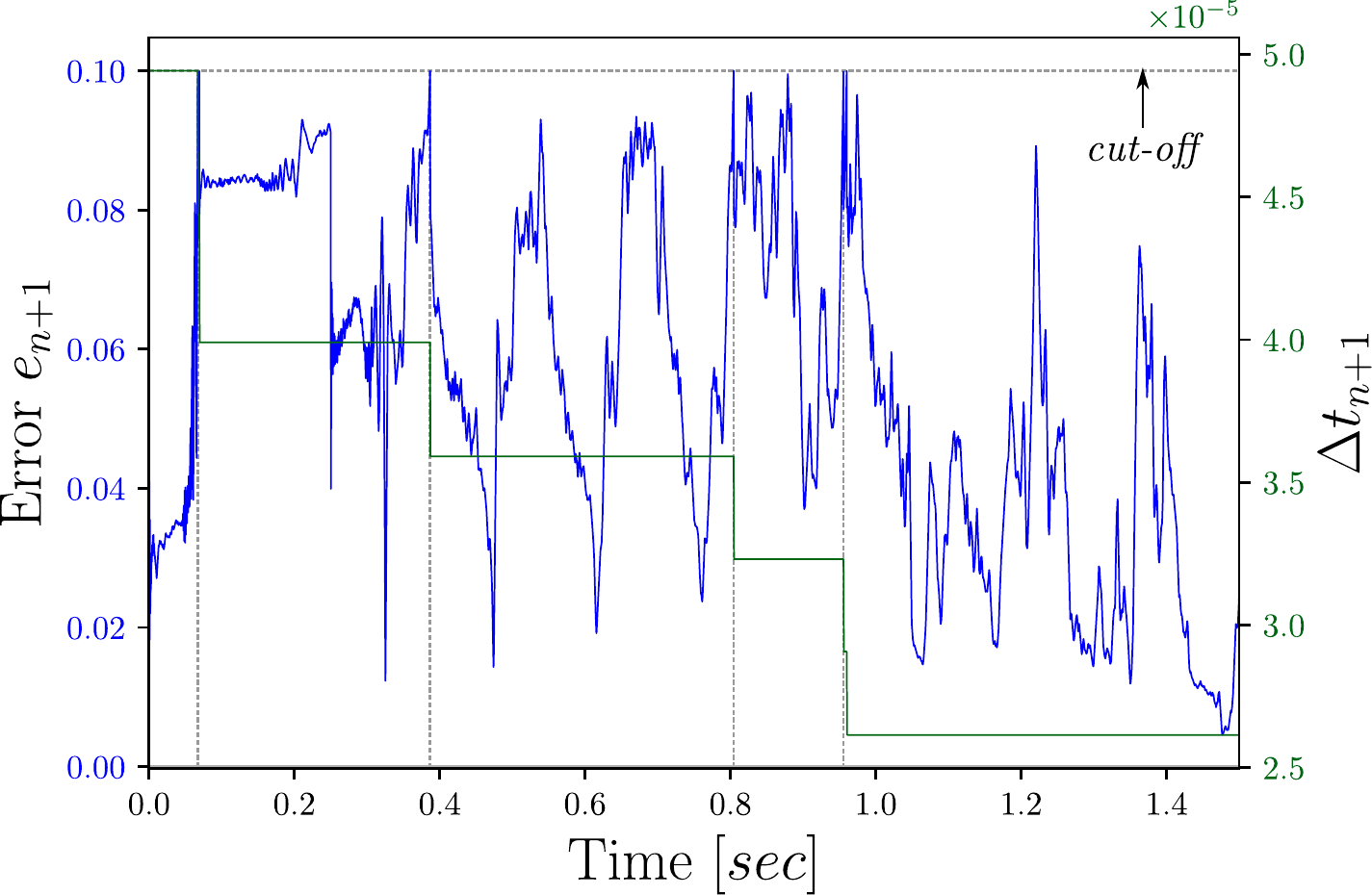}}
  \end{center}
  \caption{Adaptativity: Truncation error $e_{n+1}$ \& $\Delta t_{n+1}$ evolution vs time}
  \label{fig:Adaptivityplot}
\end{figure}

This experiment shows the robustness of our time adaptive strategy by simulating the abrupt changes in the dynamics of a hyperelastic material when boundary conditions change drastically. Figure~\ref{fig:energyplot} (a) shows the problem setup: a square-section bar of $5 \ \text{m}$ length and $1 \ \text{m}$ in each side of its section with initial density of $\rho_0 = 1100 \ \frac{\text{kg}}{\text{m}^3}$, a Young modulus $E = 17 \ \text{MPa}$, and a Poisson ratio of $\nu = 0.3$. The bar's weight compresses it under gravity loading $\rho_0 \bold{g}$,  using $\bold{g} = 9.81 \frac{\text{m}}{\text{s}^2}$,  and the parameter $\rho_{b}=1$ is assumed to avoid further numerical dissipation.   We prescribe a twisting displacement tangent to the top plane of the bar; this displacement follows the following time-dependent function
\begin{equation}
  \theta \left( t \right) = 8 t \pi \ , \ \text{for} \ t \leq 0.25 \ \text{s}.
\end{equation}
Once the bar reaches the maximum twisting of $2\pi$, we release it and simulate the force-free undamped oscillations that ensue. We perform the simulations using two meshes; Mesh $\#1$ has 90,000 DOFs, while Mesh $\#2$ has 246,960 DOFs. In this simulation, we use a cut-off tolerance for error $e_{n+1}$ of $tol = 0.1$,  and a safety factor of $\rho_{tol} = 0.9$ for the time adaptivity; we set the final time to $1.5 \text{s}$, and only use the first two terms of the updating formula~\eqref{eq:Galphau7} in each time increment. This update strategy is followed by standard explicit algorithms, where the nonlinear residual of the system is not considered before marching to the next step.

\begin{figure}[h!]
  \begin{center}
    {\includegraphics[width=.9\linewidth]{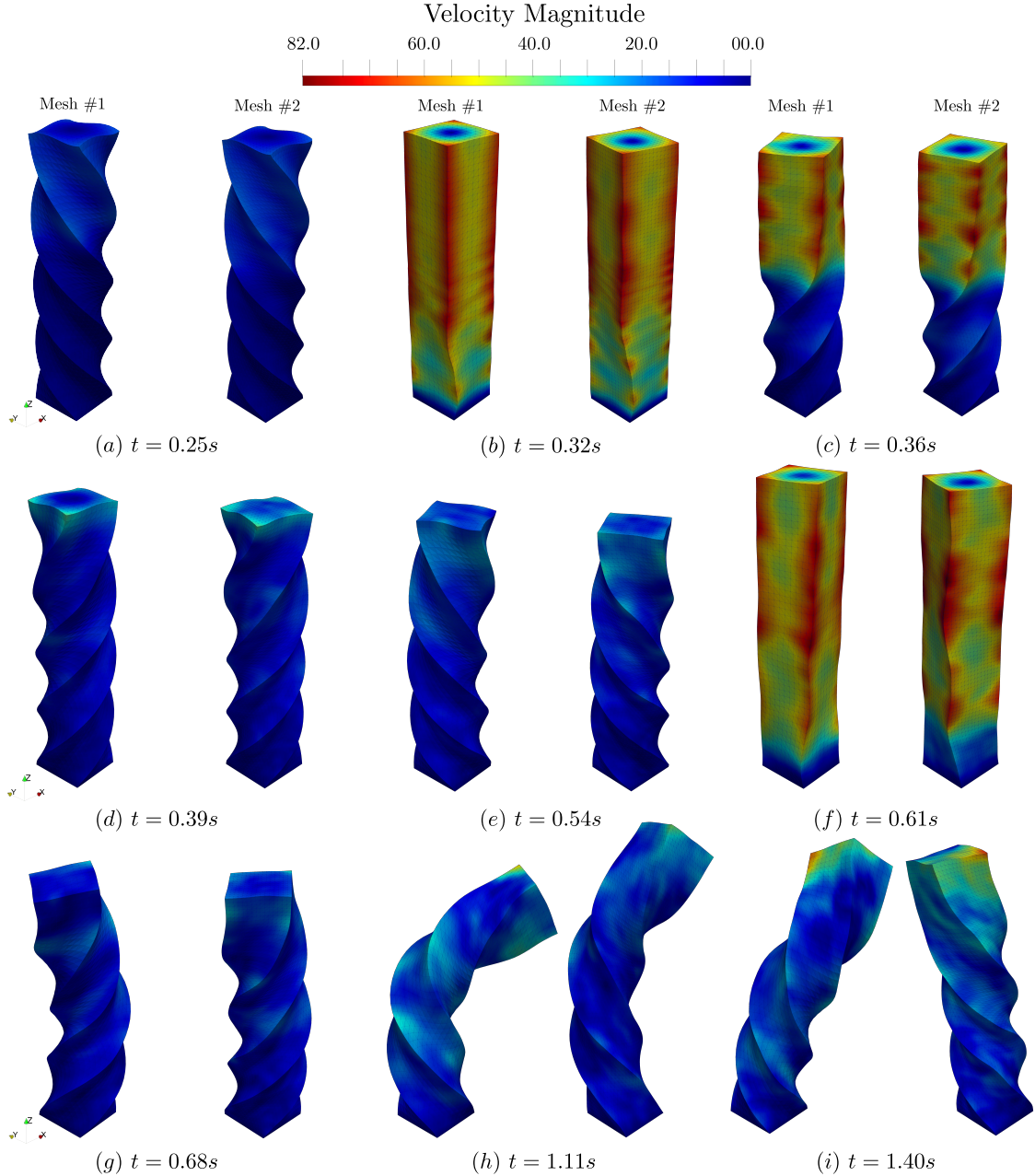}}
  \end{center}
  \caption{Bar deformation snapshots for Meshes $\# 1\ \&\ \# 2$}
  \label{fig:barsnapshots}
\end{figure}

Figure~\ref{fig:energyplot}~(b) shows the kinetic $K$, potential $P$, and total energies $T=K+P$ as functions of time for both meshes. We calculate the kinetic energy as
\begin{equation}
  {K} = \int_{\mathcal{B}_{0}} \frac{1}{2} \rho_0
  \dot{\boldsymbol{u}} \cdot \dot{\boldsymbol{u}} \ d\mathcal{B}_{0}, 
\end{equation}
and the potential energy as
\begin{equation}
 {P} = \int_{\mathcal{B}_{0}} \mathcal{W} \left( \mathcal{F}
  \right) \ d\mathcal{B}_{0}, 
\end{equation}
while the total energy is their sum. Until $t=0.25 \ \text{s}$, Figure~\ref{fig:energyplot} (b) shows that the twisting boundary condition increases the potential energy in the bar. After the twisting release, the system exchanges kinetic and potential energies, while their sum remains almost constant throughout the simulation.  It is worth to note that oscillations on the total energy are produced for computing just the first two terms on the updating formula,  reducing the time-marching proposal in a simple predictor-corrector method.  In following examples we will show how this issue is solved when further terms on the sum are included,  controlling the propagated error throughout the simulation.

Figure~\ref{fig:Adaptivityplot} shows the effectiveness of our time-adaptive approach for these meshes, where the blue line represents the truncation error $e_{n+1}$ for the specific time while the red one represents the time-step size $\Delta t_{n+1}$ for Mesh $\# 1$, while green line shows the time-step size for Mesh $\# 2$ in a different scale.  Both meshes have several sudden error increases; these rises reduce the time-step size to avoid the sum divergence. Without this adaptive time-stepping strategy, the simulation is unable to progress beyond $t=0.07 \ \text{s}$ in during the initial twisting regime.

Figure~\ref{fig:barsnapshots} displays snapshots of the deformed bar at different instants of the free undamped oscillation regime. Between $t=0.25 \ \text{s}$ to $t=0.68 \ \text{s}$, Figures~\ref{fig:barsnapshots}~(a) to~(g) show that the bar deformation follows a periodic twisting movement, combined with a compressing-stretching wave. Compression and stretching increase the error since the constitutive definition contains a logarithmic function that penalizes the system as $J \rightarrow 0$, avoiding element flipping.  From time $t=0.68 \ \text{s}$ onwards, the bar buckles, compressing the side elements and, thus, reducing further the time-step size.

\subsubsection{Impact dynamics of a hyperelastic tube}

\begin{figure} [h!]
  \begin{center}
    \subfigure[Problem setup]{\includegraphics[height=6.0cm]{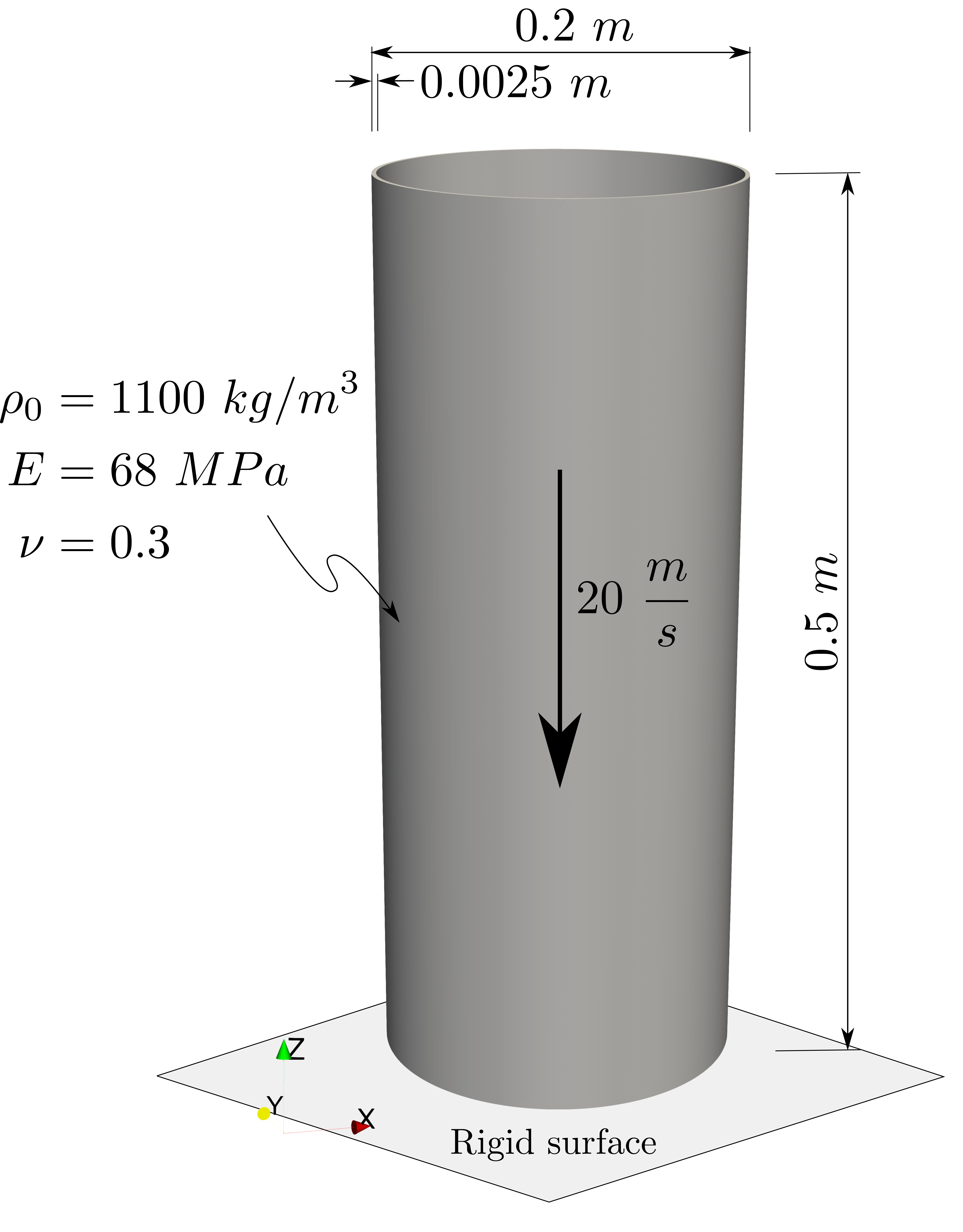}}
    \subfigure[$\sim350K$ tetrahedra]{\includegraphics[height=6.24cm, trim={0.0cm 0.0cm 0.0cm 0.0cm},clip]{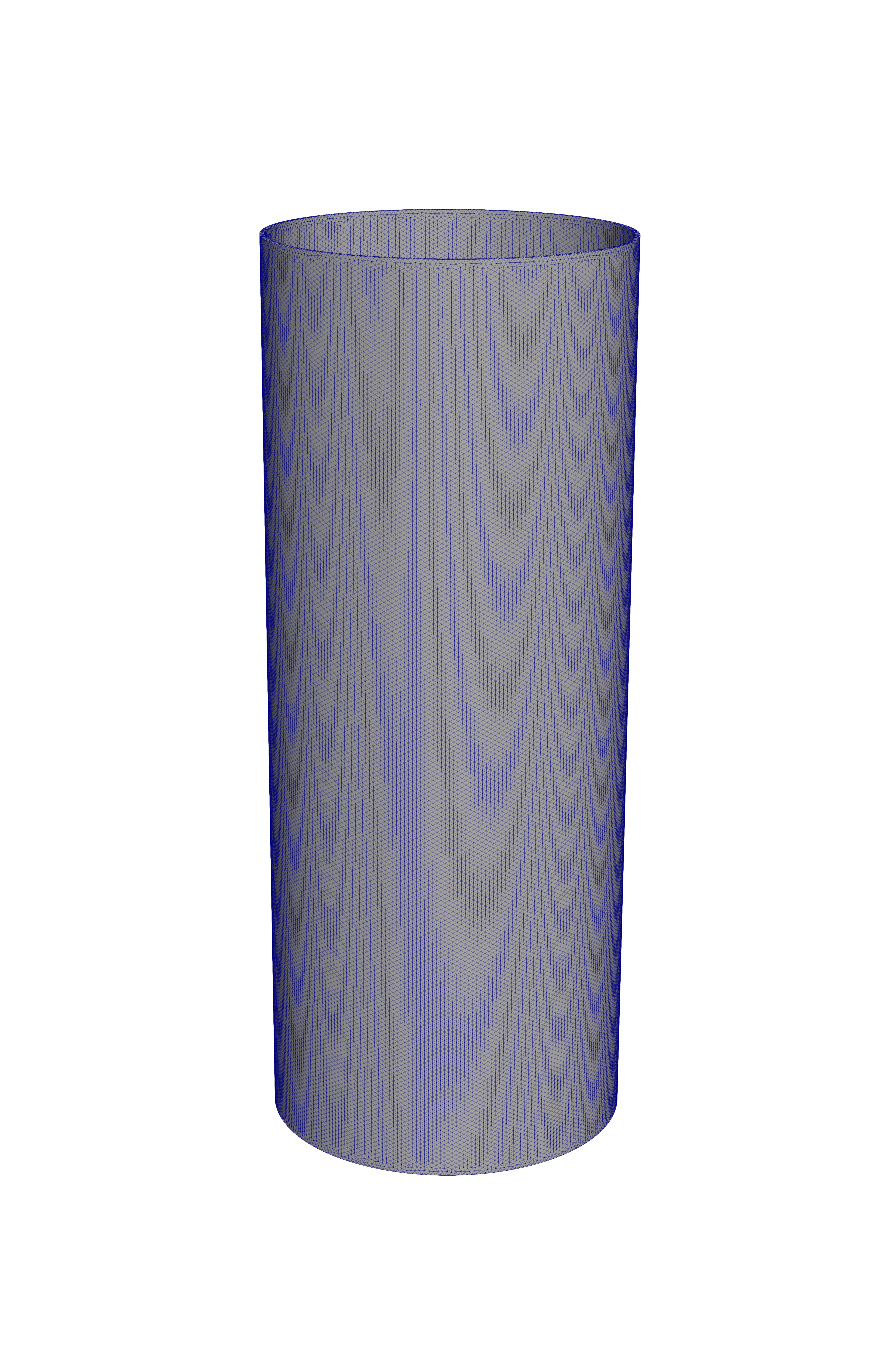}}
  \end{center}
  \caption{Impact of hyperelastic tube: Problem setup and mesh}
  \label{fig:Impactsetup}

  \begin{center}
    \subfigure[Potential, kinetic and total energies]
    {\includegraphics[width=.485\linewidth]{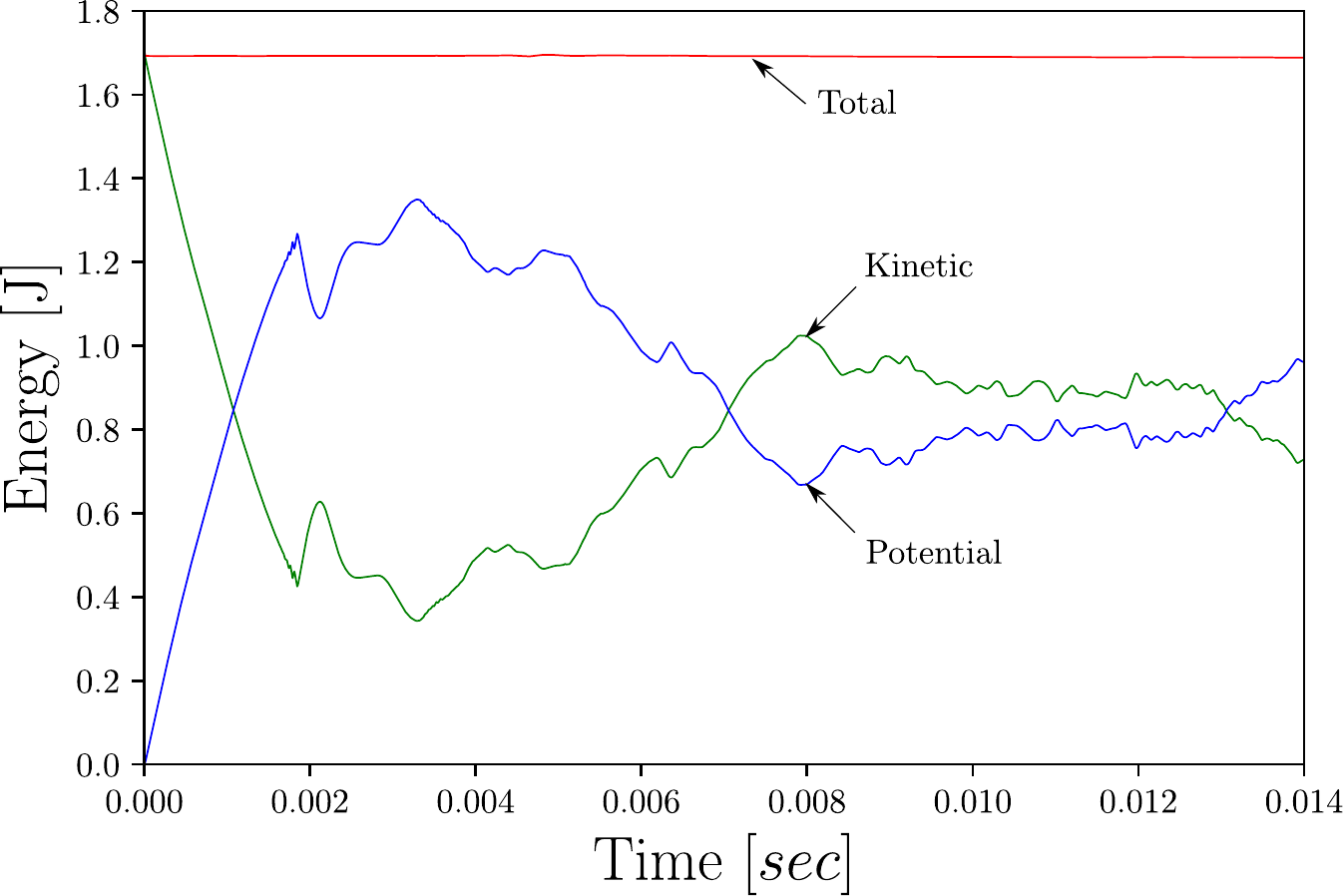}}
    \subfigure[Error and time-step size]
    {\includegraphics[width=.505\linewidth]{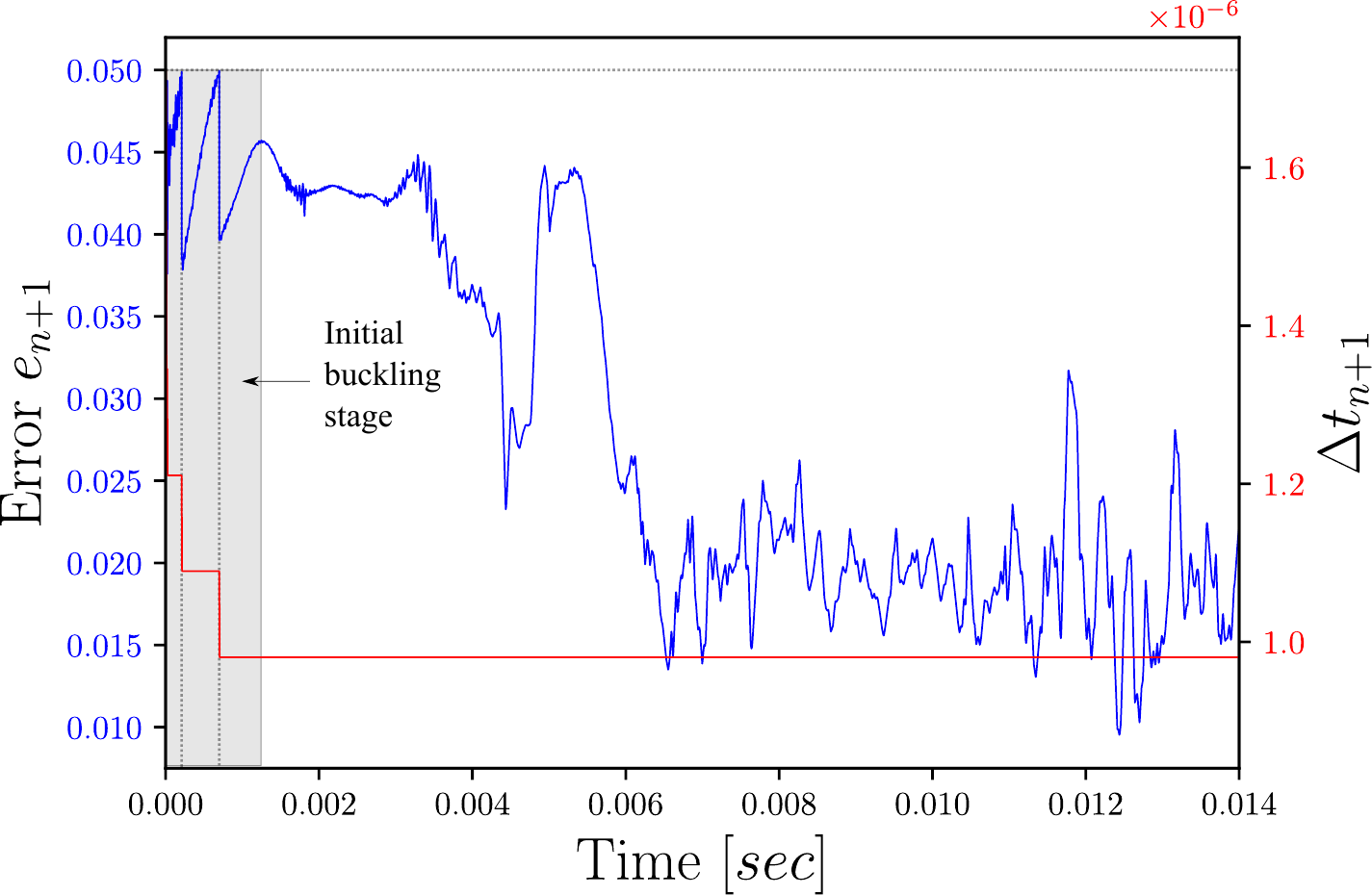}}
  \end{center}
  \caption{Energy stability, error and time-step size evolution vs time}
  \label{fig:ResImpactMesh1}
\end{figure}

\begin{figure} [h!]
  \begin{center}
    {\includegraphics[width=.75\linewidth]{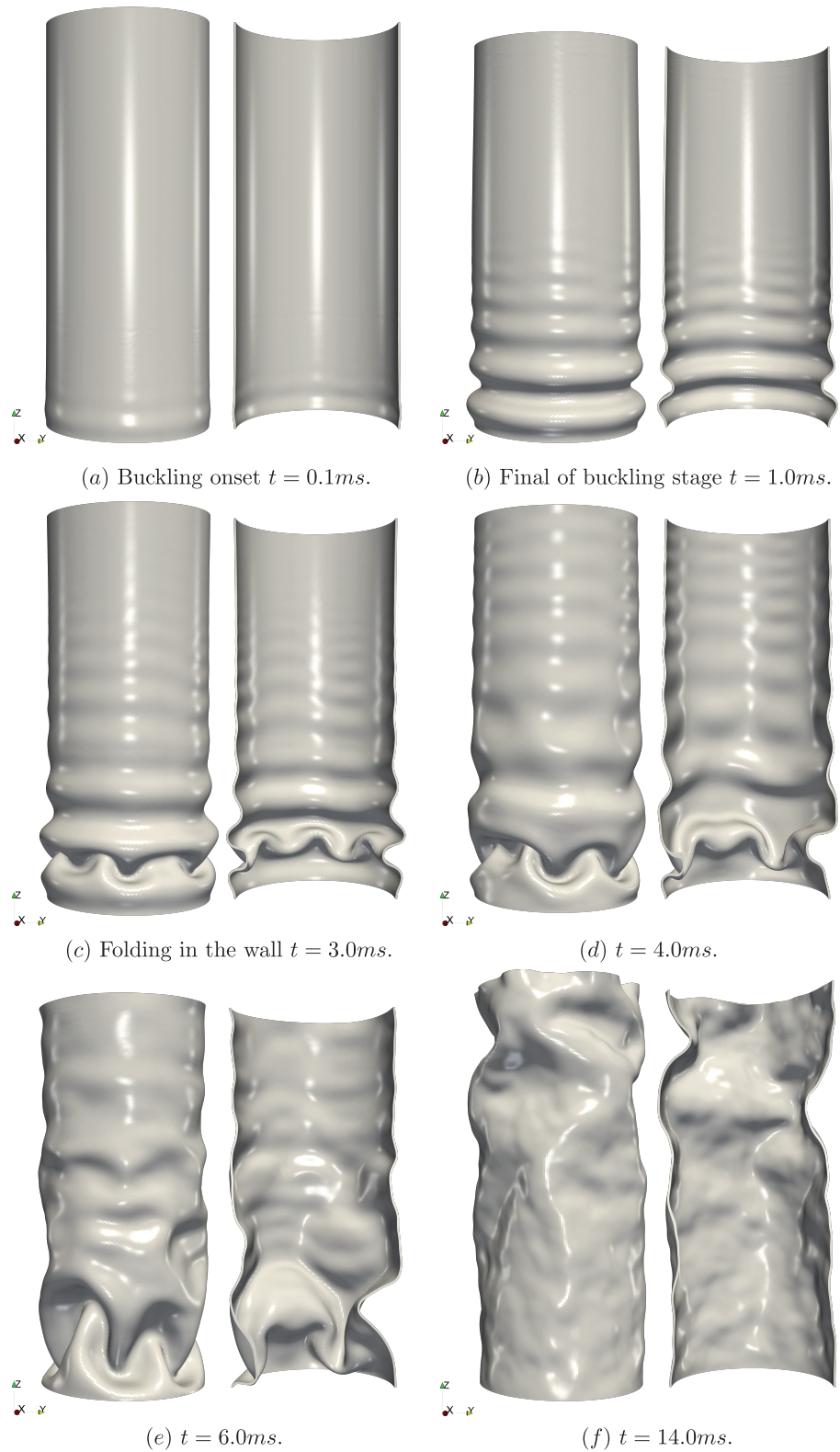}}
  \end{center}
  \caption{Deformed configuration snapshots}
  \label{fig:ResImpactMesh2}
\end{figure}

We analyse a thin-walled tube that impacts on a rigid surface. This test shows our method's ability to adapt the time-step size to robustly deal with simultaneous buckling processes.  Figure~\ref{fig:Impactsetup}~(a) sketches the problem, a cylindrical hyperelastic tube impacts a rigid surface with an initial velocity of $\dot{\boldsymbol{u}}_{0} = 20 \ \frac{\text{m}}{\text{s}}$, at impact we arrest the lower surface and its velocity becomes zero. The tube dimensions are $0.5 \ \text{m}$ length, $0.2 \ \text{m}$ diameter and $2.5 \ \text{mm}$ wall thickness,  material parameters are $E = 68 \ \text{MPa}$ and $\nu = 0.3$ for the Young modulus and Poisson's ratio,  respectively, while the initial density is $\rho_{0}= 1100  \ \frac{\text{kg}}{\text{m}^3}$.   Analogously to previous test, the parameter $\rho_{b}=1$ is assumed to avoid further numerical dissipation. We use an unstructured mesh with $\sim 350,000$ linear tetrahedra, with one element within the thickness; see Figure~\ref{fig:Impactsetup}~(b). 

Figure~\ref{fig:ResImpactMesh1}~(a) displays the total, kinetic, and potential energies; these plots show that the error control delivers energy stable solutions as the problem evolves Figure~\ref{fig:ResImpactMesh1}~(b) presents the error evolution and the time-step size versus time.; these plots attest that during the first $1.2 \ \text{ms}$ error accumulation occurs due to the initial buckling when the tube contacts the rigid surface. Meantime, the algorithm reduces the time-step size to ensure convergence, bounding the maximum eigenvalue of the amplification matrix.  

It is well-known that low-order elements produce poorer resolution results on stresses or any other gradient-related when compared with higher-order elements.  Thus,  denser meshes are usually recommended across the thickness if those quantities want to be assessed on specific spots.  Nonetheless,  the kinetic energy is showing indirectly that despite we used a coarser mesh,  overall the time-integrator is not losing information on gradient-dependent quantities.

Figure~\ref{fig:ResImpactMesh2}  shows representative simulation snapshots; each subfigure includes a front view and an internal view of the back wall. Figures~\ref{fig:ResImpactMesh2}~(a) and~(b) are snapshots during the initial evolution, where buckling generates a wave pattern. Figures~\ref{fig:ResImpactMesh2}~(c) and~(d) show the folding that occurs in the first buckling wave while beyond this point, the wall vibrates in an irregular pattern as Figures~\ref{fig:ResImpactMesh2}~(e) and~(f) show.  

\subsection{Impact of mass matrix formation: Space toss of a ruler}

\begin{figure}[h!]
  \begin{center}
    {\includegraphics[width=.8\linewidth]{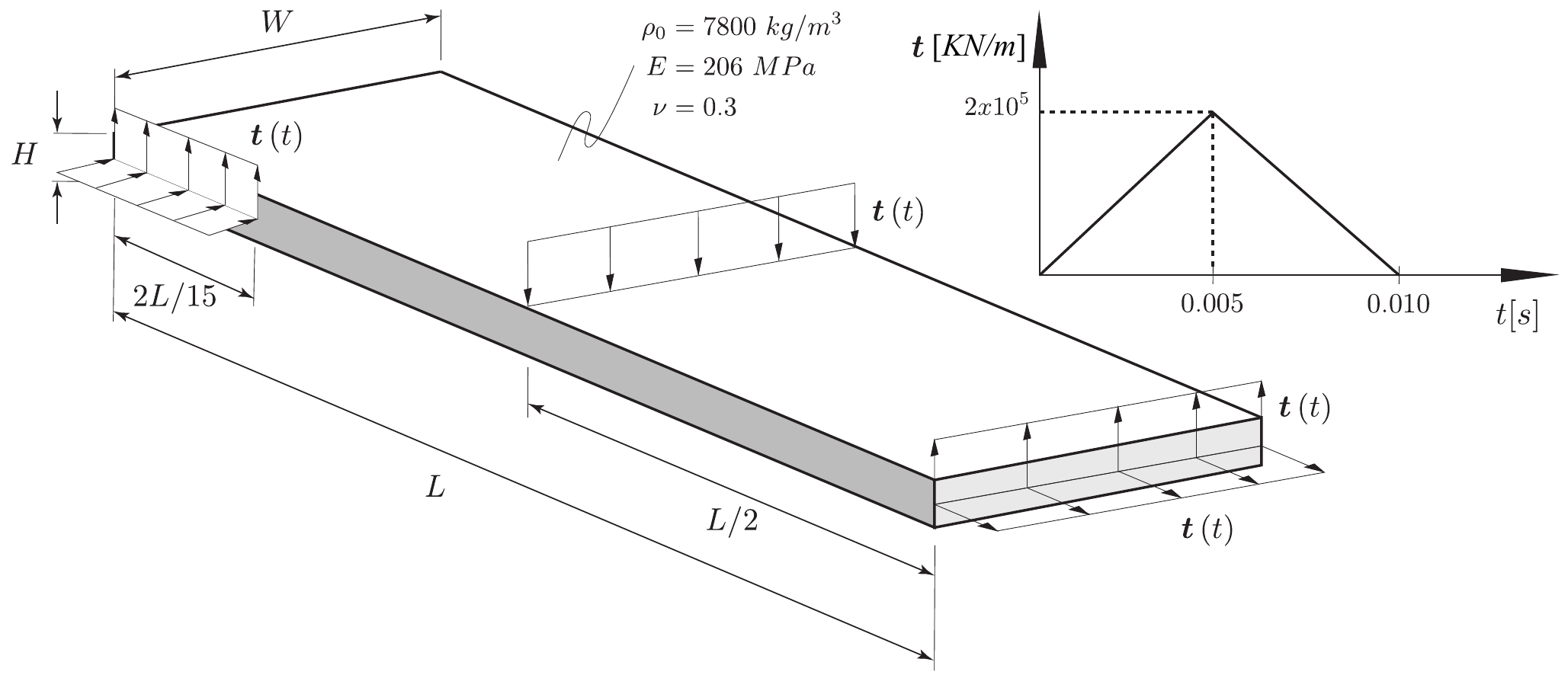}}
  \end{center}
  \caption{Tossed ruler: Problem setup}
  \label{fig:TossRuleSetup}
\end{figure}
This numerical experiment demonstrates the value of considering the implicit residual in the dynamical evolution of a nonlinear system.  For example, we turn a non-conservative explicit time integrator into a conservative solver, even when we use mass lumping. In~\cite{ KUHL1999343, Espath2015}, the authors discussed a similar example using different algorithms and time-integration schemes; they analyzed its conservation and stability properties. We investigate the space movement of a tossed ruler using our predictor/multicorrector method.  We simulate a geometrically nonlinear problem with evolving large displacements and finite strains.

We study two additional kinematic quantities in this example to better understand the conservation structure of the problem. These are the linear momentum, which we compute as
\begin{equation}
  \textbf{L} = \int_{\mathcal{B}_{0}} \rho_{0}   \dot{\boldsymbol{u}}
  \ d\mathcal{B}_{0}, 
  \label{eq:linearmomentum}
\end{equation}
and the angular momentum, which we compute as
\begin{equation}
  \textbf{J} = \int_{\mathcal{B}_{0}} \rho_{0}  \boldsymbol{x}  \times
  \dot{\boldsymbol{u}} \ d\mathcal{B}_{0}. 
  \label{eq:angularmomentum}
\end{equation}
\begin{figure}[h!]  
  \begin{center}
    {\includegraphics[width=.9\linewidth]{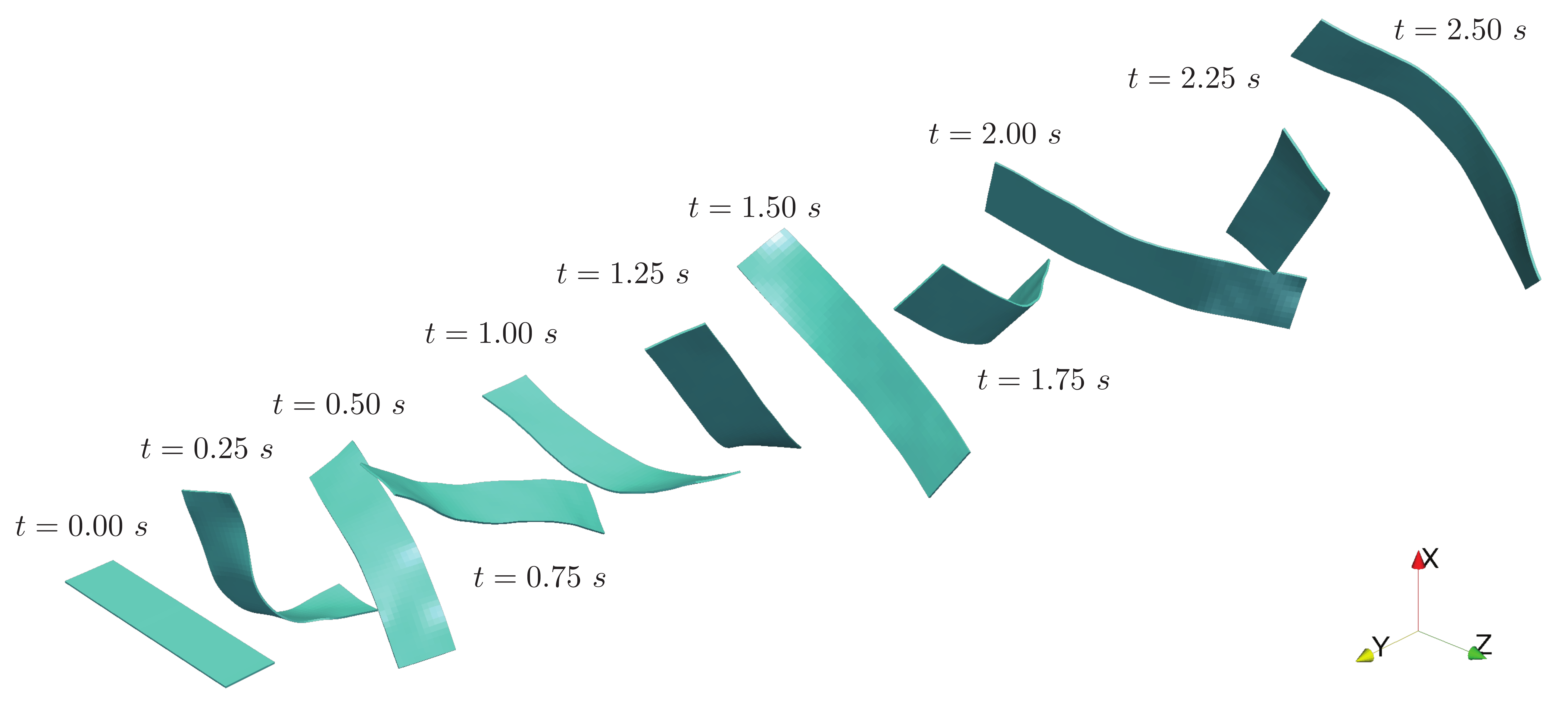}}
  \end{center}
  \caption{Space toss rule snapshots for consistent mass matrix.}
  \label{fig:TossRuleSnap}
\end{figure}

Figure~\ref{fig:TossRuleSetup} shows the problem setup of the numerical experiment, together with the material parameters and the applied force $\boldsymbol{t}$. The height $H$, width $W$ and length $L$ of the ruler are $0.002 \ \text{m}$, $0.06 \ \text{m}$ and $0.3 \ \text{m}$ respectively. We apply a set of linear tractions on ruler's edge and surface, using the same triangular variation in time,
\begin{equation}
  \boldsymbol{t} \left( t \right) =   \left\{
    \begin{aligned}
      &\frac{  t}{0.005} 200,000 \ \text{KN/m} ,
      &&   \text{if} \ \ t \leq 0.005 \ \text{s} , \\ 
      &\frac{ 0.010 - t}{0.005}    200,000 \ \text{KN/m} , \qquad
      & &  \text{if} \ \  0.005 \ \text{s} < t \leq 0.010 \ \text{s} ,\\
      &0.0 \ \text{KN/m},  &&  \text{if} \ \  t\geq  0.010 \ \text{s} ,
    \end{aligned}\right.
\end{equation}
where the maximum force is $200,000 \ \text{KN/m}$ at $0.005 \ \text{s}$. On the other hand, material parameters are $E = 206 \ \text{MPa}$ for the Young modulus, $\nu = 0.3 $ for the Poisson modulus and $\rho_{0} = 7800 \ {\frac{\text{kg}}{\text{m}^3}} $.  Since the purpose of this experiment is to prove that our setting is energy stable in long-time simulations, the spectral dissipation will be fixed on $\rho_{b} = 1.0$ in order to avoid dissipation of high frequencies.

\subsubsection{Consistent mass matrix}
\begin{figure} 
  \begin{center}
    \subfigure[Linear momentum.]
    {\includegraphics[height=5.0cm]{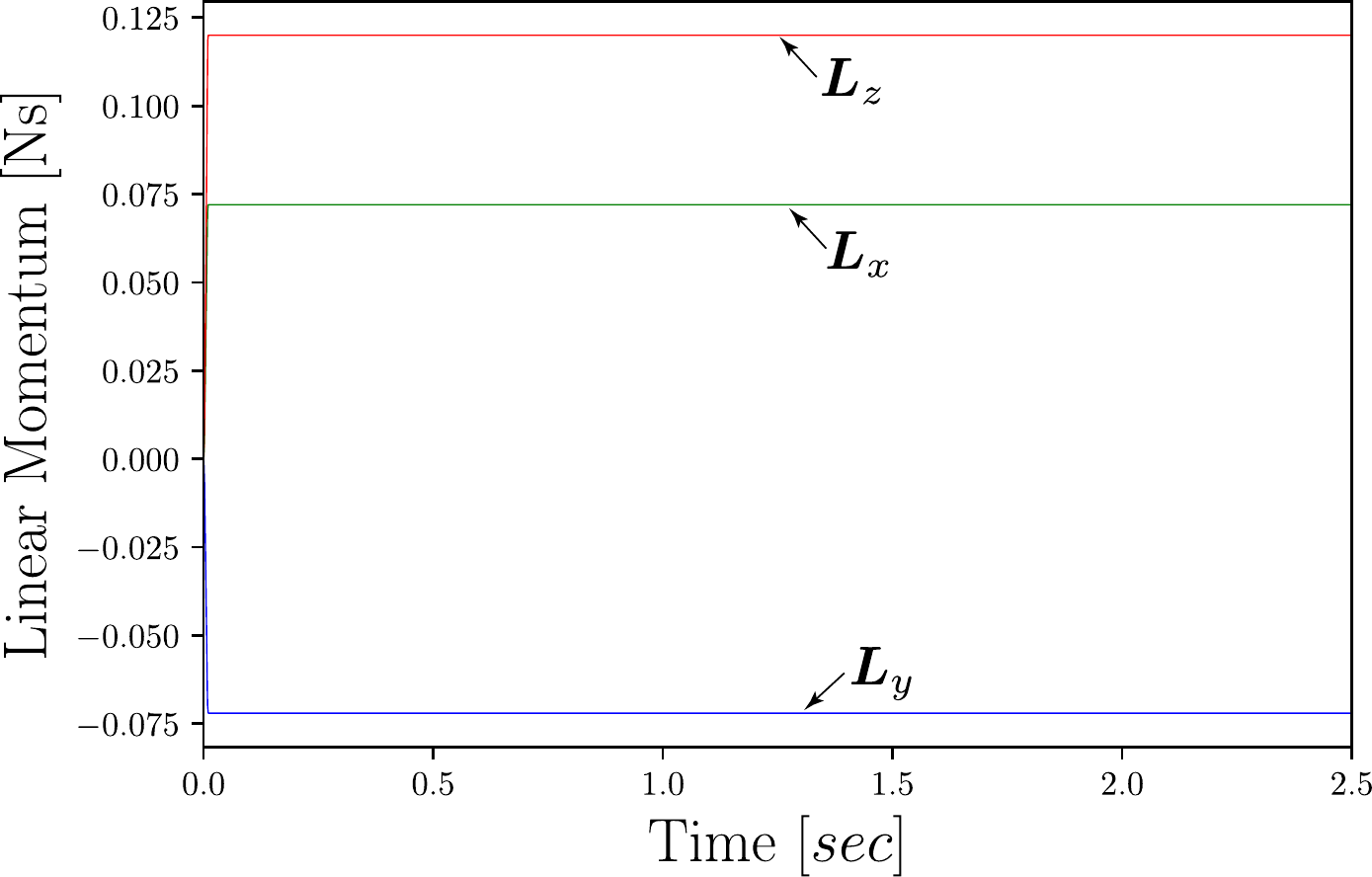}}
    \subfigure[Angular momentum.]
    {\includegraphics[height=5.0cm]{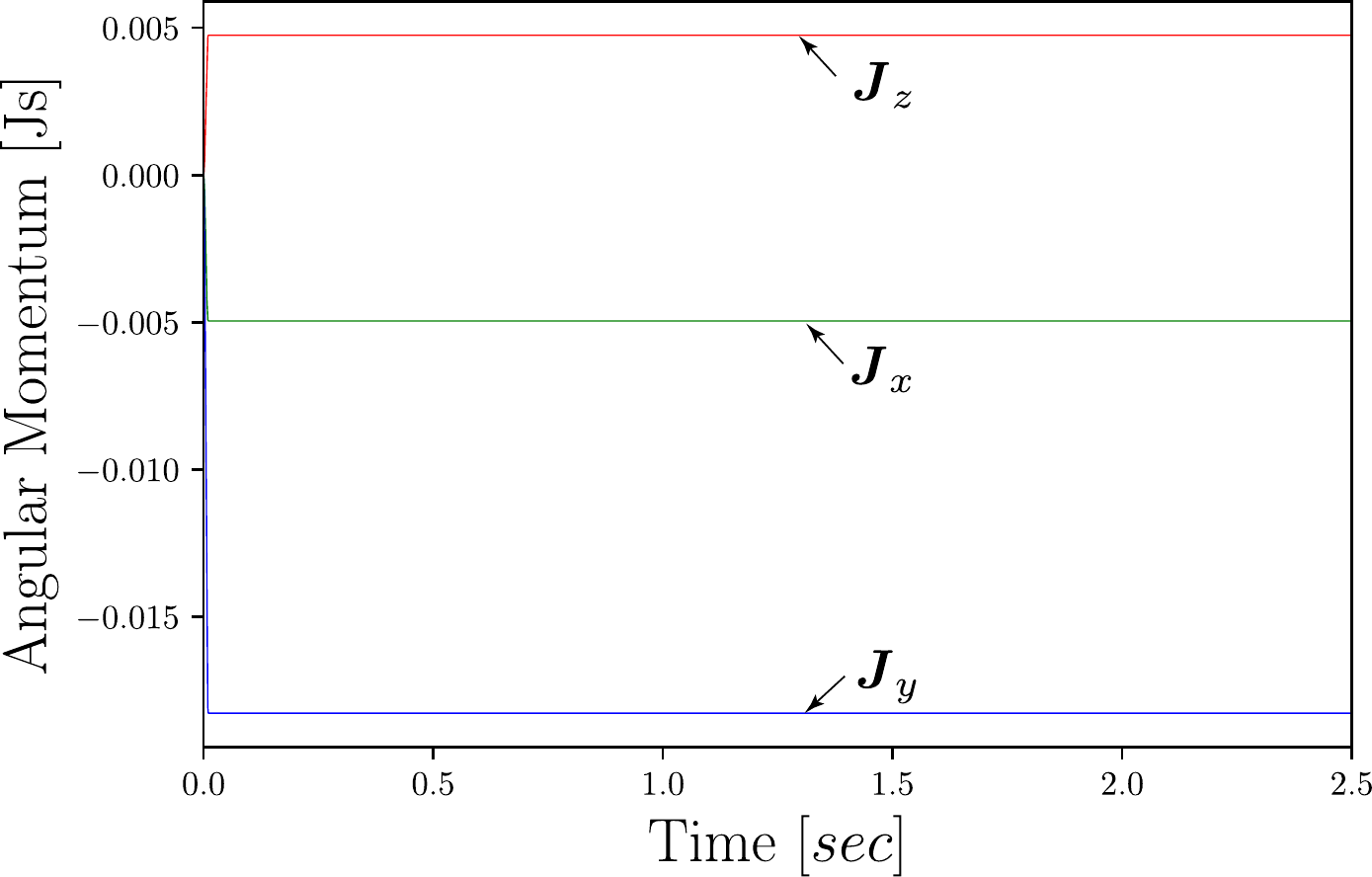}} \\
     \subfigure[Energy balance.]
     {\includegraphics[height=5.0cm]{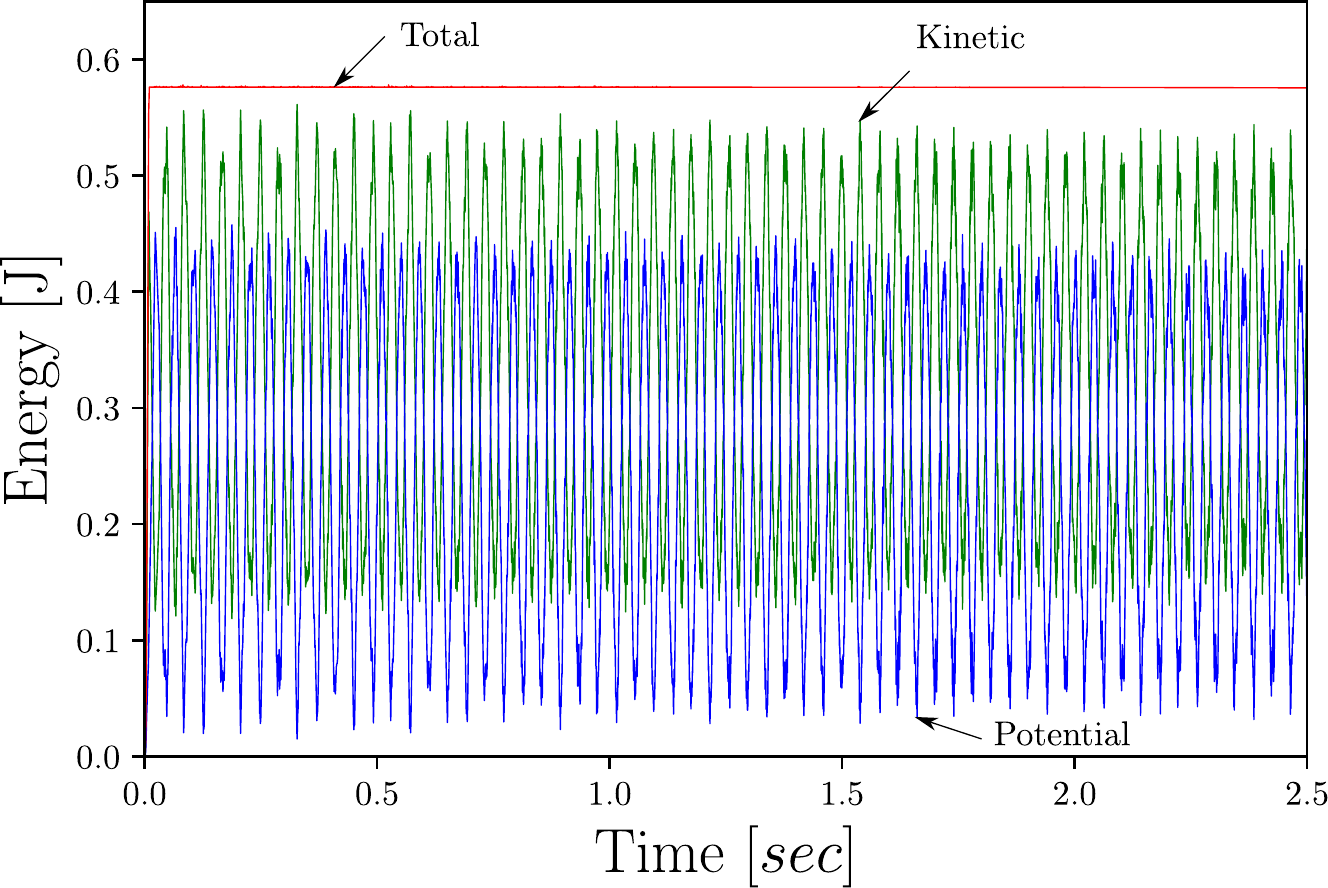}}
     \subfigure[Time-step adaptation.]
     {\includegraphics[height=5.0cm]{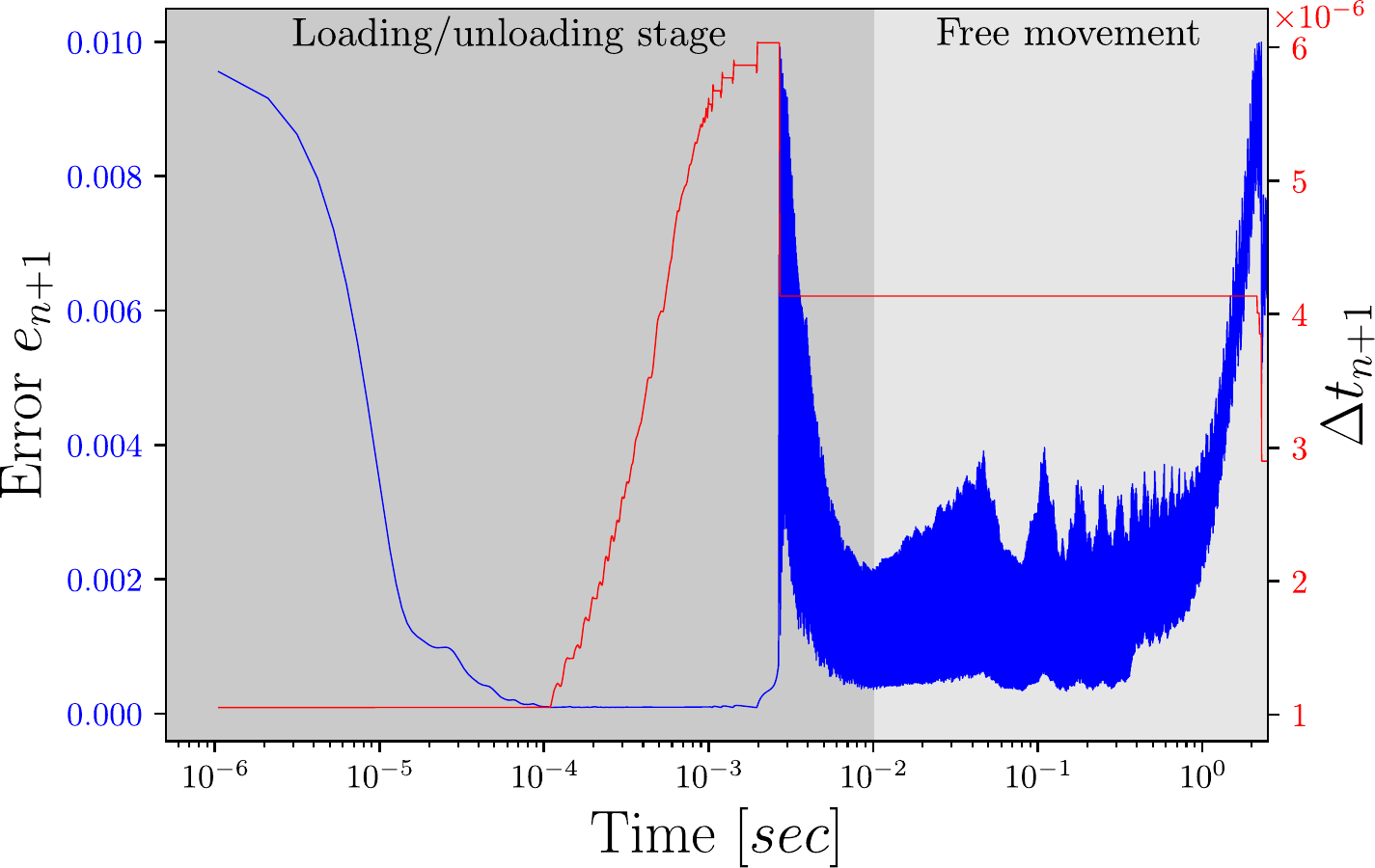}}
  \end{center}
  \caption{Tetrahedral mesh with consistent mass matrix (second term correction in~\eqref{eq:Galphau7})}
  \label{fig:ResConsistentTol5}
\end{figure}
We define a structured grid of points $50\text{x}10\text{x}1$ over the ruler's geometry and define a linear tetrahedral mesh with $3,000$ elements and $3,366$ degrees of freedom over this grid. The set the final time at $t_f = 2.5 \ \text{s}$ for all simulations to have a sufficiently long simulation to demonstrate the effects of energy dissipation on the results.  Figure~\ref{fig:TossRuleSnap} shows a set of solution snapshots of the problem, showing how the rule undergoes both large displacements and strains. The body translates and rotates over the three coordinate directions, while at the same time, it undergoes significant structural vibrations mainly in terms of bending moments.

Figure~\ref{fig:ResConsistentTol5} shows the linear momentum, angular momentum, energy balance and time-step adaptation for this first simulation.  With full quadrature of the mass matrix integration, the simulation conserves the total energy as well as the linear and angular momenta. Furthermore, Figure ~\ref{fig:ResConsistentTol5} (d) shows a typical time adaptivity we observe for the space tossed ruler problem. During the loading/unloading stage, the time-step size $\Delta t$ progressively grows from $10^{-6}$ to $6*10^{-6}$; once the ruler flies in free motion, the time-step size reduces to $\Delta t = 4.2e^{-6} \ \text{s}$ for the rest of the simulation.   Previous examples showed a time-step adaptive trend that progressively reduces its size due to the propagation of compressive waves,  increasing the maximum eigenvalue of the stiffness matrix due to the presence of a logarithm in the constitutive equation.  Tossed rule example also shows that,  when the problem allows it,  the time-marching method can increase the time step if truncation error $e_{n+1}$ drops below a minimum tolerance.
  
\subsubsection{Lumped mass matrix}
\begin{figure} [h!]
  \begin{center}
    \subfigure[Linear momentum.]
    {\includegraphics[height=5.0cm]{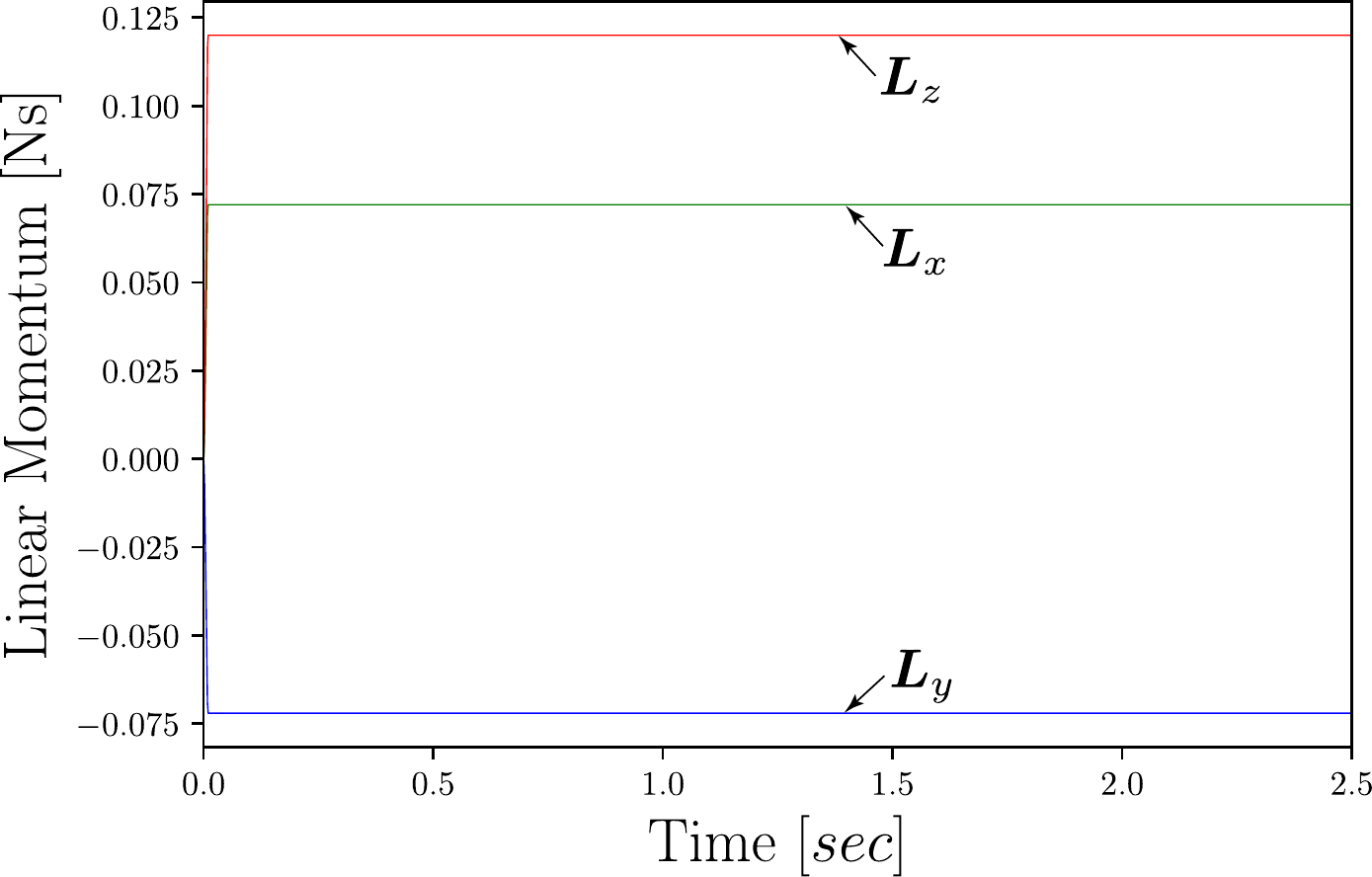}} 
  \subfigure[Angular momentum.]
  {\includegraphics[height=5.0cm]{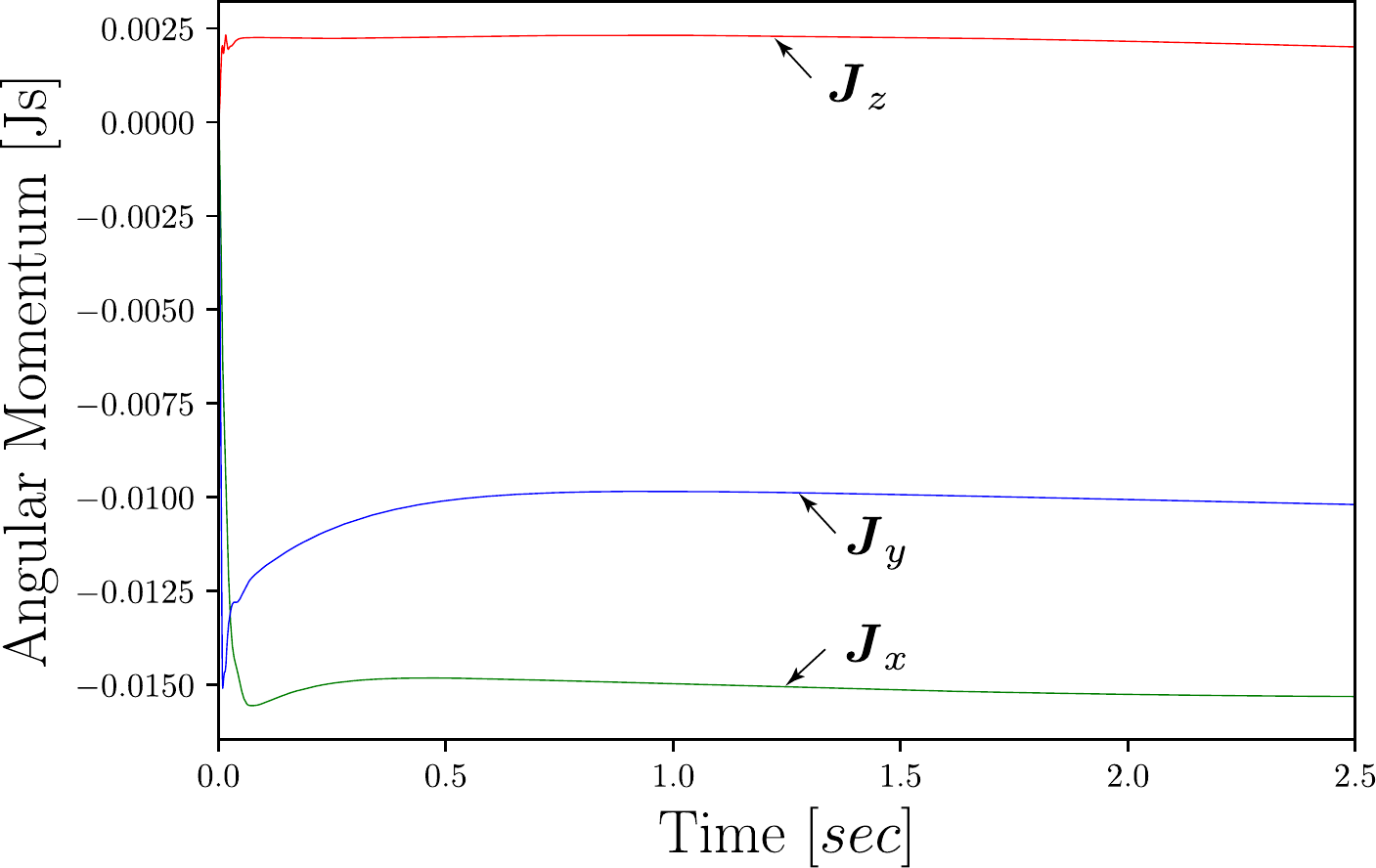}} \\
    \subfigure[Energy balance.]
    {\includegraphics[height=5.0cm]{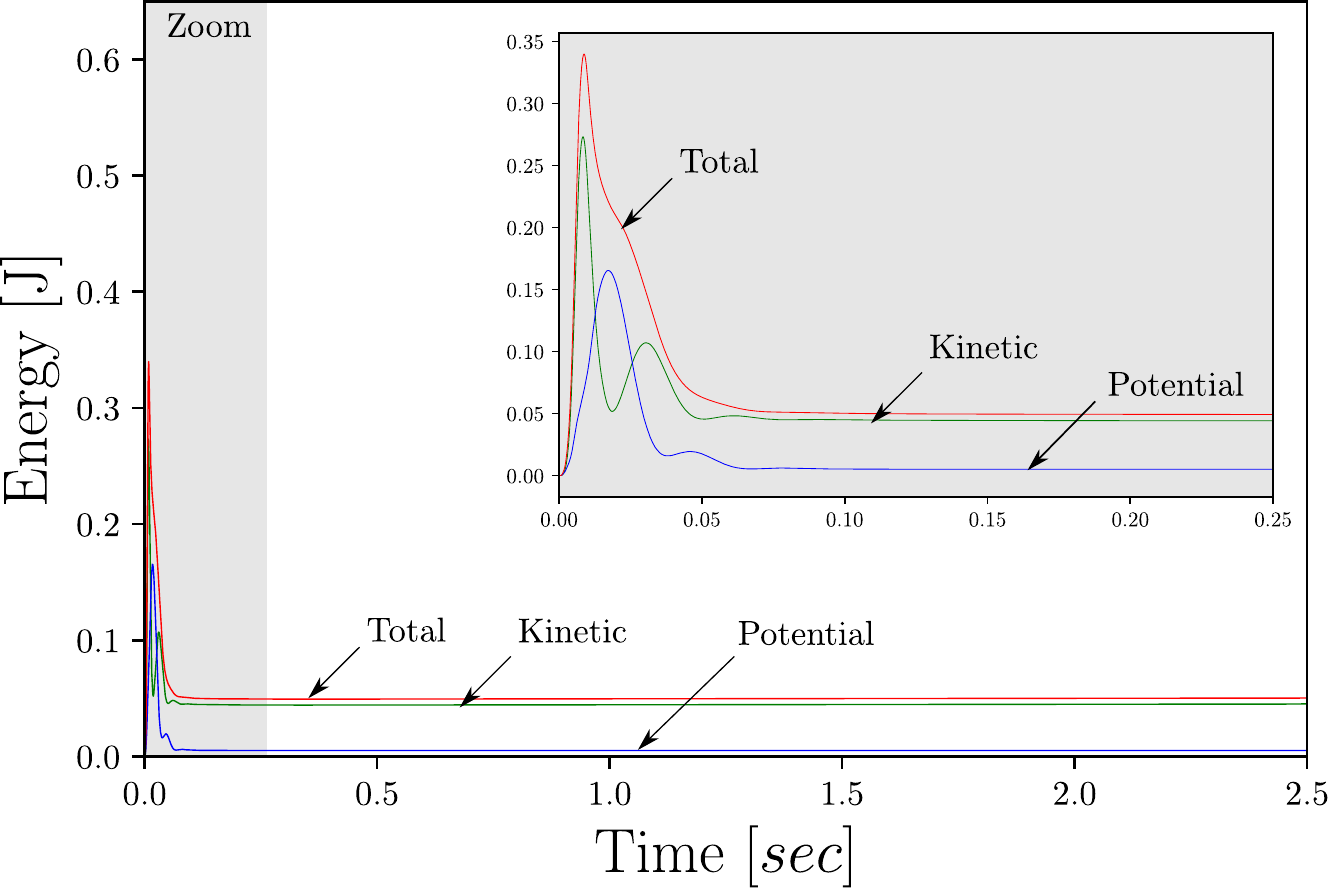}}
  \end{center}
  \caption{Tetrahedral mesh with lumped mass matrix (second term correction in~\eqref{eq:Galphau7})}
  \label{fig:ReslumpedTol2}
\end{figure}
Until now, all simulations use consistent integration (Gauss quadrature) to form the mass matrix, while we only use the first two terms of the update formula~\eqref{eq:Galphau7}, reducing our method to a predictor-corrector integrator. Effectively, this strategy is standard in explicit multistep methods. Nevertheless, this kind of simulation is common practice to use lumping by a Gauss-Lobatto quadrature, collocating the degrees of freedom with the quadrature points. This methodology results in a lumped (diagonal) mass matrix.  We study the effect of lumping on the simulation results and address its impact.

\begin{figure}  
  \begin{center}
    {\includegraphics[width=.9\linewidth]{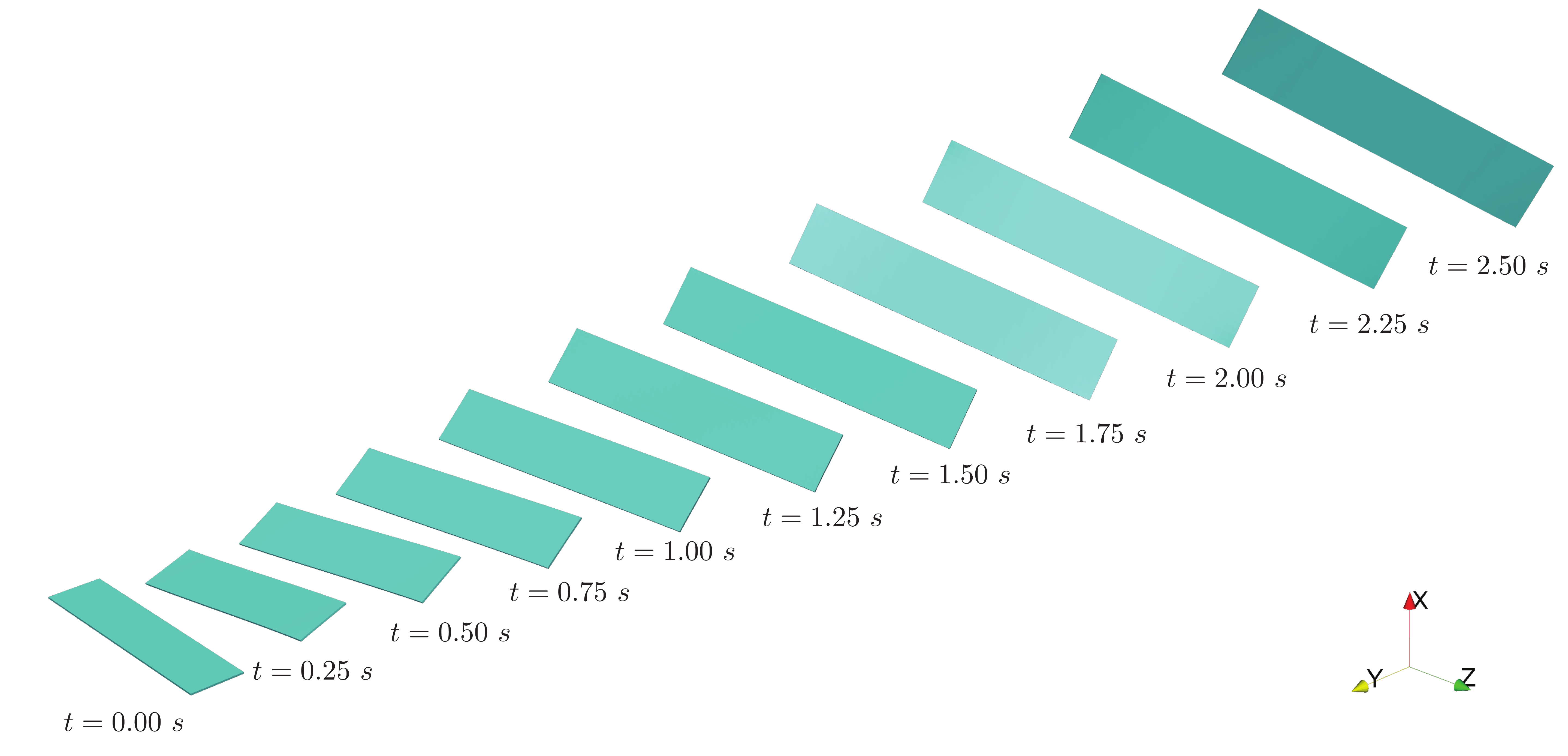}}
  \end{center}
  \caption{Tossed ruler: solution snapshots for lumped mass matrix associated to Figure~\ref{fig:ReslumpedTol2}}
  \label{fig:TossRuleSnapLumped}
\end{figure}

Figure~\ref{fig:ReslumpedTol2} demonstrates how just using the first two terms of the update formula~\eqref{eq:Galphau7} with a lumped mass matrix loses the conservation properties on the time marching process. In particular, Figure~\ref{fig:ReslumpedTol2}~(c) shows a significant loss of energy during the first $0.05 \ \text{s}$ of the simulation. Similarly, Figure~\ref{fig:ReslumpedTol2}~(b) shows that the simulation does not conserve angular momentum. Effectively, mass lumping with explicit time marching conserves the linear momentum as Figure~\ref{fig:ReslumpedTol2}~(a) shows. 

\begin{figure} 
  \begin{center}
  \subfigure[Truncation error, $ TOL =1e^{-10}$]
  {\includegraphics[width=.48\linewidth]{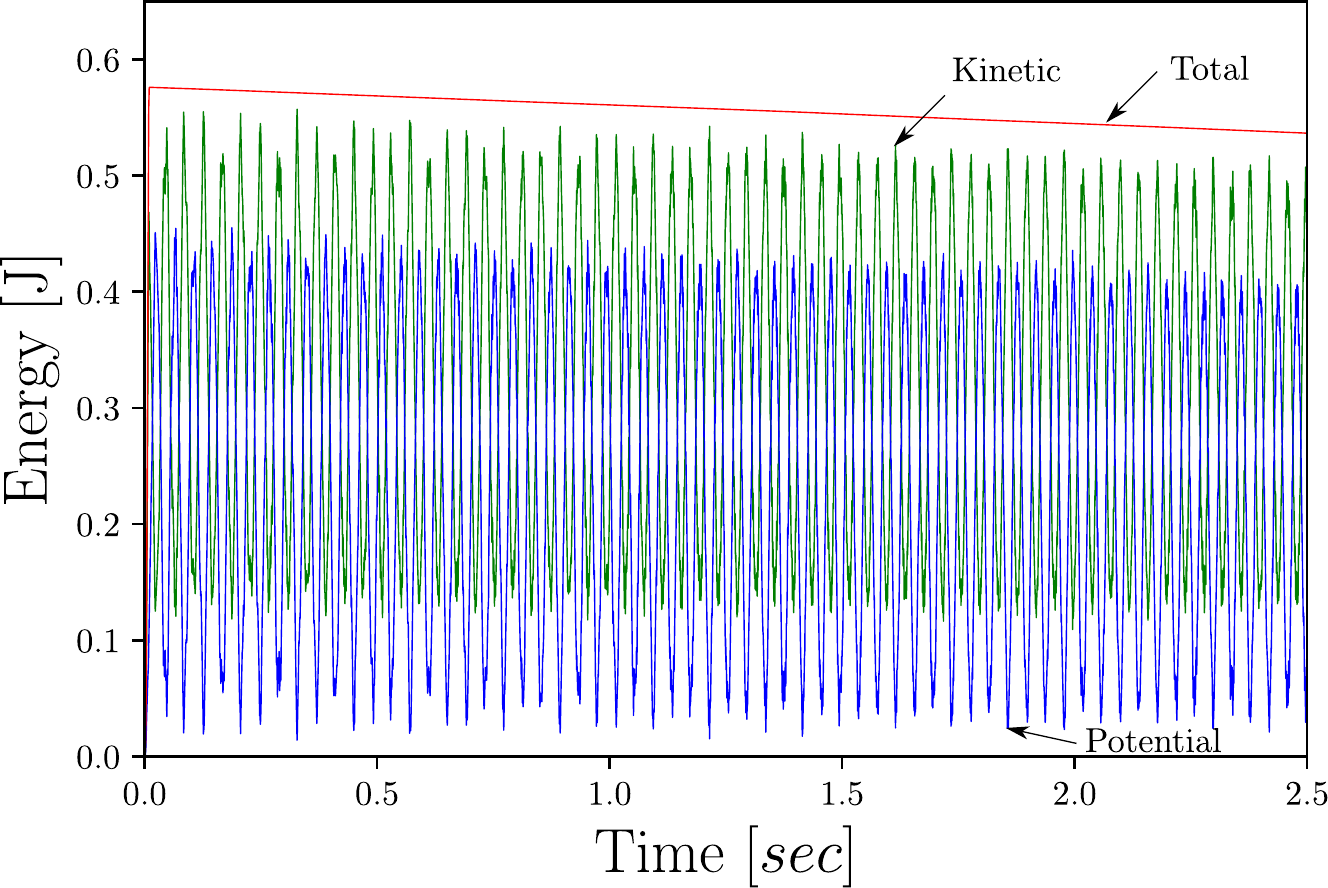}}
  \subfigure[Truncation error, $ TOL= 1e^{-14}$]
  {\includegraphics[width=.48\linewidth]{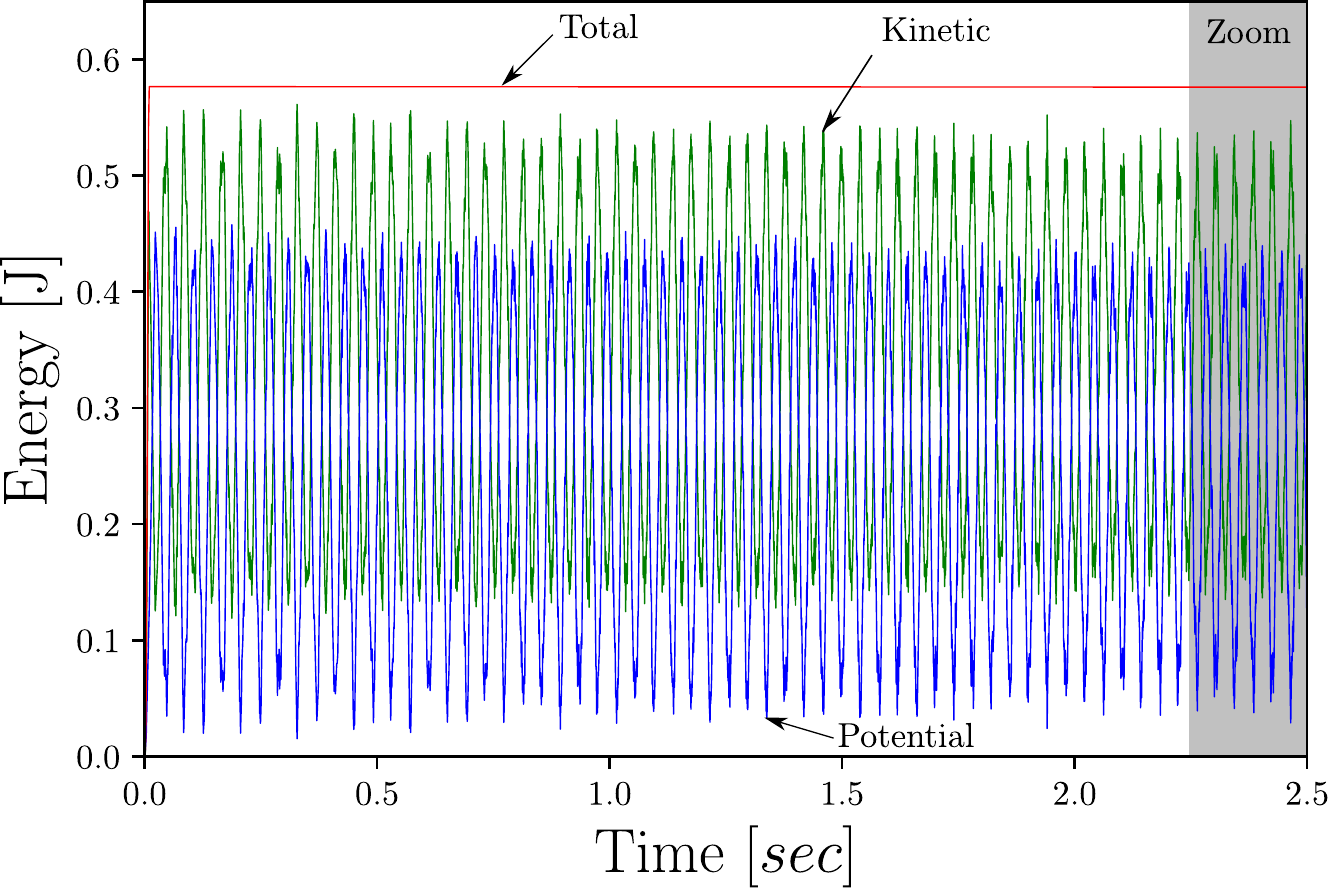}}
  \subfigure[Lumped vs consistent mass, zoom of grey zone in~(b)]
  {\includegraphics[width=.48\linewidth]{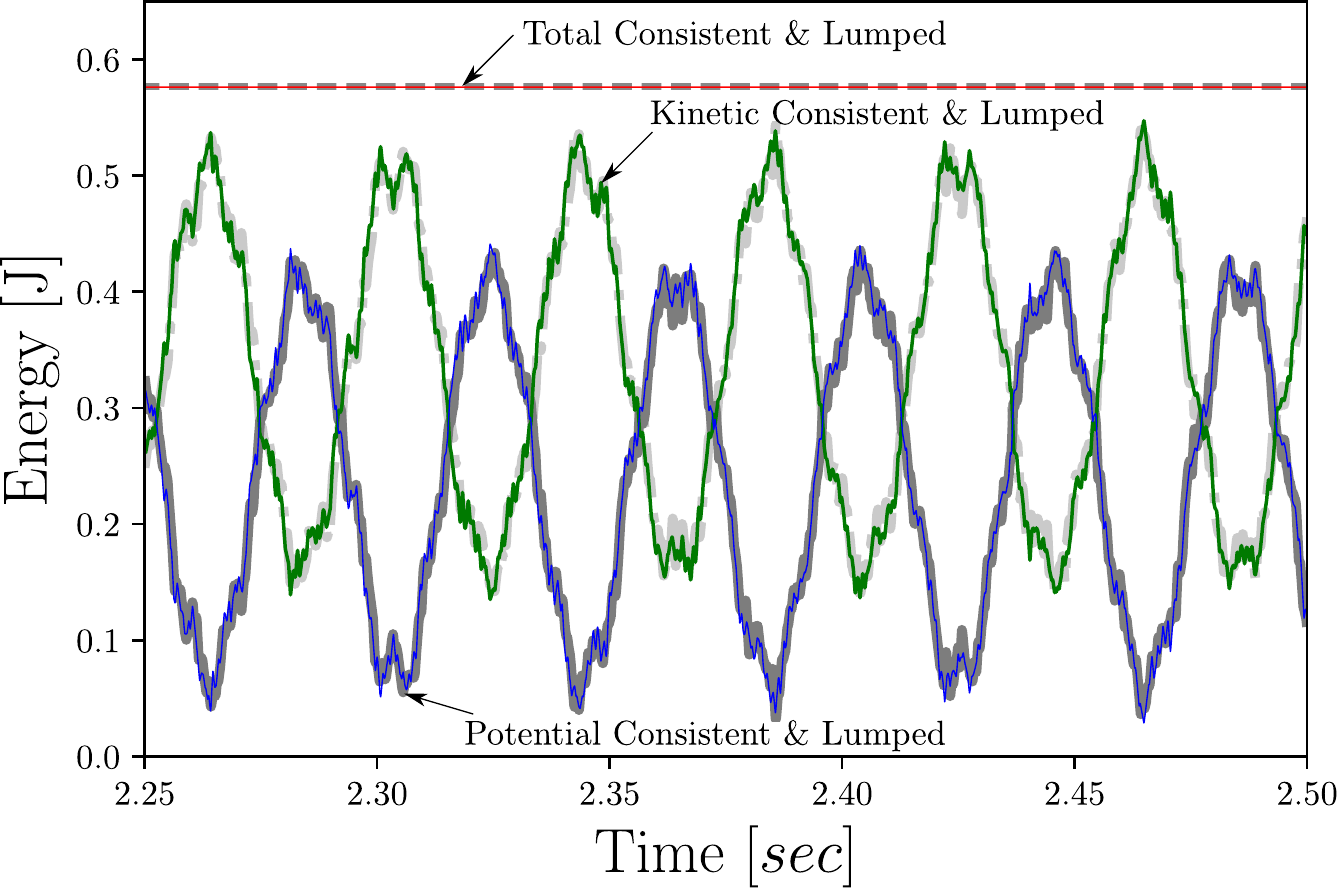}}   
  \subfigure[Angular momentum:  $ TOL = 1e^{-10}\ \& \ 1e^{-14}$]
  {\includegraphics[width=.48\linewidth]{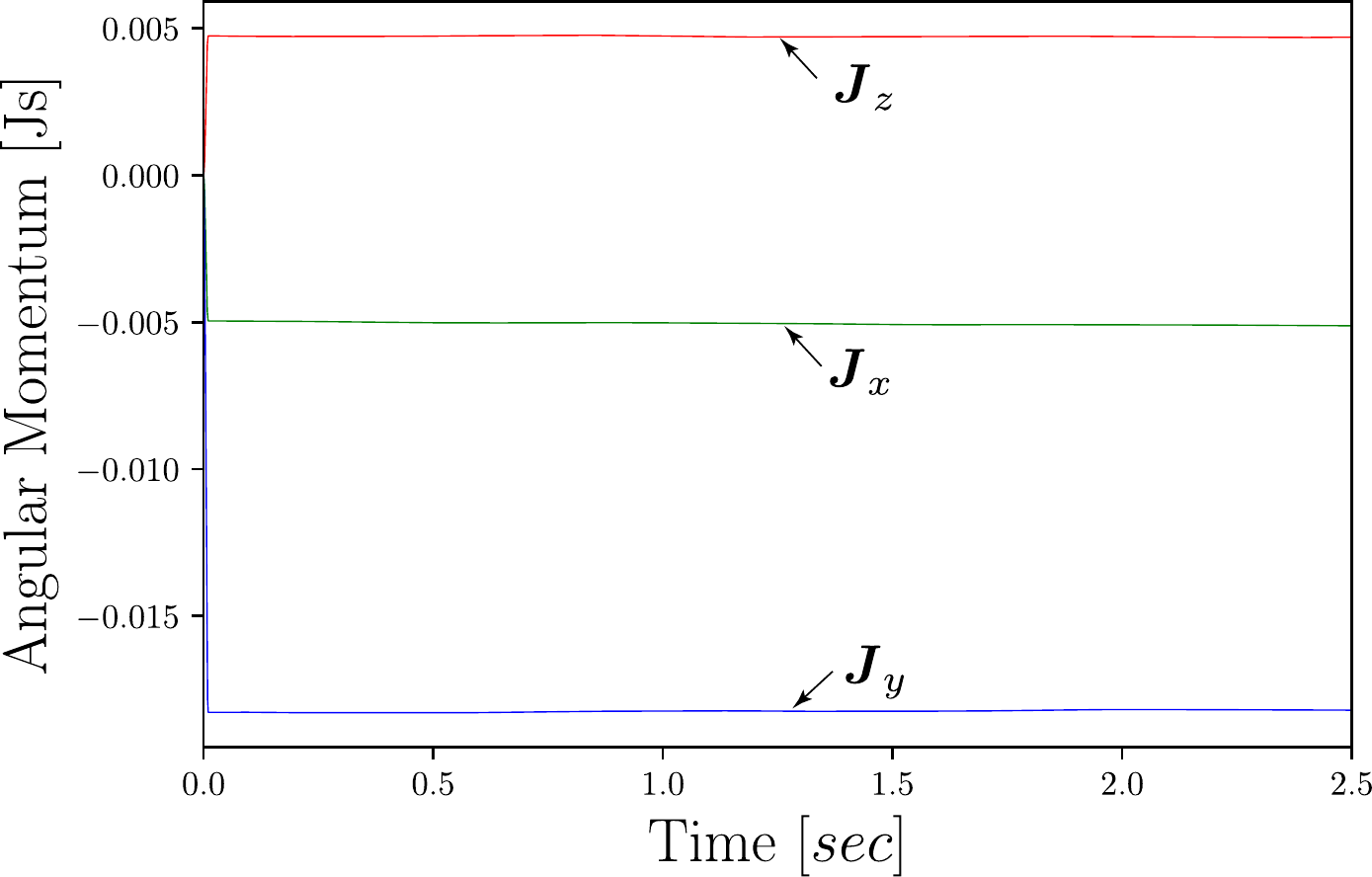}} 
  \end{center}
  \caption{Tetrahedral mesh with lumped mass matrix (truncation control in~\eqref{eq:Galphau7})}
  \label{fig:ResLumpedTol10}
\end{figure}

Figure~\ref{fig:TossRuleSnapLumped} shows a set of snapshots of the simulation computed with the lumped matrix, using just the first two terms of~\eqref{eq:Galphau7}. Under these conditions, the motion of the ruler becomes almost pure translation combined with a slight rotation over the z-axis. At the same time, the deformation is negligible due to the spurious numerical damping.

The mass matrix resulting from using a Gauss-Lobatto quadrature is diagonal, but the resulting simulation loses conservation when using standard update processes; it becomes unphysical. We circumvent this limitation by considering more terms in the summation we use in the update formula~\eqref{eq:Galphau7} for our predictor-multicorrector method. Figure~\ref{fig:ResLumpedTol10} shows the energy balance and the angular momentum we obtain when we control truncation error; that is, we iterate to enforce $TOL = \| \delta \boldsymbol{u}^{k}_{n+1} \| \leq 1e^{-10}$. Figure~\ref{fig:ResLumpedTol10}~(a) shows a improvement in the energy conservation (compare with Figure~\ref{fig:ReslumpedTol2}~(c)); nonetheless, the simulation still dissipates energy as the simulation progresses. On average, the simulation requires about eight iterations to reach the tolerance. Now, to emphasize the importance of converging the nonlinear implicit residual, we use truncation tolerance to  $\| \delta \boldsymbol{u}^{k}_{n+1} \| \leq 1e^{-14}$. The stringent convergence requirement delivers a conservative simulating, as Figure~\ref{fig:ResLumpedTol10}~(b) shows. Effectively, mass lumping recovers the response of the consistent mass matrix as Figure~\ref{fig:ResConsistentTol5}~(c) shows the overlap in the energy patterns displayed by both methods during the last few cycles of the evolution. The dashed grey lines correspond to the consistent mass, while the colored correspond to the lumped mass. In this case, the average iteration number in the update formula is twenty. Lastly, both tolerances conserve the angular momentum, see Figure~\ref{fig:ResLumpedTol10} (d). Thus, our predictor-multicorrector method is conservative for consistent and lumped matrices for nonlinear problems with large displacements and deformations for long integration times.  

%

 \section{Conclusions}
\label{section:Conclusions}

We introduce a new class of predictor/multicorrector time-marching techniques for explicit generalized-$\alpha$ methods. Our proposal has a cheap and robust error estimator that can guide adaptability even in the presence of sudden shocks to the system. Another key advantage of our predictor/multicorrector method is that our variable updates are explicit from $n$ to $n+1$ while evaluating the residual at $n+1$; we only factor the mass matrix to compute the explicit update, while we can control the residual at time $n+1$ as in implicit methods. We estimate the method's error by comparing the truncation terms of the update formula. This time-error estimator provides a robust quantity to guide the time-step size adaptivity, circumventing the lack of convergence of the nonlinear residual that characterizes explicit methods. 

Our explicit generalized-$\alpha$ method controls the high-frequency numerical dissipation. We analyze the method's stability regimes and determine the corresponding parameters to guarantee second-order accuracy in time and to maximize stability.  We demonstrate the method's performance (accuracy and efficiency) on various numerical examples. The first two confirm that our method delivers the optimal second-order convergence for linear and nonlinear hyperbolic equations. Next, we study large-deformation problems subject to sudden shocks that induce instantaneous energy releases. We study a twisted bar subject to gravitational loading that we release suddenly and the impact of a cylindrical tube onto a rigid wall. These examples show our time adaptability's ability to describe the dynamics of rubber-like hyperelastic material subject to sudden releases of potential or kinetic energies. The final example is a classical benchmark used in the community to analyze the momentum and energy conservation. The \textit{tossed ruler} problem shows that our method is energy stable when using consistent and lumping mass matrices to compute the variable update. We show that the predictor/multicorrector approach described is a suitable tool to correct the energy unbalance produced when Gauss-Lobatto quadrature, which is commonly used to obtain diagonal mass matrices in dynamic problems. 

In conclusion, our method is a computationally efficient and robust predictor/multicorrector method for real-life applications, where convergence is always guaranteed regardless of the problems' dynamics while accounting for shocks within the solution as well as from external loadings,  providing energy and momentum stable solutions. This is a key advantage when it is compared with other explicit time-marching methods, where the error is not a built-in quantity while it has to be computed using methods with different convergence order,  at the same time they dissipate energy and momentum when lumped mass matrices are used.


%



\bibliographystyle{unsrt}
\bibliography{mybibfile}

\begin{thebibliography}{10}

\bibitem{hughes2012finite}
T.~J.R. Hughes.
\newblock {\em The finite element method: linear static and dynamic finite
  element analysis}.
\newblock Courier Corporation, 2012.

\bibitem{Newmark1959}
N.~M. Newmark.
\newblock A method of computation for structural dynamics.
\newblock {\em Journal of the Engineering Mechanics Division}, 85(3):67--94,
  1959.

\bibitem{wilson1968}
E.L. Wilson.
\newblock {\em A computer program for the dynamic stress analysis of
  underground structures}.
\newblock SESM Report No. 68-1, Division of Structural Engineering and
  Structural Mechanics, University of California, Berkeley, CA,, 1968.

\bibitem{Hilber1977}
H.~M. Hilber, T.~J.~R. Hughes, and R.~L. Taylor.
\newblock Improved numerical dissipation for time integration algorithms in
  structural dynamics.
\newblock {\em Earthquake Engineering \& Structural Dynamics}, 5(3):283--292,
  1977.

\bibitem{wood1980}
W.~L. Wood, M.~Bossak, and O.~C. Zienkiewicz.
\newblock An $\alpha$ modification of {N}ewmark's method.
\newblock {\em International Journal for Numerical Methods in Engineering},
  15(10):1562--1566, 1980.

\bibitem{Bazzi1982}
G.~Bazzi and E.~Anderheggen.
\newblock The $\rho$-family of algorithms for time-step integration with
  improved numerical dissipation.
\newblock {\em Earthquake Engineering \& Structural Dynamics}, 10(4):537--550,
  1982.

\bibitem{HOFF1988367}
C.~Hoff and P.J. Pahl.
\newblock Development of an implicit method with numerical dissipation from a
  generalized single-step algorithm for structural dynamics.
\newblock {\em Computer Methods in Applied Mechanics and Engineering},
  67(3):367 -- 385, 1988.

\bibitem{HOFF198887}
C.~Hoff and P.J. Pahl.
\newblock Practical performance of the $\theta$1-method and comparison with
  other dissipative algorithms in structural dynamics.
\newblock {\em Computer Methods in Applied Mechanics and Engineering}, 67(1):87
  -- 110, 1988.

\bibitem{HOFF198987}
C.~Hoff, T.J.R. Hughes, G.~Hulbert, and P.J. Pahl.
\newblock Extended comparison of the {H}ilber-{H}ughes-{T}aylor $\alpha$-method
  and the $\theta$1-method.
\newblock {\em Computer Methods in Applied Mechanics and Engineering}, 76(1):87
  -- 93, 1989.

\bibitem{chung1993time}
J.~Chung and G.M. Hulbert.
\newblock A time integration algorithm for structural dynamics with improved
  numerical dissipation: the generalized-$\alpha$ method.
\newblock {\em Journal of Applied Mechanics}, 60(2):371--375, 1993.

\bibitem{Bazilevs2006}
Y.~Bazilevs, V.M. Calo, Zhang Y., and T.J.R. Hughes.
\newblock Isogeometric fluid–structure interaction analysis with applications
  to arterial blood flow.
\newblock {\em Computational Mechanics}, 38:310--332, 2006.

\bibitem{Bazilevs2008}
Y.~Bazilevs, V.M. Calo, T.J.R. Hughes, and Zhang Y.
\newblock Isogeometric fluid-structure interaction: theory, algorithms, and
  computations.
\newblock {\em Computational Mechanics}, 43:3--37, 2008.

\bibitem{NOELS2005358}
L.~Noels, L.~Stainier, and J.-P. Ponthot.
\newblock Simulation of complex impact problems with implicit time algorithms:
  Application to a turbo-engine blade loss problem.
\newblock {\em International Journal of Impact Engineering}, 32(1):358--386,
  2005.
\newblock Fifth International Symposium on Impact Engineering.

\bibitem{ROSSI2016208}
S.~Rossi, N.~Abboud, and G.~Scovazzi.
\newblock Implicit finite incompressible elastodynamics with linear finite
  elements: A stabilized method in rate form.
\newblock {\em Computer Methods in Applied Mechanics and Engineering},
  311:208--249, 2016.

\bibitem{ScovazziIJMNE2016}
G.~Scovazzi, B.~Carnes, X.~Zeng, and S.~Rossi.
\newblock A simple, stable, and accurate linear tetrahedral finite element for
  transient, nearly, and fully incompressible solid dynamics: a dynamic
  variational multiscale approach.
\newblock {\em International Journal for Numerical Methods in Engineering},
  106(10):799--839, 2016.

\bibitem{BONET2015689}
J.~Bonet, A.~J. Gil, C.~H. Lee, M.~Aguirre, and R.~Ortigosa.
\newblock A first order hyperbolic framework for large strain computational
  solid dynamics. {P}art {I}: Total lagrangian isothermal elasticity.
\newblock {\em Computer Methods in Applied Mechanics and Engineering},
  283:689--732, 2015.

\bibitem{GIL2016146}
A.~J. Gil, C.~H. Lee, J.~Bonet, and R.~Ortigosa.
\newblock A first order hyperbolic framework for large strain computational
  solid dynamics. {P}art {II}: Total lagrangian compressible, nearly
  incompressible and truly incompressible elasticity.
\newblock {\em Computer Methods in Applied Mechanics and Engineering},
  300:146--181, 2016.

\bibitem{LAVRENCIC2020107075}
M.~Lavrenčič and B.~Brank.
\newblock Comparison of numerically dissipative schemes for structural
  dynamics: Generalized-alpha versus energy-decaying methods.
\newblock {\em Thin-Walled Structures}, 157:107075, 2020.

\bibitem{BEHNOUDFAR2020100021}
P.~Behnoudfar, V.~M. Calo, Q.~Deng, and P.~D. Minev.
\newblock A variationally separable splitting for the generalized-$\alpha$
  method for parabolic equations.
\newblock {\em International Journal for Numerical Methods in Engineering},
  121(5):828--841, 2020.

\bibitem{behnoudfar2019higher}
P.~Behnoudfar, Q.~Deng, and V.~M. Calo.
\newblock High-order generalized-$\alpha$ method.
\newblock {\em Applications in Engineering Science}, 4:100021, 2020.

\bibitem{BEHNOUDFAR2021113725}
P.~Behnoudfar, Q.~Deng, and V.~M. Calo.
\newblock Higher-order generalized-$\alpha$ methods for hyperbolic problems.
\newblock {\em Computer Methods in Applied Mechanics and Engineering},
  378:113725, 2021.

\bibitem{LOS2020109}
M.~Łoś, P.~Behnoudfar, M.~Paszyński, and V.M. Calo.
\newblock Fast isogeometric solvers for hyperbolic wave propagation problems.
\newblock {\em Computers and Mathematics with Applications}, 80(1):109--120,
  2020.

\bibitem{Miranda1989}
I.~Miranda, R.~M. Ferencz, and T.~J.~R. Hughes.
\newblock An improved implicit-explicit time integration method for structural
  dynamics.
\newblock {\em Earthquake Engineering \& Structural Dynamics}, 18(5):643--653,
  1989.

\bibitem{Hughes1978}
T.~J.~R. Hughes and W.~K. Liu.
\newblock {Implicit-Explicit Finite Elements in Transient Analysis: Stability
  Theory}.
\newblock {\em Journal of Applied Mechanics}, 45(2):371--374, 06 1978.

\bibitem{Hughes1978b}
T.~J.~R. Hughes and W.~K. Liu.
\newblock {Implicit-Explicit Finite Elements in Transient Analysis:
  Implementation and Numerical Examples}.
\newblock {\em Journal of Applied Mechanics}, 45(2):375--378, 06 1978.

\bibitem{HULBERT1996175}
G.~M. Hulbert and J.~Chung.
\newblock Explicit time integration algorithms for structural dynamics with
  optimal numerical dissipation.
\newblock {\em Computer Methods in Applied Mechanics and Engineering},
  137(2):175 -- 188, 1996.

\bibitem{Daniel2003}
W.~J.~T. Daniel.
\newblock Explicit/implicit partitioning and a new explicit form of the
  generalized $\alpha$ method.
\newblock {\em Communications in Numerical Methods in Engineering},
  19(11):909--920, 2003.

\bibitem{BONELLI2002695}
A.~Bonelli, O.S. Bursi, and M.~Mancuso.
\newblock Explicit predictor–multicorrector time discontinuous galerkin
  methods for non-linear dynamics.
\newblock {\em Journal of Sound and Vibration}, 256(4):695 -- 724, 2002.

\bibitem{Bonelli2005}
A.~Bonelli and O.S. Bursi.
\newblock Predictor‐corrector procedures for pseudo‐dynamic tests.
\newblock {\em Engineering Computations}, 22(7):783--834, 2005.

\bibitem{Tripodi2016}
E.~Tripodi, A.~Musolino, R.~Rizzo, and M.~Raugi.
\newblock A new predictor–corrector approach for the numerical integration of
  coupled electromechanical equations.
\newblock {\em International Journal for Numerical Methods in Engineering},
  105(4):261--285, 2016.

\bibitem{lopez2020}
S.~Lopez.
\newblock A predictor–corrector time integration algorithm for dynamic
  analysis of nonlinear systems.
\newblock {\em Nonlinear Dynamics}, 101:1365–1381, 2020.

\bibitem{Thomas1973}
G.~R. Thomas.
\newblock A variable step incremental procedure.
\newblock {\em International Journal for Numerical Methods in Engineering},
  7(4):563--566, 1973.

\bibitem{Schmidt1978}
W.~F. Schmidt.
\newblock Adaptive step size selection for use with the continuation method.
\newblock {\em International Journal for Numerical Methods in Engineering},
  12(4):677--694, 1978.

\bibitem{PADOVAN1982365}
J.~Padovan and S.~Tovichakchaikul.
\newblock Self-adaptive predictor-corrector algorithms for static nonlinear
  structural analysis.
\newblock {\em Computers and Structures}, 15(4):365--377, 1982.

\bibitem{Belytschko2014}
T.~Belytschko, W.K. Liu, and B.~Moran.
\newblock {\em Nonlinear Finite Elements for Continua and Structures}.
\newblock Wiley, 2014.

\bibitem{Gustafsson1991}
K.~Gustafsson.
\newblock Control theoretic techniques for stepsize selection in explicit
  runge-kutta methods.
\newblock {\em ACM Trans. Math. Softw.}, 17(4):533–554, December 1991.

\bibitem{SODERLIND2006225}
G.~Söderlind and L.~Wang.
\newblock Adaptive time-stepping and computational stability.
\newblock {\em Journal of Computational and Applied Mathematics},
  185(2):225--243, 2006.
\newblock Special Issue: International Workshop on the Technological Aspects of
  Mathematics.

\bibitem{Butcher2016}
J.~C. Butcher.
\newblock {\em Runge–Kutta Methods}, chapter~3, pages 143--331.
\newblock John Wiley and Sons, Ltd, 2016.

\bibitem{PARK1980241}
K.C. Park and P.G. Underwood.
\newblock A variable-step central difference method for structural dynamics
  analysis — part 1. theoretical aspects.
\newblock {\em Computer Methods in Applied Mechanics and Engineering},
  22(2):241--258, 1980.

\bibitem{UNDERWOOD1980259}
P.G. Underwood and K.C. Park.
\newblock A variable-step central difference method for structural dynamics
  analysis- part 2. implementation and performance evaluation.
\newblock {\em Computer Methods in Applied Mechanics and Engineering},
  23(3):259--279, 1980.

\bibitem{Zienkiewicz1991}
O.~C. Zienkiewicz and Y.~M. Xie.
\newblock A simple error estimator and adaptive time stepping procedure for
  dynamic analysis.
\newblock {\em Earthquake Engineering \& Structural Dynamics}, 20(9):871--887,
  1991.

\bibitem{Zeng1992}
L.~F. Zeng, N.-E. Wiberg, X.~D. Li, and Y.~M. Xie.
\newblock A posteriori local error estimation and adaptive time-stepping for
  newmark integration in dynamic analysis.
\newblock {\em Earthquake Engineering \& Structural Dynamics}, 21(7):555--571,
  1992.

\bibitem{Wiberg1993}
N.-E. Wiberg and X.~D. Li.
\newblock A post-processing technique and an a posteriori error estimate for
  the newmark method in dynamic analysis.
\newblock {\em Earthquake Engineering \& Structural Dynamics}, 22(6):465--489,
  1993.

\bibitem{Romero2006}
I.~Romero and Luis~M. Lacoma.
\newblock A methodology for the formulation of error estimators for time
  integration in linear solid and structural dynamics.
\newblock {\em International Journal for Numerical Methods in Engineering},
  66(4):635--660, 2006.

\bibitem{HULBERT1995155}
G.~M. Hulbert and I.~Jang.
\newblock Automatic time step control algorithms for structural dynamics.
\newblock {\em Computer Methods in Applied Mechanics and Engineering},
  126(1):155--178, 1995.

\bibitem{Courant56}
R.~Courant, K.~Friedrichs, and H.~Lewy.
\newblock On the partial difference equations of mathematical physics.
\newblock {\em AEC Research and Development Report, NYO-7689, New York: AEC
  Computing and Applied Mathematics Centre}, 76, 1956.

\bibitem{behnoudfar2018variationally}
P.~Behnoudfar, V.~M Calo, Q.~Deng, and P.~D. Minev.
\newblock A variationally separable splitting for the generalized-$\alpha $
  method for parabolic equations.
\newblock {\em arXiv preprint arXiv:1811.09351}, 2018.

\bibitem{BEHNOUDFAR2021113656}
P.~Behnoudfar, Q.~Deng, and V.~M. Calo.
\newblock Split generalized-$\alpha$ method: A linear-cost solver for
  multi-dimensional second-order hyperbolic systems.
\newblock {\em Computer Methods in Applied Mechanics and Engineering},
  376:113656, 2021.

\bibitem{hairer2010solving}
E.~Hairer and G.~Wanner.
\newblock {\em Solving Ordinary Differential Equations II: Stiff and
  Differential-Algebraic Problems}, volume~14.
\newblock Springer, 2010.

\bibitem{SimoandHughes1998}
J.C. Simo and T.J.R. Hughes.
\newblock {\em Computational Inelasticity}, volume~7.
\newblock Springer-Verlag New York, 1998.

\bibitem{HUGHES1978391}
T.~J.R. Hughes and K.~S. Pister.
\newblock Consistent linearization in mechanics of solids and structures.
\newblock {\em Computers and Structures}, 8(3):391 -- 397, 1978.

\bibitem{BatheRamm1975}
K.-J. Bathe, E.~Ramm, and E.~L. Wilson.
\newblock Finite element formulations for large deformation dynamic analysis.
\newblock {\em International Journal for Numerical Methods in Engineering},
  9(2):353--386, 1975.

\bibitem{alnaes2015fenics}
M.~S. Aln{\ae}s, J.~Blechta, J.~Hake, A.~Johansson, B.~Kehlet, A.~Logg,
  C.~Richardson, J.~Ring, M.~E. Rognes, and G.~N. Wells.
\newblock {The FEniCS project version 1.5}.
\newblock {\em Archive of Numerical Software}, 3(100):9--23, 2015.

\bibitem{Kirby2004}
R.~C. Kirby.
\newblock Algorithm 839: Fiat, a new paradigm for computing finite element
  basis functions.
\newblock {\em ACM Trans. Math. Softw.}, 30(4):502–516, December 2004.

\bibitem{GOMEZ20084333}
H.~Gómez, V.~M. Calo, Y.~Bazilevs, and T.~J.R. Hughes.
\newblock Isogeometric analysis of the {C}ahn–{H}illiard phase-field model.
\newblock {\em Computer Methods in Applied Mechanics and Engineering},
  197(49):4333 -- 4352, 2008.

\bibitem{LANG1995223}
J.~Lang.
\newblock Two-dimensional fully adaptive solutions of reaction-diffusion
  equations.
\newblock {\em Applied Numerical Mathematics}, 18(1):223 -- 240, 1995.

\bibitem{KUHL1999343}
D.~Kuhl and E.~Ramm.
\newblock Generalized energy–momentum method for non-linear adaptive shell
  dynamics.
\newblock {\em Computer Methods in Applied Mechanics and Engineering},
  178(3):343--366, 1999.

\bibitem{Espath2015}
L.~F.~R. Espath, A.~L. Braun, A.~M. Awruch, and L.~D. Dalcin.
\newblock A nurbs-based finite element model applied to geometrically nonlinear
  elastodynamics using a corotational approach.
\newblock {\em International Journal for Numerical Methods in Engineering},
  102(13):1839--1868, 2015.

\end{thebibliography}

\end{document}